\documentclass[a4paper,twoside,12pt]{report}
\usepackage{graphicx}
\usepackage{makeidx}
\def\CaixaPreta{\vrule Depth0pt height5pt width5pt}
\def\inicproof{\noindent {\bf Proof:} \hspace{3mm}}
\def\fimproof{\phantom{.}\hfill \CaixaPreta \vspace{5mm}}
\newcommand{\n}{\noindent}
\newcommand{\noi}{\noindent}

\newcommand{\dy}{\displaystyle}

\newcommand{\aw}{\langle}
\newcommand{\fw}{\rangle}
\newcommand{\mod}{\mbox{mod}}

\newcommand{\dw}{{\downarrow}\kern -0.21cm\raise 0.10cm
\hbox{$\downarrow$}}
\newcommand{\du}{{\uparrow}\kern -0.21cm\raise 0.10cm
\hbox{$\uparrow$}}
\newcommand\barnot{{\mid}\kern -0.18cm\raise -0.015cm
\hbox{$\not$}\,\,\,\, }
\def\free#1{\smash{\mathop{\hbox to
2cm{\rightarrowfill}}\limits_{\displaystyle #1}^{}}}

\def\mapdown#1{{#1}\kern -0.05cm\raise 0.02cm
\hbox{$\Big\downarrow$}}
\def\mapup#1{{#1}\kern -0.05cm\raise 0.02cm
\hbox{$\Big\uparrow$}}

\def\mapright#1#2{\smash{\mathop{\hbox to 0.80cm{\rightarrowfill}}\limits^{#1}_{#2}}}
\def\mapleft#1#2{\smash{\mathop{\hbox to 0.80cm{\leftarrowfill}}\limits^{#1}_{#2}}}
\def\mapleftright#1#2{\smash{\mathop{\hbox to 0.80cm{\leftarrowfill \rightarrowfill}}\limits^{#1}_{#2}}}

\newfont{\fontsete}{cmr7 scaled\magstep0}
\newfont{\fontoito}{cmr8 scaled\magstep0}
\newfont{\fontnove}{cmr9 scaled\magstep0}
\newfont{\fontdez}{cmr10 scaled\magstep0}
\newcommand{\jtt}{\vspace{-1pc}}
\def\runninghead#1#2{\pagestyle{myheadings}
\markboth{{\protect\footnotesize\it{\hfill #1}}\hfill}
{\hfill{\protect\footnotesize\it{#2\hfill}}}}
\makeindex
\title{GRAPH OF MAPS}
\author{by \\ S\'ostenes Lins\\ \\ \\ \\ \\ \\ \\ \\
{\footnotesize A thesis presented to the University of Waterloo}\\
 {\footnotesize in partial fulfillment of the }\\
 {\footnotesize requirements for the degree of}\\
 {\footnotesize Doctor of philosophy in Mathematics}\\
 \\ \\ \\ \\ \\ \\ \\
 {\footnotesize Waterloo, Ontario, Canada, 1980}\\
 {\footnotesize \copyright{S\'ostenes Lins}}}

\date{}

\runninghead{S\'ostenes Lins}{Graph of Maps}

\begin{document}
\bibliographystyle{alpha}
\maketitle
\thispagestyle{empty}
\baselineskip18pt
\abovedisplayskip=0.3cm
\belowdisplayskip=0.6cm

\vfill\eject

I hereby declare that I am the sole author of this thesis.

I authorize the University of Waterloo to lend this thesis to
other institutions or individuals for the purpose of scholarly
research.

Signature: \vspace{3cm}

I further authorize the University of Waterloo to reproduce this
thesis by photocopying or by other means, in total or in part, at
the request of other institutions or individuals for the purpose
of scholarly research.

Signature: \vspace{3cm}

The University of Waterloo requires the signatures of all persons
using or photocopying this thesis. Please sign below, and give
address and date.

\vfill\eject \vglue19cm

\phantom{.}\hfill A meus Pais

\vfill\eject

\begin{abstract}
This work studies certain aspects of graphs embedded on surfaces.
Initially, a colored graph model for a map of a graph on a surface
is developed.  Then, a concept analogous to (and extending) planar
graph is introduced in the same spirit as planar abstract duality,
and is characterized topologically.  An extension of the Gauss
code problem treating together the cases in which the surface
involved is the plane or the real projective plane is established.
The problem of finding a minimum transversal of
orientation-reversing circuits in graphs on arbitrary surfaces is
proved to be NP-complete and is algorithmically solved for the
special case where the surface is the real projective plane.
\end{abstract}

\vfill\eject

\begin{center}
{\large\bf ACKNOWLEDGEMENTS}
\end{center}
\bigskip\bigskip

I am indebted to Bernardete, my wife, whose moral support was
invaluable in all phases of this work; to Lauro, my son, whose
birth during this period filled my life with joy and thus greatly
contributed to the relief of the unavoidable tensions.

I am grateful to my supervisor, Professor Daniel H. Younger. He
sparked my curiosity towards the subject of graphs on surfaces.
Our discussions, his suggestions and his always firm, patient, and
friendly criticism, greatly influenced not only the thesis, but my
attitude towards mathematics in general. I learned from every
aspect of our relationship.

I want to thank my colleagues and friends Arnaldo Mandel and
Ephraim Korach for innumerable stimulating mathematical
discussions. I also thank Professor Herb Shank: some of the key
ideas in this work arose from his theory on left-right paths and
from our personal interaction. To Professor U.S.R. Murty I express
my appreciation for sharing with me the enthusiasm for the graph
colored model.  His feedback was important to me.

I thank UFPE and to CAPES, in Brazil, for financial support.  The
CIDA/COMBRA project also provided some support. Finally, I wish to
thank Mrs. Sue Embro for her excellent typing of the thesis, as
well as for her patience with the corrections.

\vfill\eject \setcounter{chapter}{-1}

\chapter{INTRODUCTION}
This thesis is about graphs embedded in compact surfaces
satisfying the property that the complement of the graph in the
surface\index{surface} is (topologically) a collection of disjoint
open discs, called faces. We summarize the latter restriction by
saying that the graph is well-embedded in the surface. A pair
$(G,S)$ where graph $G$ is well-embedded in surface $S$ is called
a topological map\index{map} of a graph or, for short, a
$t$-map\index{$t$-map}.

In Tutte [25] an algebraic abstract model for the above object is
introduced and called simply ``map''.  In Chapter 1 of our work we
define a map as a finite, cubic, 3-edge-colored graph in which two
of the three colors induce squares.

The main point of the first part of Chapter 1 is to show that the
three objects are equivalent. Indeed, we present a natural
bijection between our maps and the algebraic maps of Tutte, which
we denote for short as $a$-maps.  We present also a natural
bijection between our maps and the $t$-maps. With these
bijections, an intuitive geometrical interpretation of Tutte's
axioms is established.

The correspondence between the maps and the a-maps is based on the
construction of the Cayley color graph of the action of three
fixed point free involutions, which are slightly different from
the three permutations which constitute Tutte's $a$-maps.

The correspondence between the maps and $t$-maps is based on the
dual of the barycentric division of the $t$-map, a very simple
topological operation. Specifically it consists of the thickening
of each edge, which is replaced by a bounding digon on the
surface, followed by the expansion of each vertex, yielding a
cubic graph\index{cubic graph} in which the edges of the original
are replaced by squares.

A suitable 3-coloration of the edges of the expanded graph can be
used to reconstruct the $t$-map; this is the basic idea: to encode
the surface and the cyclic sequences of edges around vertices
using three colors.

Fundamental concepts of $t$-maps can be directly interpreted in
the 3-edge-colored cubic graph. For instance, the Euler
characteristic is defined as a signed sum of the numbers of the
three types of 2-colored polygons. The orientability of the
surface\index{orientability of the surface}, as we show in Theorem
1.6, is equivalent to the bipartiteness of the cubic graph.
Therefore, since it is simpler, we take the latter as our
definition of orientability. We use these definitions to give a
new proof of the known fact that the Euler characteristic of an
orientable surface is even (Theorem 1.7). The interpretation of
the topological concept of an orientation-reversing circuit (a
circuit in the graph\index{circuit in the graph} of a $t$-map
which has, in the surface, a neighborhood homeomorphic to the
Moebius atrip) is the simpler concept of a circuit with an odd
number of edges. This is important in the proof of Lemma 3.4, a
central fact in Chapter 3.

Theorem 1.8 is a combinatorial counterpart, in the 3-colored
model, of the topological fact that the union of a set of
vertex-disjoint orientation-reversing circuits in a $t$-map is not
the boundary of a subset of faces\index{subset of faces}. This
fact is used to establish the result of Chapter 4. It is also used
at the end of Chapter 2 (Theorem 2.14) to derive an upper bound
for the number of independent vertices in certain graphs.

In the final part of Chapter 1 we introduce the structure of maps
which forms the basis of Chapter 2. Some topological consequences
of this structure are also considered.

The reentrant paths in a planar $t$-map which alternate using the
rightmost and leftmost edges at each vertex have been studied by
Shank [24], who names them left-right paths. They are also known
as Petrie polygons; see Coxeter [3].  These objects have natural
generalization for arbitrary surfaces, and we name them, in this
context, zigzag paths.

In our model the abstract equivalence between the following three
cyclic sequences arising from $t$-maps becomes apparent:\jtt
\begin{itemize}
\item[(a)] the cyclic sequence of edges around a vertex;
\item[(b)] the cyclic sequence of the edges in a boundary of a
face\index{boundary of a face}; \item[(c)] the cyclic sequence of
the edges in a zigzag path.
\end{itemize}

The symmetry between (a) and (b) is widely known.  It is given by
the dual $t$-map which interchanges vertices and faces.  We
observe that the dual maintains the zigzags.  We introduce the
``phial\index{phial}'' map, whose corresponding $t$-map
interchanges (relative to the original $t$-map) the vertices and
the zigzags, maintaining the faces.  Also, in the same way, we
introduce the ``antimap\index{antimap}'', which interchanges the
faces and the zigzags, maintaining the vertices.

The operations of taking the dual and taking the phial generate,
from a given map, a set of six maps in correspondence with the
elements of the symmetric group on three elements.  This structure
is used (Theorem 1.9) to prove that every graph in which the edges
are doubled can be embedded so as to form a self-dual $t$-map.

Chapter 2 introduces and investigates a concept analogous to
planar duality.

A theorem by Whitney [27] states that a graph is planar if and
only if it has an abstract dual.  This dual can be taken to mean
another graph with the same edges as the first graph such that a
cycle in one of them corresponds to a coboundary in the other.

We introduce the following concept: a graph is called (cycle) rich
if there are two other graphs with the same edges as the first
graph and such that a cycle in one of the three graphs corresponds
to the symmetric difference of two coboundaries, one for each of
the other two graphs. Thus, we have an abstract correspondence
between cycles and coboundaries of three graphs instead of two.

The main result of the first part of the chapter shows (Theorem
2.3) that triples of graphs sharing the above relation are induced
naturally from certain classes of $t$-maps. If a rich graph arises
in this way, then we say that the graph is map-rich. We do not
know whether there are rich graphs which are not map-rich. If not,
this would complete the analogy with Whitney's Theorem.

It follows, by a corollary of Proposition 2.6, that every planar
graph is rich (in fact map-rich). Hence, the concept of rich graph
is not only analogous to the concept of planar graph, but is in
fact a generalization.

The central fact, which makes the theory symmetric in the three
graphs, Theorem 2.5. It states that for arbitrary $t$-maps, the
intersection of the vector space (over $GF(2)$) of the
coboundaries of its graph and the vector space generated by the
boundaries of its faces is contained in the vector space generated
by the cycles induced by its zigzag paths. In each one of the six
maps mentioned before, the above spaces occur differently
permuted. The six maps induce the three graphs, each one occurring
twice. The symmetry follows from properties.

The main results of the second part of the chapter are of map-rich
graphs (Theorem 2.9). The first of these states that a graph is
map-rich if it has a special kind of embedding (where the vertices
are represented by closed curves and the edges by intersection of
the curves) in a surface whose Euler characteristic and
orientability are specified functions of the graph. The second
characterization states that a graph is map-rich if and only if
the special kind of embedding induces in some surface, a (usual)
$t$-map which is 2-face-colorable and such that the cycle space of
its (4-regular) graph is generated by the boundaries of its faces
and the curves which represent the vertices of the original graph.

In the third part of the chapter we use the second
characterization mentioned above to prove that every complete
graph is map-rich (Theorems 2.11 and 2.12). Moreover, from this
characterization is evident an easy topological ``recipe'' to
produce map-rich graphs from $t$-maps with arbitrary surface, in
the same sense that Whitney's theorem permits the presentation of
examples of abstract duality simply by the drawing of planar
graphs. The first characterization of map-rich graphs in Theorem
2.9 is used (Lemma 2.13) to prove that some graphs are not
map-rich, for instance, the graph obtained from the Petersen graph
by doubling all its edges.

Trying to get more insight in the membership decision problem for
map-rich graphs we were led to study the Gauss code problem, which
is the topic of our next chapter.

The Gauss code problem is the following: given a cyclic sequence
of symbols in which each symbol occurs twice, find necessary and
sufficient conditions to embed a closed curve in the plane with
the restriction that the cyclic sequence of self-intersections of
the curve reproduces the given cyclic sequence, when the symbols
are identified with the self-intersection points.

A brief survey of this problem can be found in Gr\"unbaum
[10-pp.72-73). Lovasz and Marx [17] present a solution in terms of
a forbidden substructure. Rosenstiehl [19] presents a solution in
terms of a graph called an interlace graph\index{interlace graph},
which is easily computable from the given cyclic sequence.  This
solution includes a polynomial algorithm to test and realize the
embedding (if it exists).

A generalized version of this problem can be considered: dropping
the planar restriction, we ask for a surface of minimum
connectivity for which there is an embedding of a curve satisfying
the prescribed sequence of self-intersections.

The embedding of the curve induces a $t$-map with a 4-regular
graph, since each crossing point is visited twice. Consider the
special case of the generalized problem, where the restriction of
2-face-colorability of the induced $t$-map is imposed.

Note that this restricted form is yet a generalization of the
original problem, since every 4-regular graph embedded in the
plane is 2-face-colorable. The motivation for considering this
special case is that an algorithmic solution to it would enable us
to determine whether an arbitrary cubic graph is map-rich.

The main result of Chapter 3 is an algorithmic answer for the
question of embeddability satisfying the 2-colorability
restriction in the case of surfaces of connectivity at most one.
Our approach treats together the cases of the plane and projective
plane. Particularized to the plane, our conditions reproduce those
of Rosenstiehl [19].

This result uses some algebraic properties of $t$-maps with one
vertex, which we establish. These properties are based on the
generalization for arbitrary surfaces, of some homomorphisms,
which we call crossing functions, introduced for the plane by
Shank [24] and also used by Rosenstiehl [19].

The $t$-maps induced by the-restricted form of the generalized
Gauss code problem are in 1-1 correspondence with the $t$-maps
with one vertex. A $t$-map with one vertex induces naturally a
crossing function. The central fact, from which we specialize the
conditions for the plane and the projective plane, is the equality
(Theorem 3.3(a)) between the dimension of the image of the
composition of the crossing functions induced by a map and its
antimap (whose associated $t$-maps have one vertex), and the
connectivity of the surface associated with the phial map.

Chapter 4 presents an algorithmic solution for the problem of
finding a minimum transversal of orientation reversing circuits in
graphs embedded in the real projective plane.

Our interest in this problem arises from the fact that for graphs
embedded in arbitrary surfaces, it includes the problem of finding
a maxinum coboundary in an arbitrary graph (Proposition 4.2).
Since the latter problem is known to be $NP$-complete, (Garey and
Johnson [6]) so is the problem that we treat, if considered for
graphs in arbitrary surfaces.

Certain subsets of faces of a $t$-map which has an eulerian graph,
called reducers, are recognized. By means of a suitable operation
which eliminates a reducer\index{reducer}, a $t$-map is
transformed to a ``smaller'' one, without change in the
cardinality of a minimum transversal of orientation reversing
circuits. For the case of graphs embedded in the projective plane,
the non-existence of a reducer implies that the $t$-map has a
simple form, looking like a system of lines in the real projective
plane. From this fact we deduce the following minimax relation
(Theorem 4.5): for an eulerian graph $G$ well-embedded in the real
projective plane, the cardinality of a minimum transversal of
orientation-reversing circuits in $G$ is equal to the cardinality
of a maximum disjoint collection of orientation-reversing circuits
in $G$.

An important fact, which permits the use of the reducers, is that
for projective $t$-maps a minimal transversal of
orientation-reversing circuits is an orientation-reversing circuit
in the dual $t$-map. We generalize this relation by showing in
Theorem 4.3 that the minimal transversals of orientation-reversing
circuits are cycles of a specific homology (mod 2) in the dual
$t$-map. This theorem implies that the elimination of a reducer
does not alter the cardinality of a transversal of
orientation-reversing circuits. The reducers are defined for
orientable surfaces as well. We show with an example (fig. 4.6)
that $t$-maps free of reducers in surfaces of higher connectivity
(even the torus) can be very complicated.
\chapter{Cubic Graphs, Three Colors, and Maps}
\section*{1.A\quad Maps and $a$-maps}
The graphs which we use may contain loops and multiple edges. The
set of edges of graph $G$ is denoted by $eG$; its set of vertices
is denoted by $VG$.  The two vertices which are incident to an
edge are called the {\it ends} of the edge.  If the two ends
coincide the edge is called a {\it loop}\index{\it loop}.

The {\it valency} (or {\it degree}) of a vertex\index{{\it degree}
of a vertex} $v$ of graph $G$, $val_G(v)$, is the number of edges
incident to $v$ where loops are counted twice. A subset of edges
in $e$G is called a {\it perfect matching}\index{\it perfect
matching} if it consists of a set of non-incident edges which
covers every vertex of $G$. A graph $G$ is called {\it cubic}, as
usual, if $val_G(v)=3$ for every vertex $v$ in $VG$. A graph $G$
is a {\it polygon}\index{\it polygon} if it has at least one edge,
is connected, and $val_G(v)=2$ for all $v$ in $VG$.

A map $M$ is an ordered triple $(C_M,v_M,f_M)$ where:\jtt
\begin{itemize}
\item[(i)] $C_m$ is a finite cubic graph; \item[(ii)] $v_M$ and
$f_M$ are disjoint perfect matchings in $C_M$, such that the
subgraph of $C_M$ induced by $v_M\cup f_M$ is a collection of
polygons with four edges or {\it squares}.
\end{itemize}

From the above definition, it follows that $C_M$ may contain
double edges but not loops.

A third perfect matching in $C_M$ is $eC_M-(v_M\cup f_M)$ and is
denoted by $a_M$. The set of diagonals of the squares, denoted by
$z_M$ is a perfect matching in the complement of $C_M$.

The edges in $v_M,f_M,a_M$ are called respectively $v_M$-edges,
$f_M$-edges, $a_M$-edges.

In this section we show a natural 1-1 correspondence between the
above objects and the algebraic maps introduced by Tutte in [25].
As we show, they are essentially the same; however, the definition
of orientability in our model is simpler. In the next section we
relate our maps with the usual topological model.

Tutte's paper is a permutation representation for graphs embedded
in arbitrary surfaces. Here, we call it an algebraic map. In this
section we do not give the topological interpretation of the
objects.

An {\it algebraic map} $A$, or {\it $a$-map}\index{\it $a$-map}
$A$, is an ordered triple of permutations,
$$A=(R,\Theta,\Phi),$$
acting on a finite set $B$ satisfying the following
axioms:\bigskip
$$\begin{array}{ll}
\mbox{(am1)}\qquad &
\Theta\Phi=\Phi\Theta,\quad\Theta^2=\Phi^2=\mbox{identity};\\
[0.5cm] \mbox{(am2)} & x,\Theta x,\ \Phi x,\ \Theta\Phi
x\quad\mbox{are different for all}\quad x\in B;\\ [0.5cm]
\mbox{(am3)} & R\Theta=\Theta R^{-1};\\ [0.5cm] \mbox{(am4)} &
R^n(x)\ne\Theta x,\quad\mbox{for all integer}\quad n\quad\mbox{and
all}\quad x\quad\mbox{in}\quad B.
\end{array}$$
If, in addition, the following property is satisfied, $A$ is a
{\it connected} $a$-map\index{{\it connected} $a$-map}:
$$\mbox{(amS)}\qquad \aw R,\Theta,\Phi\fw\quad\mbox{is transitive in}\quad B,$$
where $\aw R,\Theta,\Phi\fw$ represents the group of permutations
generated by $R,\Theta$, and $\Phi$.

In [25] axiom (am5) is imposed. If it does not hold, the object
there is called a pre-map. Paper [25] and an observation due to
V.A. Liskovets [30, pg.14], which consists essentially in the
replacement of $R$ by $R\Theta$ in the above axioms, were the
basis for the definition of map that we give.

We define a function $\psi^a$, from the set of maps into the set
of $a$-maps as follows: given map $M=(C_M,V_M,f_M)$, the $a$-map
$\psi^a(M)$, which we simplify to $M^a$, is $M^a=(a_Mv_M,v_M,f_M)$
acting on the vertices of $C_M$. For this purpose, we identify the
perfect matchings in $C_M,V_M,f_M$, and $a_M$ as fixed-point-free
involutions on $VC_M$ in the natural way: if $b$ is a perfect
matching in $(a_M,v_M,f_M)$ and $x\in VC_M$, then $b(x)$ is the
other vertex in $VC_M$ such that $x$ and $b(x)$ are the ends of an
edge\index{ends of an edge} in $b$.

\medskip\n{\bf (1.1) Proposition:} $\psi^a$ is a bijection from the set of maps
onto the set of a-maps.
\medskip

\inicproof Initially we show that $M^a$ is indeed an $a$-map, by
verifying Tutte's axioms for it.  The squares formed by the
perfect matchings $v_M$ and $f_M$, and the interpretation of these
as fixed-point-free involutions account for (am1) and (am2). Axiom
(am3) also follows,
$$\begin{array}{rcl}
(a_M\circ v_M)\circ v_M &=& a_M\\
&=& v_M\circ(v_M\circ a_M)\\
&=&v_M\circ(a_M\circ v_M)^{-1}
\end{array}$$

We prove (am4) for $M^a$ indirectly, assuming that it fails. Since
$VC_M$ is finite we can restrict ourselves to non-negative
exponents in (am4). Choose the pair $(n,x)$ such that $n$ is the
smallest possible non-negative integer with (am4) failing. First
we show that $n>1$. If $n=0$, then $x=v_M(x)$, contradicting the
fact that $v_M$ is a matching in $C_M$. If $n=1$, then $a_Mv_M(x)=
v_M(x)$, or $a_M(y)=y$, with $y=v_M(x)$; this contradicts the fact
that $a_M$ is a matching in $C_M$. Thus $n>1$. We then have
$$(a_M,v_M)^n(x)=v_M(x);$$
equivalently
$$(a_Mv_M)^{n-1}(a_Mv_M)x=v_M(x),$$
or $$v_Ma_M(a_Mv_M)^{n-1}(a_Mv_M(x))=v_Ma_Mv_M(x).$$ Hence,
$(a_Mv_M)^{n-2}(y)=v_M(y)$ where $y=a_Mv_M(y)$. This contradicts
the choice of $(n,x)$; thus, (am4) holds for $M^a$. The proof that
$M^a$ is an $a$-map is complete.

To prove that $\psi^a$ is a bijection, we define a function
$\overline{\psi}^a$ from the set of $a$-maps into the set of maps
such that, functions $\overline{\psi}^a$ and $\psi^a$ are
inverses. Given an $a$-map $A=(R,\Theta,\Phi)$ acting on $B$ we
define $\overline{\psi}^a(A)$ as $M=(C_M,v_M,f_M)$, where
$VC_M=B,\ v_M=\Theta,\ f_M=\Phi$, and $a_M=R\Theta$. The edges of
$C_M$ are defined by $v_M,f_M,a_M$. It follows that $C_M$ is a
finite cubic graph. Axioms (am1) and (am2) imply that the subgraph
of $C_M$ induced by $v_M\cup f_M$ is a collection of disjoint
squares. Therefore, $M=\overline{\psi}^a(A)$ is a map. Observe
that $\overline{\psi}^a$ can be seen as just the construction of a
Cayley color graph for the action of $(R\Theta,\Theta,\Phi)$ over
$B$, according to (9); we drop the directions on the edges,
replacing a pair of opposite oriented parallel edges by a single
unoriented edge. This is possible without losing information
because $R\Theta,\Theta$, and $\psi$ are involutions. It is
evident from the definitions of $\psi^a$ and $\overline{\psi}^a$
that they are inverses. This concludes the proof of the
proposition. \fimproof

Observe that (amS) holds for an $a$-map $A=(R,\Theta,\Phi)$ iff
$C_M$ is connected, where $M=\overline{\psi}^{a}(A)$. In general,
the orbits of $\aw R,\Theta,\phi\fw$ acting over $VC_M$ are the
components of $C_M$.

The {\it components} of a map\index{{\it components} of a map} $M$
are the components of $C_M$. $M$ is {\it connected} if $C_M$ has
one component.

The $v_Mf_M$-polygons in $C_M$ are called the {\it squares} of
$M$, or the {\it $M$-squares}\index{\it $M$-squares}. The set of
$M$-squares is denoted by $SQ(M)$. The $v_Ma_M$-polygons in $C_M$
are called the $v$-gons of $M$. The set of $v$-gons of $M$ is
denoted by $vG(M)$.

The $f_Ma_M$-polygons in $C$ are called the {\it
$f$-gons}\index{\it $f$-gons} of $M$.

The set of $f$-gons of $M$ is denoted by $fG(M)$.

The $z_Ma_M$-polygons in $C_M\cup z_M$ are called the {\it
$z$-gons} of $M$. The set of $z$-gons of $M$ is denoted by
$zG(M)$. The vertices of $C_M$ are also called the {\it
corners}\index{\it corners} of $M$.

The {\it Euler characteristic}\index{\it Euler characteristic} of
a map $M$ is the integer
$$\chi(M)=|vG(M)|-|SQ(M)|+|fG(M)|.$$

The {\it connectivity} of a map $M$ is the integer $\xi(M)$
defined as $2-\chi(M)$.

A map $M$ is called {\it orientable} if $C_M$ is bipartite;
otherwise it is called {\it non-orientable}. The {\it
orientability} label of $M$, $OL(M)$, is defined as $1$ if $M$ is
orientable and as $-1$ otherwise.

An $a$-map $A=(R,\Theta,\Phi)$ acting on $B$ is {\it
non-orientable} (see [25]) if $\aw R,\Theta\Phi\fw$ has one orbit
for some orbit of $\aw R,\Theta,\Phi\fw$. The alternative (as we
show in the proof of the next proposition) is for $\aw
R,\Theta\Phi\fw$ to have two orbits for every orbit of $\aw
R,\Theta,\Phi\fw$. In this case, $A$ is called orientable.

The concepts of orientability for maps\index{orientability for
maps} and $a$-maps (as well as for topological maps, as we show
later) agree.

We have the following proposition:

\medskip\n{\bf (1.2) Proposition:} Map $M=(C_M,v_M,f_M)$ is orientable iff $a$-map $M^{a}$ is orientable.
\medskip

\inicproof We can, evidently, restrict the proof to one component
of $C_M$, corresponding to one orbit of $\aw a_Mv_M,v_M,f_M\fw$,
acting on $VC_M$. In other words, we may suppose $C_M$ connected.

Rephrasing the proposition according to the definitions, we want
to prove that $C_M$ is bipartite iff $\aw a_Mv_M,v_Mf_M\fw$ has
more than one orbit on $VC_M$. Observe that the actions of
$a_Mv_M$ and $v_Mf_M$ correspond to paths of length two in $C_M$.
It follows that if $C_M$ is bipartite, then $\aw a_Mv_M,v_Mf_M\fw$
is not transitive on $VC_M$. To prove the converse, we establish
as a lemma that there are at most two orbits of $\aw
a_Mv_M,v_Mf_M\fw$ acting on $VC_M$, for arbitrary connected $M$.
To this end, we first show that all the four (or less) squares
adjacent to a square $S$ by $a_M$-edges are reachable by a path
which starts in a fixed corner $x$ in $S$, defined by a member of
$\aw a_Mv_M,v_Mf_M\fw$. Take the following four paths defined
as\jtt
\begin{itemize}
\item[(i)] $(a_Mv_M)(v_Mf_M)x$; \item[(ii)] $(v_Ma_M)(v_Mf_M)x$;
\item[(iii)] $(v_Ma_M)x$; \item[(iv)] $(a_Mv_M)x.$
\end{itemize}
\noi Note that $v_Ma_M=(a_Mv_M)^{-1}\in \aw a_Mv_M\fw$ since
$VC_M$ is finite.

Observe that all the squares adjacent to $S$ by an $a_M$-edge
contain the final vertex of at least one of the above paths. By
the connectivity of $C_M$, it follows that every orbit of $\aw
a_Mv_M,v_Mf_M\fw$ contains vertices of all squares. Since opposite
vertices of a square are in the same orbit, it follows that there
are at most two orbits of $\aw a_Mv_M,v_Mf_M\fw$ on $VC_M$.
Moreover, if there are two orbits, each contains two opposite
vertices of each square.

Now, we can easily conclude the proof. If $C_M$ is not bipartite,
consider an odd polygon in $C_M$.  At least two adjacent vertices,
$x$ and $y$, of this polygon are in the same orbit, as we proved
that there are at most two orbits. We may consider that $x$ and
$y$ are in the same square. If not, then they are linked by an
$a_M$-edge and the replacement of $y$ by $v_Ma_M(y)$ fixes the
situation. Therefore, we have two adjacent vertices of a square in
the same orbit of $\aw a_Mv_M,v_Mf_M\fw$.  Hence, by the argument
above, which proves that each orbit contains two opposite vertices
of each square, there is just one orbit. This concludes the proof.
\fimproof

With the above proposition we conclude the analysis relating
$a$-maps and maps.

\section*{1.B\quad Maps and t-maps}
In this section we establish a natural 1-1 correspondence between
our maps and the usual maps of graphs embedded on surfaces. Some
properties of the latter objects are directly established in our
model.

A {\it surface (without boundary)} is a Hausdorff topological
space in which every point has a neighborhood homeomorphic to an
open 2-ball, or open {\it disc}\index{\it disc}.

It is a classical theorem in topology that the surfaces which are
connected, compact and without boundary are determined (up to
homeomorphism) by the Euler characteristic and orientability.

The surfaces which we use have no boundary, are connected and,
with the exception of the plane, are compact. In view of this we
no longer use these qualifications. The plane, in reality, is
simply a convenient way to replace the 2-sphere. Thus, we talk
about planar maps where, to be precise, we should talk about
spherical maps. We transform the sphere into the plane using a
stereographic projection.

We denote the Euler characteristic of a surface $S$ by $\chi(S)$.
Also, let $OL(S)$ be $+1$ or $-1$ according to whether or not $S$
is orientable.

Given a connected map\index{connected map} $M$, we associate with
$M$ the surface $S_M$, characterized by $\chi(S_M)=\chi(M)$ and
$OL(S_M)=OL(M)$.

In view of this, {\it planar} maps are connected maps $M$ which
have $\chi(M)=2$ and $C_M$ bipartite. We show in Theorem (2.0.2)
that the latter condition about the bipartiteness of $C_M$ is
redundant.

{\it Projective} maps\index{{\it Projective} maps} are connected
maps $M$ for which $S_M$ is the surface known as the real
projective plane. That means $\chi(M)=1$ and $C_M$ not bipartite.
The latter condition is shown to be redundant in Theorem (1.7).

{\it Toroidal} maps\index{{\it Toroidal} maps} are connected maps
$M$ for which $\chi(M)=O$ and $C_M$ is bipartite. $S_M$ in this
case is the familiar torus.

{\it Kleinian} maps\index{{\it Kleinian} maps} are connected maps
$M$ for which $\chi(M)=O$ and $C_M$ is not bipartite. $S_M$ in
this case is the Klein bottle.

By a {\it path} in a graph\index{{\it path} in a graph} $G$ we
mean a sequence $\pi=(v_0,e_1,e_2,\dots,v_{n-1},e_n,v_n)$ such
that $v_i\in VG,\ e_i\in eG$ and $v_{i-1}$ and $v_i$ are the ends
of edge-term $e_i$ for each $i$ between 1 and $n$. A path is {\it
degenerate}\index{\it degenerate} if it consists of just one
vertex-term. Vertex $v_0$ is called the {\it origin} of $\pi$ and
vertex $v_n$ is its {\it terminus}. The integer $n$ is called the
length of the path. A path is {\it reentrant} if its origin is
equal to its terminus. A {\it cyclic} path\index{{\it cyclic}
path} is a reentrant path\index{reentrant path} in which the
sequence is considered cyclic, i.e., the origin can be adjusted to
any vertex-term. If $\pi$ is a path, the set of edges which appear
in the sequence $\pi$ as edge-terms is denoted by $e\pi$;
analogously, the corresponding set of vertices is denoted by
$V\pi$. If the terminus of $\pi_1$ is the same as the origin of
$\pi_2$, then we denote by $\pi_1\circ\pi_2$ the composition of
the paths $\pi_1$ and $\pi_2$. The inverse of path $\pi$ is
denoted $\pi^{-1}$.

We proceed to show that the correspondence $M\to S_M$ is a natural
one.

Assume that $\pi=\{\pi_1,\dots,\pi_f\}$ is a collection of cyclic
paths in $G$ that uses each one of its edges twice. Suppose also
that the identifications of the ends of the edges in $G$ forced
only by the transitions in paths of $\pi$ reproduce the vertices
of $G$. Then $\pi$ is called a {\it combinatorial
embedding}\index{\it combinatorial embedding} of $G$.

For the next definition we consider a graph as a topological
object, in the usual sense.

A graph $G$ is {\it well-embedded} in a surface $S$ if $G$ is
embedded in $S$, in the topological sense, and if the difference
$S\backslash G$ is topologically equivalent to a collection of
disjoint open discs.

A $t$-map $X$ is a pair $(G,S)$ in which graph $G$ without
isolated vertices is well-embedded in surface $S$.

Associated with a combinatorial embedding of a graph $G$, there is
a natural $t$-map: let each one of the cyclic paths which defines
the embedding bound a disc; identify pairwise the boundaries of
these (disjoint) discs by ``gluing'' the two occurrences of each
edge so as to identify correctly its ends (a loop is dealt with by
the introduction of an artificial bivalent vertex), and so as to
preserve the interior of each disc as an open disc.

The final object can be seen to form a $t$-map $(G,S)$, where $S$
is a surface which has Euler characteristic $\chi(S)=v-e+f$;
$v,e$, and $f$ are, respectively, the number of vertices and edges
of $G$, and the number of cyclic paths forming the combinatorial
embedding. $S$ is {\it orientable} if it is possible to choose a
direction of traversal of the cyclic paths such that for each edge
e with ends $v$ and $u$, the traversal is once from $v$ to $u$ and
once from $u$ to $v$ (a loop is, again, dealt with by the
introduction of an artificial bivalent vertex).

The values of $\chi(S)$ and $OL(S)$, as explained above, are the
same if we take an arbitrary graph $H$ well-embedded in $S$.
Therefore, the determination of $\chi(S)$ and $OL(S)$ is effected
via some $t$-map $(H,S)$. The proof of this invariance is obtained
by transforming an arbitrary $t$-map $(H,S)$ into a canonical one,
$(H_S,S)$, by means of elementary operations such as never to
change $\chi(S)$ or $OL(S)$. This process is usually effected in
connection with the classification of surfaces. (See Ringel [18]
or Giblin [7]).

A dart\index{dart} is an oriented edge. The following proposition
is straightforward. It characterizes orientability in two slightly
different ways which are used in the proof of Theorem 1.6.

\medskip\n{\bf (1.3) Proposition:} Suppose $\chi=(G,S)$ is a $t$-map. The following
two conditions are equivalent for the orientability of $S$:\jtt
\begin{itemize}
\item[(i)] Represent each facial path\index{facial path} of $\chi$
by a polygon; let these polygons bound disjoint discs; $S$ is
orientable iff it is possible to simultaneously embed all these
discs in the plane such that for each edge of $G$ it is possible
to ``glue'' along the two occurrences (in the boundary of the
discs) of the edge ``sliding'' or ``stretching'' the discs, but
not leaving the plane. \item[(ii)] Orient arbitrarily the edges of
$G$; as in (i), represent the facial paths of $X$ by polygons and
let these polygons bound disjoint discs; $S$ is orientable iff it
is possible to embed all the discs in the plane such that for each
dart of $G$ one occurrence in the boundary of the disc is
clockwise and the other is counterclockwise.
\end{itemize}
\medskip

For a $t$-map $X=(G,S)$ the disjoint open discs in $S\backslash G$
are called the {\it faces} of $X$. Observe that the boundary of
each such disc corresponds in $G$ to a cyclic path which we call a
{\it facial path}. Note that an edge may occur twice in a facial
path.

The {\it (combinatorial) boundary} of a face of $X$ is the set of
edges occurring once in the corresponding facial path. Observe
that this notion agrees with the topological counterpart, and so,
we drop the adjectives. A boundary in $X$ is the $\mod\,2$ sum of
the boundaries of faces. If $F$ is a set of faces of $X$, the
$\mod\,2$ sum of the boundaries of the faces in $F$ is denoted by
$\partial_x(F)$, or by $\partial(F)$ if $X$ is understood.

Consider a cubic graph $G$ with a proper 3-coloring of its edges
in colors $a,b,c$. (This has the usual meaning of a partition of
$eG$ into three classes such that for each vertex the three edges
incident to it are in different classes.)

The subgraph of $G$ formed by edges of any two of the three colors
consists of a number of disjoint polygons of even length, called
$ab$-, $bc$-, $ca$-polygons.

The {\it faithful} embedding of $G$ relative to a 3-coloring in
colors $a,b,c$ is given by the cyclic paths induced by the $ab$-,
the $bc$-, and the $ca$-polygons. (It is evident that they form a
combinatorial embedding.)

Given a map $M=(C_M,V_M,f_M)$ there is a natural 3-coloring of the
edges of $C_M:\{v_M,f_M,a_M\}$. (Recall that $a_M$ is the third
perfect matching in $C_M$.)

The {\it embedding of a map} $M$ is the faithful
embedding\index{faithful embedding} of $C_M$ relative to
$\{v_M,f_M,a_M\}$.

Now we define a function $\psi$ which is a bijection from the set
of maps onto the set of $t$-maps. We denote $\psi(M)$ by $M^t$.
Given a map $M$, to obtain $M^t$ we proceed as follows. Consider
the $t$-map $(C_M,S(M))$ given by the faithful embedding of $M$.
The $v$-gons, the $f$-gons and the $M$-squares are boundaries of
(closed, in this case) discs embedded (and forming) the surface
$S(M)$. Shrink to a point the disjoint closed discs bounded by
$v$-gons. The $M$-squares, then, become bounding digons. Shrink
each such bounding digon to a line, maintaining unaffected its
vertices. With these contractions, effected in $S(M)$, $t$-map
$(C_M,S(M))$ becomes, by definition, $M^t=(G_M,S(M))$. Graph $G_M$
is called the {\it graph induced} by $M$. A combinatorial
description of $G_M$ can be given as follows: the vertices of
$G_M$ are the $v$-gons of $M$; its edges are the squares of $M$;
the two ends of an edge of $G_M$ are the two $v$-gons (which may
be the same) that contain the $v_M$-edges of the corresponding
square.

It is evident that $\psi$ is inversible: given a $t$-map we
replace each edge by a bounding digon in its surface, and then
expand each vertex to a disc in order to obtain the embedding of a
map. Therefore, $\psi^{-1}$ is well-defined; in fact, it is the
dual of a useful construction in topology, namely, barycentric
division. Thus, $\psi$ is a bijection from the set of maps onto
the set of $t$-maps.

It can be observed that $\psi$ induces a bijection from the set of
$M$-squares onto the set of edges of $G_M$. We frequently use this
bijection to identify the sets $SQ(M)$ and $eG_M$.

A subset $W$ of edges in a graph $G$ is said to be a {\it cycle}
in $G$ if the subgraph induced by $W$ has even valency at every
vertex. If $W=\phi$, we have the null cycle. A circuit is a
minimal non-null cycle. A circuit is, then, the set of edges of a
polygon. The set of cycles of a graph $G$ forms a vector space
under $\mod\,2$ sum. This space is denoted by $CS(G)$.

Bijection $\psi$ induces for each map $M$ a correspondence,
denoted $\psi^M_c$, from the cycle apace of $C_M,CS(C_M)$, onto
the cycle space of $G_M,CS(G_M)$. The correspondence $\psi^M_c$ is
defined as follows: for $S\in CS(C_M)$, an edge $s\in eG_M$ is in
$\psi^M_c(S)$ if square $s\in SQ(M)$ meets $S$ in exactly one
$f_M$-edge. With this definition, it is evident that $\psi^M_c(S)$
is a cycle in $G_M$\index{cycle in $G_M$} and that $\psi^M_c$ is
surjective.

Every vector space that we use is over $GF(2)$. Therefore, we do
not distinguish explicitly between a subset and its characteristic
vector, nor between symmetric difference and $\mod\,2$ sum of
characteristic vectors.

The following proposition permits the combinatorial definition of
orienta\-tion-reversing circuits\index{orientation-reversing
circuits}. An orientation-reversing circuit in a $t$-map $M^t$ is
a circuit in $G_M$ which has in $S_M$ a neighborhood homeomorphic
to the Moebius strip.

\medskip\n{\bf (1.4) Proposition:} $\psi^M_c$ is a homomorphism, for every map $M$.
Its kernel is the subspace of $CS(C_M)$ generated by the edges in
the $v$-gons and in the squares of $M$.
\medskip

\inicproof Let $S_1$ and $S_2$ be cycles in $C_M$. We must show
that $\psi^M_c(S_1+S_2)=\psi^M_c(S_1)+\psi^M_c(s_2)$. An edge $e$
of $G_M$ is in $\psi^M_c(S_1+S_2)$ iff exactly one $f_M$-edge of
square $e\in SQ(M)$ is in $S_1+S_2$. Therefore,
$e\in\psi^M_c(S_1+S_2)$ iff an even number of $f_M$-edges of
square $e$ belongs to one of $S_1$ and $S_2$, and one $f_M$-edge
of square $e$ belongs to the other. The latter statement is
equivalent to $e\in\psi^M_c(S_1)+\psi^M_c(S_2)$. This proves that
$\psi^M_c$ is homomorphism.

The image under $\psi^M_c$ of the edges in an arbitrary square or
$v$-gon of $M$ is the null cycle in $G_M$. Thus, the space
generated by the squares and $v$-gons is contained in ${\rm
Ker}(\psi^M_c)$. Conversely, suppose that $S\in{\rm
Ker}(\psi^M_c)$. The intersection of $S$ with the edges of an
arbitrary square has zero or two $f_M$-edges. Denote by $T$ the
cycle formed by the edges of the squares $Q\in SQ(M)$ such that
$eQ\cap S$ contains the two $f_M$-edges of $Q$. It follows that
$S+T$ has no $f_M$-edges. Since $S+T$ is the edge set of a
collection of polygons, ($C_M$ is cubic), it follows that $S+T$ is
formed by a collection of $v$-gons, whose edge-set is denoted by
$U$. Hence, $S=T+U$, with $T$ induced by squares, $U$ induced by
$v$-gons. Therefore, we conclude that ${\rm Ker}(\psi^M_c)$ is
contained in the space generated by the edge sets of squares and
$v$-gons of $M$. The proof of the proposition is complete.
\fimproof

Since an element of the kernel of $\psi^M_c$ has an even number of
edges of $C_M$, it follows that if $\psi^M_c(S_1)=\psi^M_c(S_2)$,
then $|S_1|\equiv|S_2|$ $\mod\,2$. This observation makes the
following definition meaningful.

A cycle $S$ in $G_M$ is called an {\it $r$-cycle}\index{\it
$r$-cycle} in $M^t$ if $\psi^M_c(S')=S$ and $|S'|$ is odd, for
some cycle $S'$ in $C_M$. If $|S'|$ is even and $\psi^M_c(S')=S$,
then we say that $S$ is an {\it $s$-cycle}\index{\it $s$-cycle} in
$M^t$.

We observe that the $r$-circuits in $M^t$ are precisely the
orientation-reversing circuits in $M^t$. This topological notion
is not used; we work with our parity definition of $r$-cycle.

An {\it $r$-loop}\index{\it $r$-loop} in $M^t$ is a loop in $G_M$
which is the (only) edge of an $r$-circuit in $M^t$. The
definition of {\it $s$-loop}\index{\it $s$-loop} is similar,
replacing $r$-circuít by $s$-circuit.

Observe that a subset $T\subseteq eG_M$ is a boundary in $M^t$ iff
there exists a cycle $T'$ of $C_M$, such that $\psi^M_c(T')=T$,
and $T'$ can be written as the $\mod\,2$ sum of nome subsets of
$v$-gons, $f$-gons, and squares.

For a map $M$, two cycles in $G_M$ are {\it homologous}\index{\it
homologous} $(\mod\;2)$ if their symmetric difference is a
boundary in $M^t$. Homology $(\mod\,2)$ is, thus, an equivalence
relation.

The following proposition shows that $S_M$ and $S(M)$ are the same
surface.

\medskip\n{\bf (1.5) Proposition:} For any map $M$ we have
$$S_M=S(M).$$

\inicproof First we prove that $\chi(S_M)=\chi(S(M))$. By
definition $\chi(S_M)=\chi(M)=|vG(M)|-|SQ(M)|+|fG(M)|$. Since
$(C_M,S(M))$ is a $t$-map, we use $C_M$ to get $\chi(S(M))$. Graph
$C_M$ has $6|SQ(M)|$ edges and $4|SQ(M)|$ vertices. The number of
cyclic paths that define the combinatorial embedding is
$|SQ(M)|+|vG(M)|+|fG(M)|$. Using the definition of $\chi(S(M))$
via $(C_M,S(M))$ and simplifying we get
$$\chi(S(M))-|VG(M)|-|SQ(M)|+|fG(M)|=\chi(S_M).$$
By Theorem 1.6, proved below, the faithful embedding of $C_M$
relative to \break $\{v_M,f_M,a_M\}$ is orientable iff $C_M$ is
bipartite. Thus, Proposition (1.5) is complete, provided Theorem
1.6 is proved.\fimproof

\n{\bf (1.6) Theorem:} Assume we are given a cubic graph $G$ and a
proper $3$-coloring of its edges. The faithful embedding of $G$
relative to the 3-coloring defines an orientable surface if and
only if $G$ is bipartite.

\vfill\eject

\inicproof Call the three colors $a,b,c$. Assume first that the
faithful embedding in orientable. Consider disjoint discs which
are bounded by the $ab$-, the $bc$-, and the $ca$-polygons. By
Proposition 1.3 (i) it is possible to embed the discs in the
plane, such that for an arbitrary edge we can ``slide'' one of the
discs whose boundary contains one Occurrence of the edge, so as to
``glue'' it with its second occurrence. Note that they have
correct orientation, without the need of ``turning over'' one of
the discs. Therefore, after identification along the edge, it has
an end in which the cyclic sequence of edges around it is of type
$a-b-c$ clockwise, and an end where it is $a-b-c$
counterclockwise. Also, if two edges are incident to the same
vertex, then this vertex is of type clockwise or counterclockwise
for both edges. Thus, we can partition $VG$ into clockwise and
counterclockwise vertices. Moreover, each edge links a clockwise
vertex to a counterclockwise one. It follows that $G$ is
bipartite.

Conversely, assume that $G$ is bipartite. Orient the edges
consistently from one class to the other, such that every vertex
of one class is a source and every vertex in the other class is a
sink. Embed the discs bounded by $ab$-polygons such that the
$a$-darts are clockwise and consequently the $b$-darts are
counterclockwise. Embed the discs bounded by $bc$-polygons such
that the $b$-darts are clockwise and the $c$-darts are
counterclockwise. Finally, embed the discs bounded by
$ca$-polygons such that the $c$-darts are clockwise and the
$a$-darts are counterclockwise. Since every dart appears once
clockwise and once counterclockwise, it follows from Proposition
1.3 (ii) that the surface of the faithful embedding is orientable.
\fimproof

The next theorem establishes, using our definitions, the known
fact that $\chi(M)$ is even for orientable maps. We need this
result in Proposition 3.2.2.

\medskip\n{\bf (1.7) Theorem:} Suppose that $G$ is a cubic 3-edge-colored
bipartite graph with $2k$ vertices. Then the parity of the number
of faces of the faithful embedding of $G$ relative to the
3-coloration is the same as the parity of $k$. Consequently, if
$M$ is an orientable map, then $\chi(M)$ is even.
\medskip

\inicproof Denote the colors by $a,b,$ and $c$. Let $(X,Y)$ be a
bipartition of $G$. Since $G$ is regular, $|X|=|Y|=k$. Consider
the functions $\pi_{ab},\pi_{bc}$, and $\pi_{ca}$ acting on $X$,
where $\pi_{ab}(u)$, for instance, is the terminus of the
path\index{terminus of the path} of length 2 which starts at $u$,
uses an $a$-edge and then a $b$-edge. It is clear that the image
of $u$ lies in $X$, by the bipartiteness of $G$. We prove that
$\pi_{ab},\pi_{bc},\pi_{ca}$ are permutations acting on $X$. By
the symmetry in the three colors, and since $|X|$ is finite, it is
enough to show that $\pi_{ab}$ is injective. Assume that
$\pi_{ab}(p)=\pi_{ab}(q)=r$ for $p,q,r\in X$. Vertex $r$ is the
final vertex of the $ab$-paths which start at $p$ and at $q$. But
there is only one $b$-edge incident to $r$. Thus, this $b$-edge is
used by both paths; we conclude that the other end of the $a$-edge
incident to $p$ is the same as the other end of the $a$-edge
incident to $q$. Denote this common end by $s$ (which is a vertex
in $Y$). The fact that there is only one $a$-edge incident to $s$
implies that $p=q$. We conclude that $\pi_{ab}$ is injective, and
thus a permutation acting on $X$.

Observe that the number of faces of the faithful embedding of $G$
is the sum of the numbers of cycles of the permutations
$\pi_{ab},\pi_{bc},\pi_{ca}$. Note also that the composition
$\pi_{ab}\circ\pi_{bc}\circ\pi_{ca}$ is the identical permutation,
hence an even permutation. It follows that the total number of
even cycles of the three permutations is even. (Recall that an
even cycle has odd parity, and vice-versa.) Thus, the parity of
the number of cycles of $\pi_{ab},\pi_{bc}$ and $\pi_{ca}$ is the
same as the parity of the total number of their odd cycles. Make
three copies of the vertices in $X$ and use each copy to represent
the cycles of one of the permutations, as polygons. The even
polygons use an even number of the $3k$ vertices. Therefore, the
parity of the total number of odd polygons is the same as the
parity of $3k$, that is, the same parity as $k$. Since we proved
that the parity of the number of faces in the faithful embedding
is the same as the parity of the total number of odd cycles of the
permutations $\pi_{ab},\pi_{bc}$ and $\pi_{ca}$, the first part of
the theorem is established.

Suppose $M$ is orientable. By definition, it follows that $C_M$ is
bipartite. Apply the first part of the theorem to the embedding of
$M$. Since the number of vertices of $C_M$ is a multiple of 4, it
follows that the number of faces of the faithful embedding of
$C_M$ is even. To compute $\chi(M)$ we have to add
$4|SQ(M)|-6|SQ(M)|$ to the number of faces, which we proved to be
even. The proof is complete. \fimproof

With this theorem we conclude the relations among Tutte's axioms
[25], the classical topological notion of a graph well-embedded in
a surface, and the definition that we use here. We showed that,
under the bijections $\psi^a$ and $\psi$, they are in fact
equivalent.

The next theorem can be used to give an upper bound on the number
of disjoint odd polygons which are subgraphs of $C_M$. We derive
this consequence in Lemma 4.8. It is also used in Theorem 2.14 to
give an upper bound on the maximum number of independent vertices
in special graphs.

\medskip\n{\bf (1.8) Theorem:} Assume that $G$ is a 3-edge colored cubic graph.
The $\mod\ 2$ sum of the sets of edges in an arbitrary subset $S$
of $ab$-, $bc$-, and $ca$-polygons, where $a,b,c$ are the colors
used to paint the edges of $G$, induces a disjoint set of polygons
of even length.
\medskip

\inicproof Since $G$ is cubic and the $\mod\ 2$ sum of the
edge-sets of polygons induces a cycle, it follows indeed that it
induces a disjoint set of polygons. We just have to prove that
each component induced by the $\mod\ 2$ sum of bicolored polygons
is a polygon with an even number of edges.

Consider an arbitrary such component, denoted K. Each edge in the
polygon $K$ is in an odd number of members of $S$, by definition
of $\mod\ 2$ sum. Observe that each edge of $G$ appears in
precisely two bicolored polygons. It follows that every edge of
$K$ appears in exactly one bicolored polygon in $S$. If an edge of
$K$ appears in an $ab$-polygon in $S$, let us mark it with
$\overline{c}$. Analogously, mark an edge of $K$ with
$\overline{a}(\overline{b})$ if it appears in a $bc$-polygon
($ca$-polygon) in $S$. This partitions $K$ into a number of
maximal paths whose edges have the same mark, which we say is the
label of the path. We observe that for $x,y$ different symbols in
$\{a,b,c\}$, if an edge of $K$ marked $\overline{x}$ is incident
to another marked $\overline{y}$, then the edge marked
$\overline{x}$ is an $y$-edge and the edge marked $\overline{y}$
is an $x$-edge. From this observation we derive the following
fact: the number of edges in a maximal path with the same label is
even if and only if the path that precedes it in $K$ and the path
that follows it have distinct labels.

The problem, then, is reduced to the following: if the vertices of
a polygon (with at least three edges) are properly colored in at
most three colors, then the number of vertices which have
neighbors of the some color is even.

The vertices correspond to the maximal paths. Therefore, the
vertices which have neighbors of the same color correspond to the
maximal paths with an odd number of edges. If there are just one
or two maximal paths, the theorem is trivial.  We deal with the
new problem in the next lemma.

\medskip\n{\bf (1.8.1) Lemma:} Suppose the vertices of a polygon, with at
least three vertices, are colored with at most three colors such
that adjacent vertices have different colors. Then the number of
vertices which have neighbors of the same color is even.
\medskip

\inicproof Let us denote the three colors by $0,1$, and $2$.
Orient the edges linking a vertex of color $i$ to a vertex of
color $i+1$ ($\mod\,3$) from the former to the latter. It follows
that the vertices which become a sink or a source with this
orientation are precisely the vertices which have neighbors of the
same color. The lemma follows if we observe that the number of
sources is the same as the number of sinks. \fimproof

The idea to orient the edges, which makes very simple the above
proof, is due to H. Shank in a personal communication.

\section*{1.C\quad The Phial Map and the Antimap}
In this section we use our definition of map to derive a richer
structure, which is the basis of the next chapter.

Given a map $M=(C_M,V_M,f_M)$, let us denote by $Q_M$ the
4-regular graph obtained from $C_M$ by adding to its edges the
$1$-factor $z_M$ in the complement of $C_M$, defined as the
diagonals of the squares.

By leaving out one of the three disjoint perfect matchings
$v_m,f_M,z_M$ of $Q_M$ and permuting the order of the other two,
we get a set of six maps, denoted by $\Gamma(M)$, and named as
follows:
\begin{center}
\begin{tabular}{ccl}
{\bf SYMBOL} & {\bf DEFINITION} & \multicolumn{1}{c}{\bf NAME}\\
    $M$  & $(Q_M\backslash z_M,v_M,f_M)$ & map $M$ \\
 $D_M=D$ & $(Q_M\backslash z_M,f_M,v_M)$ & dual of $M$ \\
$D^\sim_M=D^\sim$ & $(Q_M\backslash v_M,f_M,z_M)$ & antidual of $M$ \\
$P_M=P$  & $(Q_M\backslash v_M,z_M,f_M)$ & phial of $M$ \\
$P^\sim_M=p^\sim$ & $(Q_M\backslash f_M,z_M,v_M)$ & antiphial of $M$ \\
$M$      & $(Q_M\backslash f_M,V_M,z_M)$ & antimap of $M$
\end{tabular}
\end{center}

Observe that if $X\in\Gamma(M)$, then $\Gamma(X)=\Gamma(M)$.
Observe also that $D^t$ (the $t$-map associated with $D$)
corresponds to the usual dual $t$-map, which interchanges faces
and vertices relative to $M^t$.

Consider a $t$-map $M^t$. We now define some important cyclic
paths which are analogous to the facial paths. Starting at an
arbitrary point of the surface near an edge we proceed (in the
surface) Parallel to the edge, crossing it in its middle, and
continuing to go parallel to the edge and in the same direction as
before, but now in the other ``side'' of the edge. Near a vertex
we turn (without crossing points of the graph) and repeat the
process with the new edge. Iteration of the above rules ends when
we get back to the initial edge in the initial ``side''. Clearly,
we reach this stage and we reach it before we reach a ``side'' of
an edge already traversed. The edges which indicates the way and
the vertices which make the transitions between two edges must be
recorded such as to form a cyclic path in $G_M$. Such cyclic path
is called a {\it zigzag} in $M^t$.

Observe that the zigzag paths form another combinatorial embedding
of $G_M$, since each edge is used twice in the collection of
zigzag paths. This embedding corresponds to $M^\sim$. We claim
that $M^t$ and $(M^\sim)^t$ share a duality relation: their facial
paths and zigzag paths are interchanged, while the (pairs of
inverses) cyclic sequences of edges around the vertices of their
graphs are the same. This claim can be easily verified in the maps
$M$ and $M^\sim$: their $f$-gons and $z$-gons are interchanged,
while the $v$-gons are the same.

The $z$-gons of maps $M$ and $D$ are the same. This corresponds to
the fact that the cyclic sequence of edges in the zigzag paths of
$M^t$ and $D^t$ are the same.

$M$ and $P$ clarify a duality relation between $M^t$ and $P^t$.
The phial pair of maps $M$ and $P$ have $v$-gons and $z$-gons
interchanged, while the same $f$-gons. This means that $M^t$ and
$P^t$ have the cyclic sequence of edges around the vertices and
the cyclic sequence of edges in the zigzag paths interchanged,
while having the same cyclic sequence of edges in the facial
paths.

The zigzag paths are studied by Shank [24], in the restricted case
of planar $t$-maps. In this context he names these paths
``left-right paths''.  Paper [24] Influenced the definition of
zigzag paths and the subsequent development of the theory
presented in Chapters 2 and 3.

We show in Figures 1.1, 1.2, and 1.3 the elements of $\Gamma(M)$,
where $M$ is the $t$-map formed by a cube embedded in the plane.
The dual $t$-map, $D^t$, is the familiar octahedron embedded in
the plane. Observe that $S_{D^\sim}=S_P$, which is the
non-orientable surface of connectivity 4, as shown by the
non-orientable handles labelled $A$ and $B$ in Figure 1.2. Finally
$S_{P^\sim}=S_{M^\sim}$, which is the familiar torus. The reader
can observe the relations among the edges around the vertices, the
edges in the facial paths, and the edges in the zigzag paths.
\medskip
$$\includegraphics[width=59.576mm,height=58.806mm]{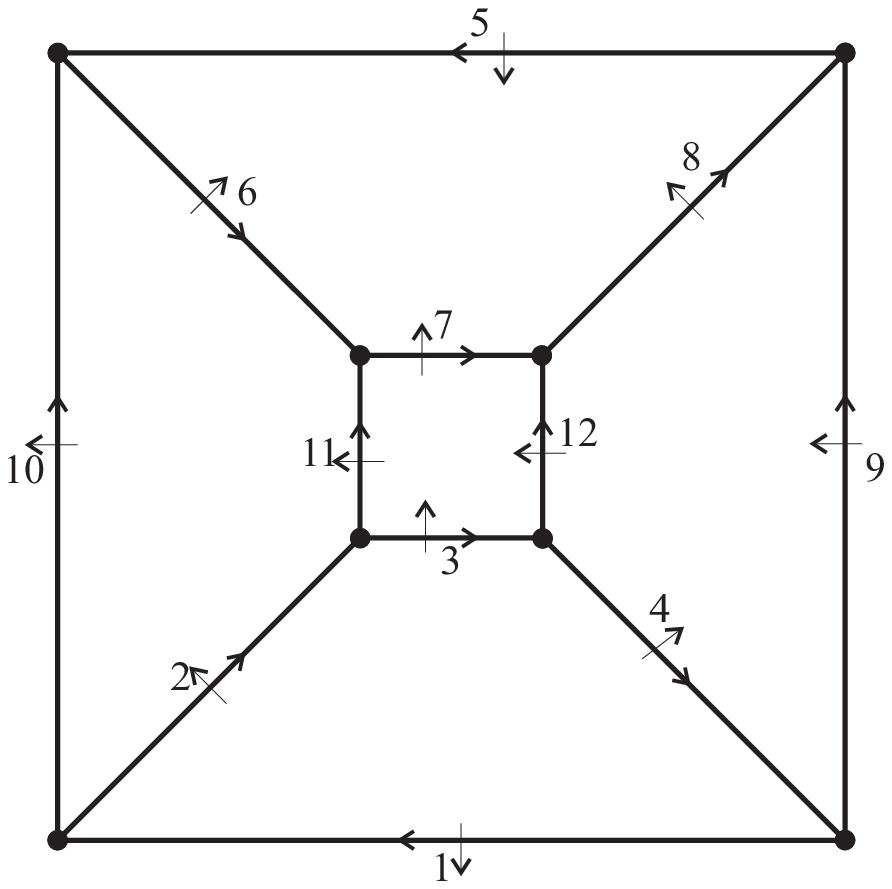}\hspace{1cm}
\includegraphics[width=61.422mm,height=54.153mm]{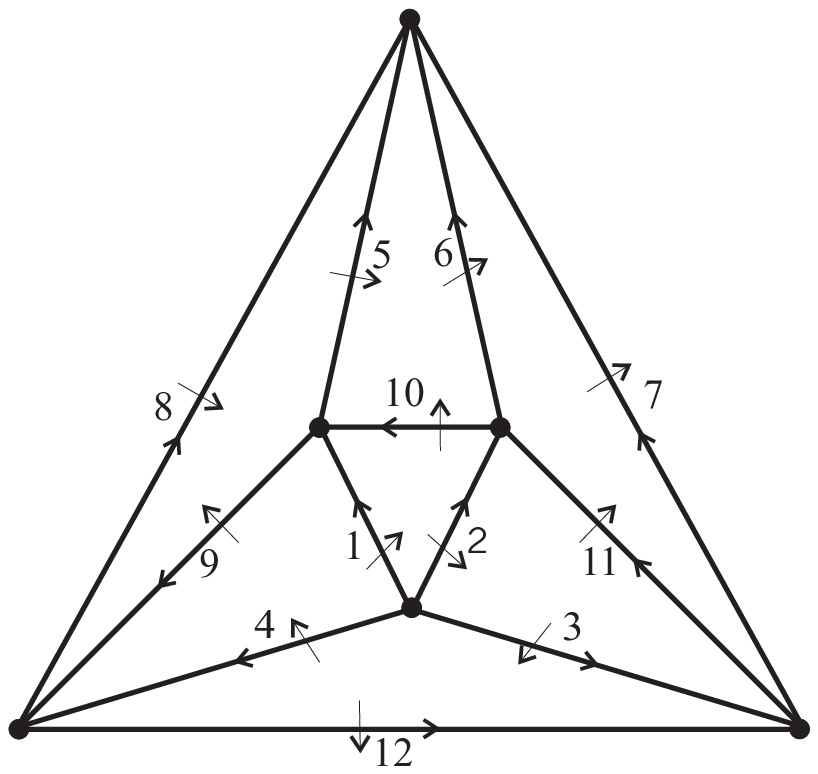}
$$
\centerline{Figure 1.1}
\bigskip

$$\includegraphics[width=85.071mm,height=59.589mm]{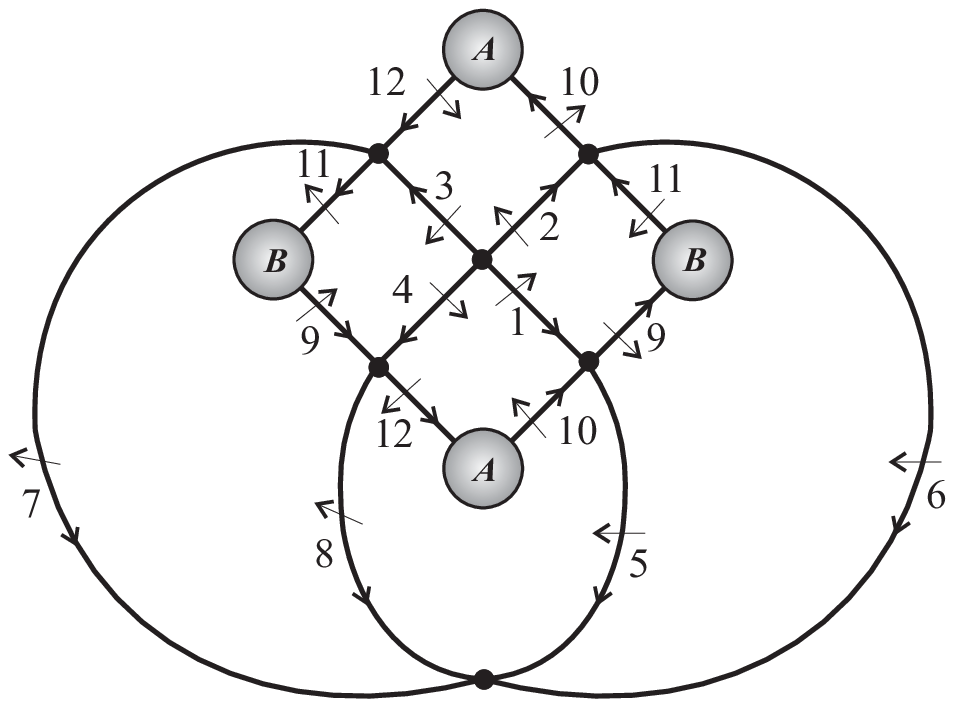}\hspace{1cm}
\includegraphics[width=54.053mm,height=54.558mm]{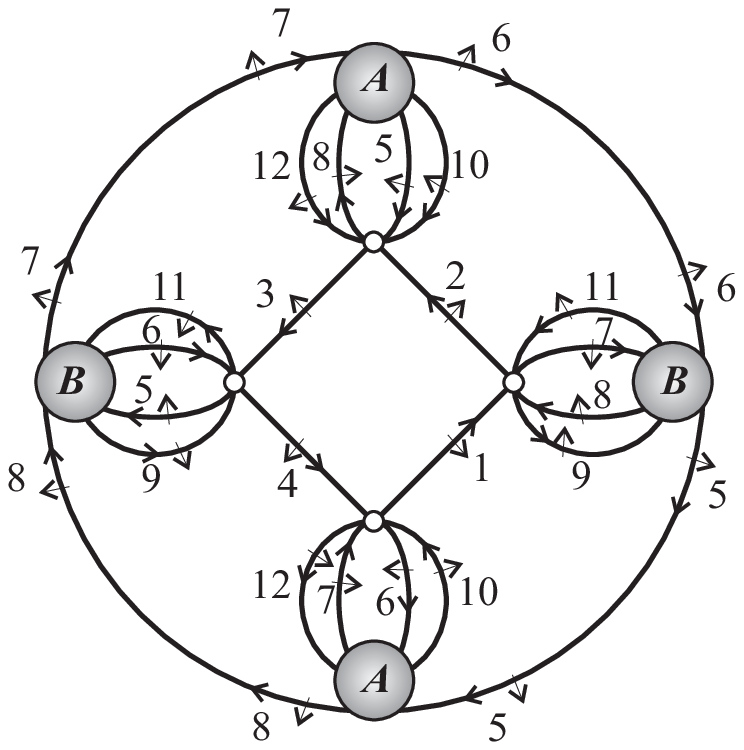}
$$
\centerline{Figure 1.2}

\eject

$$\includegraphics[width=57mm,height=57mm]{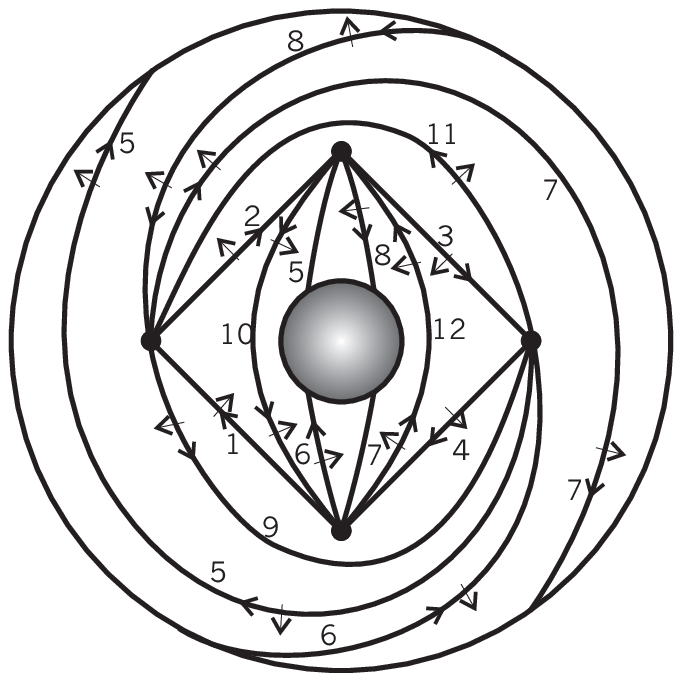}\hspace{1cm}
\includegraphics[width=57mm,height=57mm]{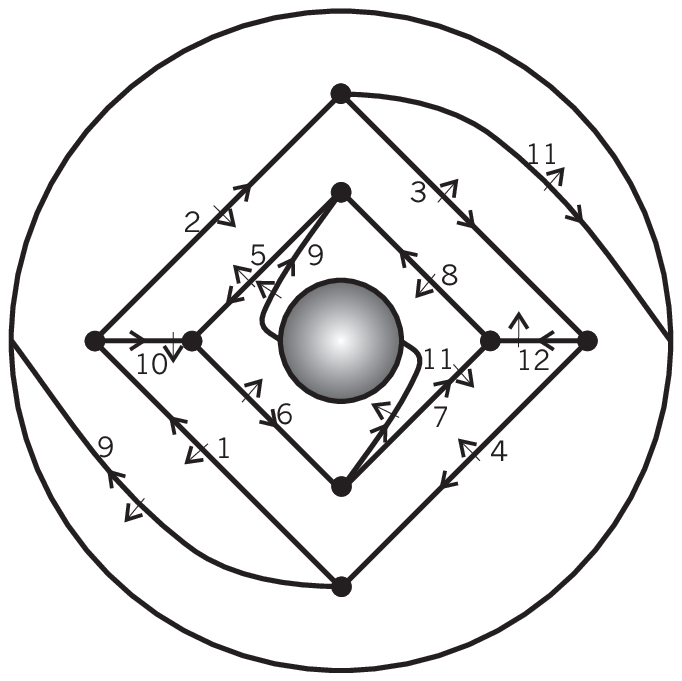}$$
\centerline{Figure 1.3}

From the definitions of its members, there exists a natural
bijective correspondence between the squares of the elements of
$\Gamma(M)$.

Two maps $M=(C_M,v_M,f_M)$ and $N=(C_N,v_N,f_N)$ are {\it
isomorphic} if there exists a bijection $b$ from $VC_N$ onto
$VC_M$ such that for every $x\in VC_M$,
$$\begin{array}{c}
bv_M(x)=v_Nb(x) \\ [8pt] bf_M(x)=f_Nb(x) \\ [8pt] ba_M(x)=a_Nb(x)
\end{array}$$

Two members of $\Gamma(M)$ are said to be {\it square-identically
isomorphic}\index{\it square-identically isomorphic} if they are
isomorphic and the isomorphism preserves the vertices of each
square.

We use the structure of $\Gamma(M)$ to prove that the graphs which
arise by doubling every edge of another graph can be embedded so
as to form a self-dual $t$-map in a strong sense, which we now
describe. Denote by $G^d$ the graph obtained from a graph $G$ by
doubling every edge of $G$. Consider an arbitrary $t$-map $E^t$
such that $G_E=G$. Denote by $D^t$ the $t$-map obtained from $E^t$
by replacing every edge of $G$ by a bounding digon. Let $M$ denote
the dual of $D$. Consider the pair of dual maps\index{dual maps}
formed by the phial $P$ of $M$ and the antidual $D^\sim$ of $M$.
Under these conditions, we show in the next theorem that $P$ and
$D^\sim$ are square-identically isomorphic. It is an easy
property, which we prove in 2.1, that $G_D\cong G_{D^\sim}$ for
every map $D$. Consequently, $G^d\cong G_D\cong G_{D^\sim}$ is
embedded such as to be {\it edge-identically}
isomorphic\index{{\it edge-identically} isomorphic}, i.e., the
isomorphism between $G_P$ and $G_{D^\sim}$ induces the same 1-1
correspondence between $eG_p$ and $eG_{D^\sim}$ as the one given
by the topological duality.

\medskip\n{\bf (1.9) Theorem:} Let $D^t$ be obtained from $E^t$ by doubling each
edge of $G_E$ so as to form a bounding digon in $S_E$. Let $M$ be
the dual of $D,P$ be the phial of $M$, and $D^\sim$ be the antimap
of $D$. Then the pair of dual maps, $P$ and $D^\sim$ are
aquare-identically isomorphic.
\medskip

\inicproof It is sufficient to prove that $M$ and $M^\sim$ are
square-identically isomorphic. This is so because $P$ is the phial
of $M$ and $D^\sim$ is the phial of $M^\sim$.

Denote by $N$ the dual of $E$. Observe that $G_M$ is equal to
$G_N$ with every edge subdivides by a bivalent vertex, called a
{\it middle} vertex of $G_M$.

It is a straightforward property of a map $M$, whose graph $G_M$
is obtainable by subdividing every edge of another graph by a
middle vertex, that $M$ and $M^\sim$ are square-identically
isomorphic. To verify this, recall that the antimap of a map $M$
replaces, in each square of $M$, the two $f_M$ edges by its
diagonals. (See Figure 1.4). To be convinced that the map and its
antimap under the hypothesis are square-identically isomorphic,
just twist the $v$-gons with four edges which correspond to the
middle vertices. \fimproof
$$\includegraphics[width=60.424mm,height=24.85mm]{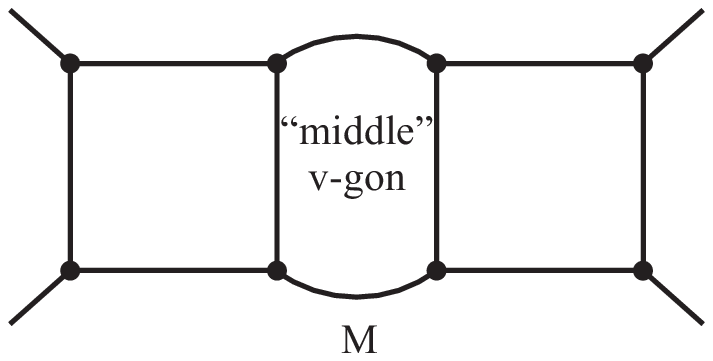}\hspace{2cm}
\includegraphics[width=60.424mm,height=24.85mm]{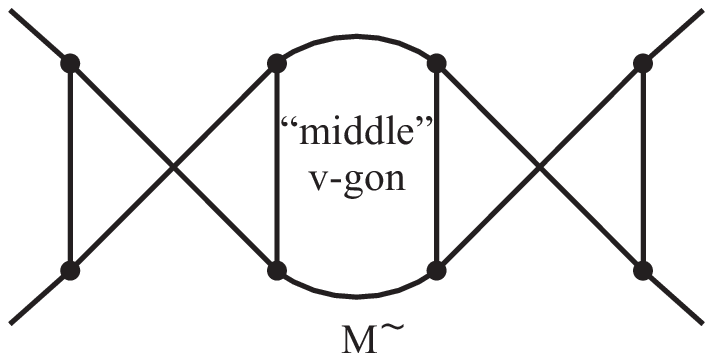}$$
\centerline{Figure 1.4}
\chapter{Rich Maps and Graphs}
\section*{2.A\quad Algebraic Preliminaries}
Initially we state the facts from linear algebra that we use.

Consider a finite dimensional vector space $W$ over a field $K$
and a basis $\{a_i\}_{i=1,\dots,n}$ of $W$. Let the symmetric
bilinear function
$$\aw\cdot,\cdot\fw:\ W\times W\rightarrow K$$
be defined by $\aw u,v\fw=\sum^n_{i=1}\lambda_i\gamma_i$, where
$\lambda_i,\gamma_i\in K$, for $i\in\{1,\dots,n\}$, and
$u=\sum^n_{i=1}\lambda_ia_i,\ v=\sum^n_{i=1}\gamma_ia_i$. Function
$\aw\cdot,\cdot\fw$ is the well-known scalar product induced by
the basis $\{a_i\}_{i=1,\dots,n}$. Relative to the scalar product,
there is the concept of orthogonality: two vectors $u,v\in W$ are
{\it orthogonal}\index{\it orthogonal} if $\aw u,v\fw=0$. If $A$
is a subspace of $W$, then the set of vectors of $W$ which are
orthogonal to all vectors in $A$ is also a subspace of $W$,
denoted by $A^\perp$.

The following is a list of elementary algebraic properties,
involving the above concepts. Proofs of these properties can be
found in Godement's book [8, Sections 19, 36].

\n{\bf(2.0.1)} For a finite dimensional vector space $W$, over an
arbitrary field, let $A$ and $B$ denote arbitrary subspaces of
$W$, and let ``$\perp$'' be defined relative to a scalar product
induced by a fixed basis for $W$. The following properties
hold:\jtt
\begin{itemize}
\item[(a)] $\dim(A+B)=\dim(A)+\dim(B)-dim(A\cap B)$; \item[(b)]
$(A^\perp)^\perp=A$; \item[(c)] $(A+B)^\perp=A^\perp\cap
B^\perp;\quad (A\cap B)^\perp=A^\perp+B^\perp$; \item[(d)]
$\dim(A)+dim(A^\perp)=\dim\ W$. \fimproof
\end{itemize}\jtt

We also use the fundamental homomorphism theorem (for vector
spaces): assume that $W_1$ and $W_2$ are vector spaces over the
same field, and $h:W_1\rightarrow W_2$ is a homomorphism. Then,
the quotient space $W_1/{\rm Ker}(h)$ is
isomorphic\index{isomorphic} to the subspace of $W_2$ defined by
image of $W_1$ under $h$. In particular,
$$\dim[{\rm Ker}(h)]+\dim[Im(h)]=\dim\ W_1.$$\jtt
\fimproof

If $W$ is a subset of vertices of a graph $G$, the {\it
coboundary} of $W$, denoted by $\delta_G(W)$ or $\delta(W)$, is
the set of edges of $G$ which have one end in $W$ and the other in
$VG-W$.

Let $M$ be a given map. Consider the graphs $G_M,G_D$, and $G_p$
The sets of coboundaries of these graphs are denoted by $V_M,F_M$,
and $Z_M$ respectively. To simplify the notation, we use $V,F$,
and $Z$, if $M$ is understood.

We give, as usual, a vector structure to the power set $2^{SQ(M)}$
by taking the sum of subsets to mean their $\mod\ 2$ sum.

Thus, the underlying field is $GF(2)$. Let the canonical basis for
$2^{SQ(M)}$ be the set of singletons representing the individual
squares of $M$.

With the choice of this basis we have a scalar product and the
notion of perpendicular\index{perpendicular} subsets of squares of
$M$, (or subsets of edges of $G_M,G_D$ or $G_P$) in this special
case two subsets are {\it perpendicular} if their intersection has
even cardinality.

Since the symmetric difference of two coboundaries in a graph $G$
is also a coboundary in $G$, the set of coboundaries of $G$,
denoted $BS(G)$, is a subspace of $2^{SQ(M)}$. Analogously, since
the symmetric difference of two cycles in $G$ is also a cycle in
$G$, the set of cycles of $G$, denoted as said before by $CS(G)$,
is a subspace of $2^{SQ(M)}$. For every graph $G$ the intersection
of any member of $BS(G)$ and any member of $CS(G)$ has even
cardinality. Indeed, $BS(G)$ and $CS(G)$ are perpendicular
subspaces, for every graph $G$. See Bondy and Murty [1]. The
dimensions of these subspaces are
$$\begin{array}{l}
\dim[CS(G)]=|eG|-|VG|+p \\ [5pt] \dim[BS(G)]=|VG|-p,
\end{array}$$
where $p$ is the number of components of $G$. See also [1].

In view of the above relation between $BS(G)$ and $CS(G)$, we
denote by $V^\perp,F^\perp$ and $Z^\perp$, respectively, the
cycles spaces of $G_M,G_D$, and $G_P$.

We use the dimensions of cycle and coboundary spaces of graphs and
the theory of 3-colored cubic graphs to establish the familiar
upper bound for $\chi(M)$, for connected maps $M$. We also show
that if the upper bound is attained, then $C_M$ is bipartite. This
shows that in the definition of planar maps, the condition that
$C_M$ is bipartite is redundant.

\medskip\n{\bf (2.0.2) Theorem:} The maximum Euler Characteristic of a
connected map is 2. If $\chi(M)=2$, then $M$ is orientable.
\medskip

\inicproof We prove a more general statement. Suppose $G$ is a
connected cubic graph and $\{a,b,c\}$ is the set of colors of a
proper 3-coloration of the edges of $G$. Let $2k$ be the number of
vertices of $G$; consequently, $3k$ be the number of its edges.
Denote by $p(x,y)$ the number of $xy$-polygons, $x,y\in\{a,b,c\}$.
We prove that $2k+p(a,b)+p(b,c)+p(c,a)-3k$ is at most two, and
that if this bound is attained, $G$ is bipartite. Specializing for
the faithful embedding of $C_M$, we have the theorem.

Observe that the dimension of the cycle space of $G$ is $3k-2k+1=
k+1$. The dimension of the subspace of $CS(G)$ generated by the
$ab-,bc-$, and $ca$-polygons is $p(a,b)+p(b,c)+p(c,a)-1$. To see
this, observe that this subspace is the space of the coboundaries
of the topological dual of the $t$-map formed by the faithful
embedding of $G$, and that this dual has its vertices in bijective
correspondence with the $ab-$, the $bc-$, and the $ca$-polygons.
Therefore, we have
$$p(a,b)+p(b,c)+p(c,a)-1\leq k+1,$$
that is,
$$p(a,b)+p(b,c)+p(c,a)-k\leq 2,$$
as we stated. If we have equality, then the three types of
polygons generate $CS(G)$. As every such generator has an even
number of edges, it follows that $G$ is bipartite. This concludes
the proof. \fimproof

The intersection of $BS(G)$ and $CS(G)$, subspace of $2^{eG}$ over
$GF(2)$, has been named the \index{\it bicycle space} of $G$ (see
[21]).

The dimension of the bicycle space of a graph is an important
parameter in the theory that follows.

\section*{2.B\quad A Ternary Relation Involving Cycles and Co\-boundaries}
In this section we introduce the concept of rich map\index{rich
map}. A rich map induces three graphs on the same set of edges
having the property that an arbitrary cycle in one of them is
expressible as the symmetric difference of two coboundaries, one
in each of the two other graphs.

\medskip\n{\bf (2.1) Proposition:} For any map $M$ we have $G\cong G_{M^\sim}$.
\medskip

\inicproof The $M$-squares, which are the edges of $G_M$ are in
natural 1-1 correspondence with the $M^\sim$-squares which are the
edges of $G_{M^\sim}$; namely, to $M$-square $S$ corresponds the
$M^\sim$-square $S^\sim$ such that $eS\cup eS^\sim$ induces in
$Q_M$ a subgraph isomorphic to $K_4$. The vertices of $G_M$ are
the $v$-gons of $M$, which are the $v$-gons of $G_{M^\sim}$.
Therefore, the sets of vertices of $G_M$ and of $G_{M^\sim}$ are
equal. Also, the ends of an $M$-square (i.e., and edge of $G_M$)
are the same $v$-gons which are the ends of the corresponding
$M^\sim$-square (corresponding edge of $G_{M^\sim}$). This
concludes the proof of the proposition. \fimproof

\medskip\n{\bf (2.2) Proposition:}
For every map $M$, the cycle space of one of $G_M,G_D,G_P$
contains the sum of the coboundary spaces (over $GF(2)$) of the
other two.
\medskip

\inicproof By the symmetry in the three graphs given by the
implication $X\in\Gamma(M)\Rightarrow\Gamma(X)=\Gamma(M)$ and by
the last proposition we simply have to prove that
$V^\perp\supseteq F+Z$. Indeed, the proof for arbitrary map $M$
that $V^\perp\supseteq F$ suffices, because for $M^\sim,\
V^\perp_{M^\sim}=V^\perp$ and $F_{M^\sim}=Z$.

To prove that $F\subseteq V^\perp$, for arbitrary map $M$, we show
that the elements of a set generating $F$ are elements of
$V^\perp$. To this aim, it is sufficient to prove that the set $T$
of squares of $M$ which have one $f_M$-edge in an arbitrary
$f$-gon corresponds to a cycle in $G_M$. This is straightforward
because $T$ is the image of the $f$-gon under $\Psi^M_c$. The
reason for the sufficiency is that $T$ is a coboundary of the
vertex in $G_D$ which corresponds to (in fact, is) an $f$-gon. The
proof of the proposition is complete. \fimproof

The above proposition makes the following definition meaningful.
The {\it (cycle) deficiency} of a map\index{deficiency of a map}
$M$ is the dimension of the quotient space $V^\perp/(F+Z)$. A map
is called {\it (cycle) rich} if
$${\rm def}(M)=\dim\left(\frac{V^\perp}{F+Z}\right)=0.$$
$$\begin{array}{c}
\includegraphics[width=51.353mm,height=56.183mm]{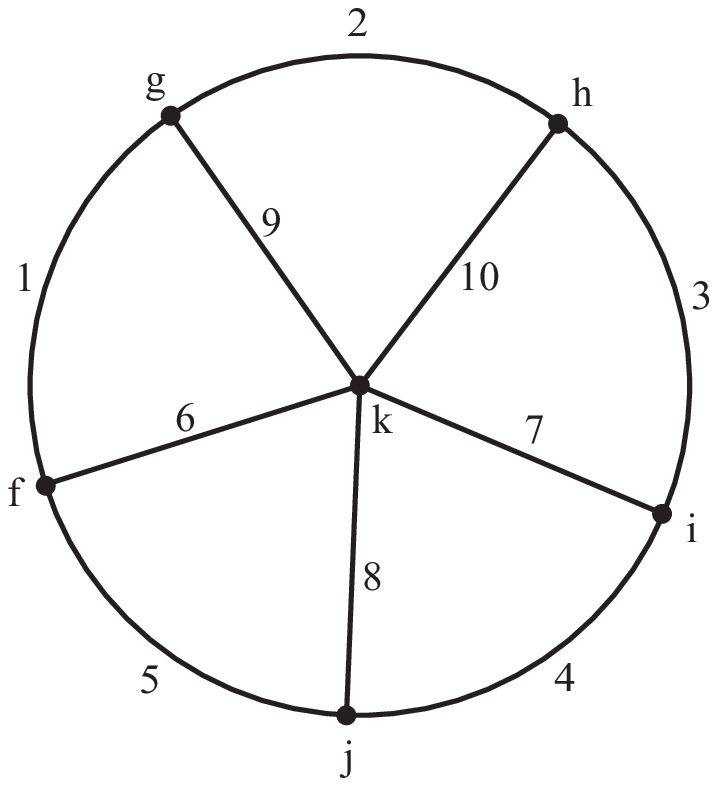}\hspace{2cm}
\includegraphics[width=55.353mm,height=55.353mm]{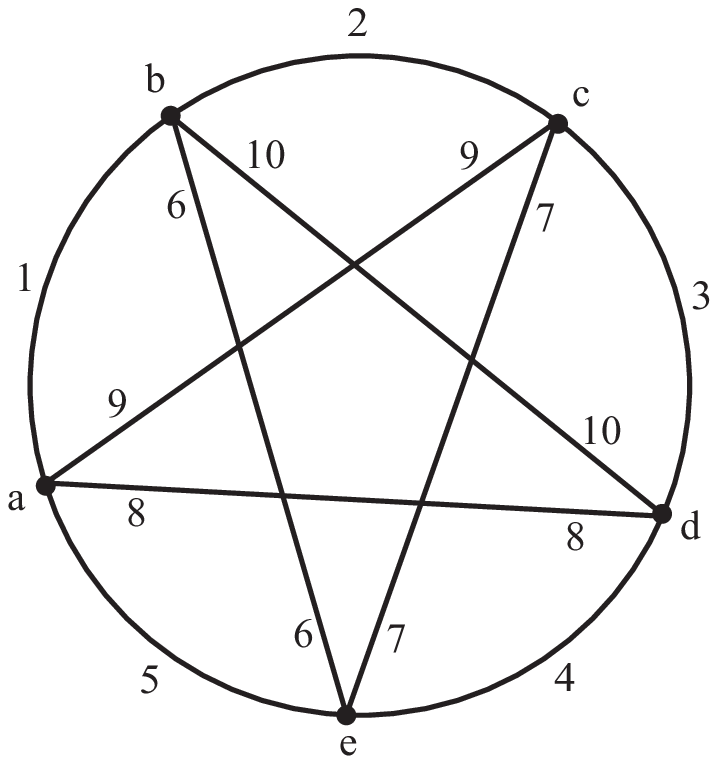}\\
\includegraphics[width=22.641mm,height=43.726mm]{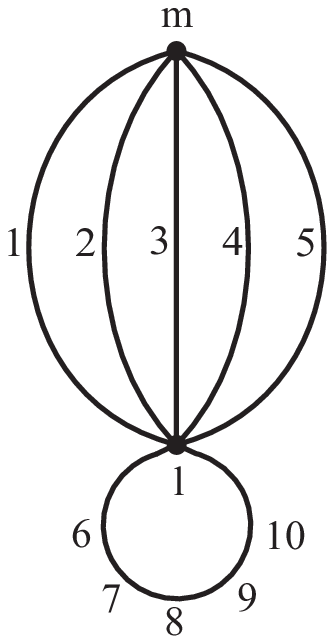}
\end{array}
$$
\centerline{Figure 2.1}

The three graphs presents in Figure 2.1 are $G_M,G_D$, and $G_P$
for a rich map $M$. A cycle in one of the graphs can be expressed
as the $\mod\ 2$ sum of coboundaries in the other two. We
exemplify this relation as follows:\jtt
\begin{itemize}
\item[(1)] the circuit $\{5,6,8\}$ in the first graph can be
written as
$$\delta(\{b,d\})+\delta(\{\ell\})=\{1,2,3,4,6,8\}+\{1,2,3,4,5\}=\{5,6,8\};$$
\item[(ii)] the circuit $\{2,6,7\}$ in the second graph can be
written as
$$\delta(\{f,i\})+\delta(\ell)=\{1,3,4,5,6,7\}+\{1,2,3,4,5\}=\{2,6,7\};$$
\item[(iii)] the circuit $\{6\}$ in the third graph can be written
as
$$\delta(\{f,h,i\})+\delta(\{a,c,d\})=\{1,2,4,5,7,10\}+\{1,2,4,5,6,7,10\}=\{6\}.$$
\end{itemize}

The explanation for this ternary symmetric relation is given by
the next theorem.

\medskip\n{\bf (2.3) Theorem.} {\bf (Cycle-coboundary symmetric ternary
relation):}\index{Cycle-coboundary symmetric ternary relation} If
$M$ is a rich map, then the three graphs $G_M,G_D$, and $G_P$
satisfy the following property: the cycle space (over $GF(2)$) of
one of them is equal to the sum of the coboundary spaces of the
other two.
\medskip

\inicproof Assume the Lemma 2.4 on equality of the
deficiencies\index{equality of the deficiencies}. The above
theorem is a direct corollary of this lemma with zero as the
common deficiencies of $M,D$, and $P$. \fimproof

\medskip\n{\bf (2.4) Lemma.} {\bf (Equality of deficiencies):} For a map $M$, every
two members of  $\Gamma(M)$ have the same deficiency.
\medskip

\inicproof By the Absorption Property, proved as Theorem 2.5
below, applied to $P$ in place of $M$, it follows that $V_P\cap
F_P\subseteq Z_P$. But $V_P=Z_M=Z,\ F_P=F_M=F,\ Z_P=V_M=V$. Thus,
$Z\cap F\subseteq V$ or $F\cap Z=V\cap F\cap Z$.
$$\begin{array}{rcl}
{\rm def}\,M&=&\dy\dim\left(\frac{V^\perp}{F+Z}\right)=\dim\,V^\perp-\dim(F+Z)\\
[8pt] &=& (e-v+p)-(f-p)-(z-p)+\gamma,
\end{array}$$
where $e$ is the number of squares of $M$; $v,f,z$ are,
respectively, the numbers of vertices of $G_M,G_D,G_P$; $p$ is the
number of components of $M$; and $\gamma$ is the dimension of
$V\cap F\cap Z$.

The above expression can be rewritten as
$${\rm def}(M)-(e+3p+\gamma)-(v+f+z),$$
which is evidently symmetric in $v,f,z$. Since every element of
$\Gamma(M)$ only induces a permutation of $v,f,z$ in the above
expression, the lemma follows, provided the Absorption Property
(Theorem 2.5) is proved. \fimproof
$$\includegraphics{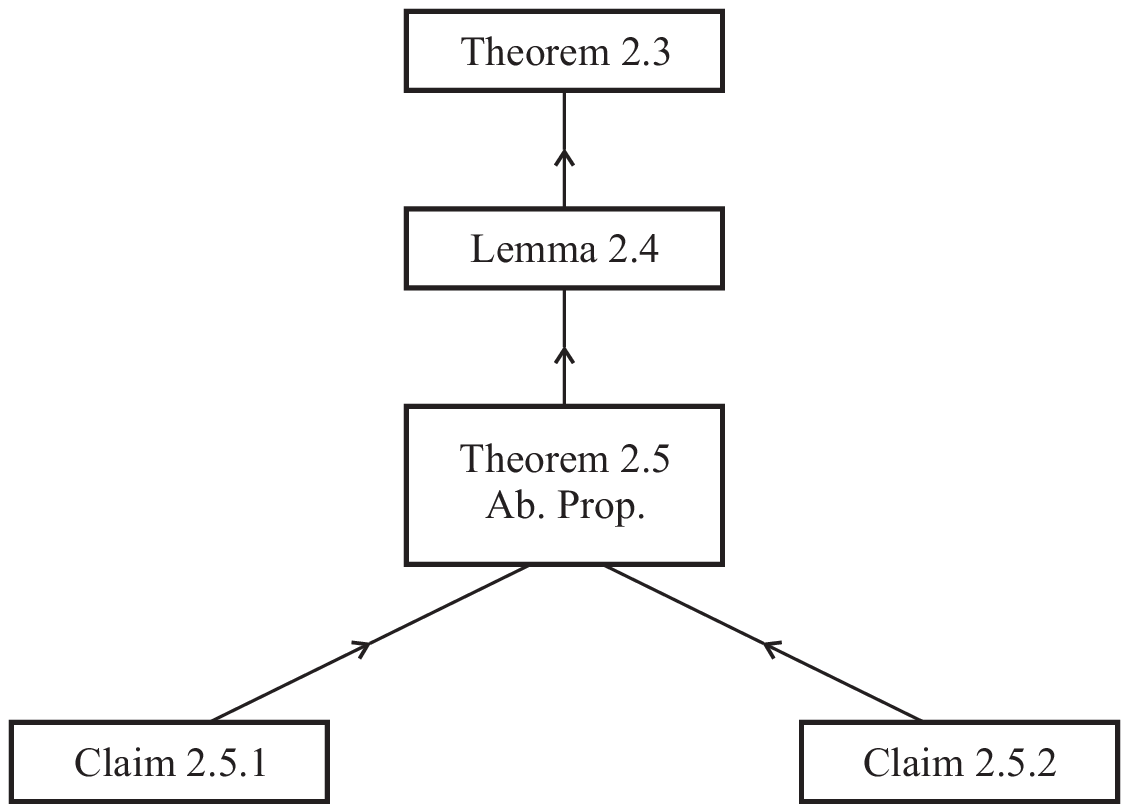}$$%%[width=7cm,height=5cm]
\centerline{Dependence Structure for Theorem 2.3}

The next theorem is the key for the symmetry in this theory.

\medskip\n{\bf (2.5) Theorem.} {\bf(Absorption Property):}\index{Absorption Property} For arbitrary map $M$ we
have $V\cap F\subseteq Z$.
\medskip

\inicproof Take an element $b\in V\cap F$. It follows that there
exist $S\subseteq vG(M)$ and $T\subseteq fG(M)$ such that $b$ is
the coboundary of $S$ in graph $G_M$ and the coboundary of $T$ in
graph $G_D$. Take the subgraph of $C_M$ which consists of the
$v$-gons in $S$, and denote it by $C_S$. Consider also the
subgraph $C_T$, defined analogously by the $f$-gons in $T$. Denote
by $R$ the subgraph of $C_M$ induced by the symmetric difference
of $eC_S$ and $eC_T$. We claim the following two facts about $R$:

\n{\bf(2.5.1)} For every corner $x$ of $M$, the $a_M$-edge
incident to $x$ is in $eR$ if and only if the $a_M$-edge incident
to $z_M(x)$ (the corner opposite to $x$ in its square) is also in
$eR$.

\n{\bf(2.5.2)} A square of $M$ has exactly one pair of opposite
$a_M$-edges which are incident to it in $eR$ iff the square
corresponds to an edge of $G_M$ (and of $G_D$) which is in $b$.

Let us assume the two claims and conclude the proof of the
Absorption Property. Following that we establish the claims.

By claim (2.5.1), graph $R$ defines naturally a set $U$ of
$v$-gons in map $P=P_M$. Specifically, the subgraph $C_U$ of $C_P$
is formed from $R$ by deleting every $v_{M^-}$ and $f_{M^-}$ edge
of $R$ and including a $z_M$-edge, whenever it links a pair of
corners such that the $a_M$-edges which are incident to them are
in $eR$. Thus, $C_U$ is formed by a subset $U$ of $VG(P)$. By
claim 2.5.2 the coboundary of $U$ in $G_p$ is $b$. Hence $b\in Z$,
concluding the proof of the theorem, subject to the proofs of the
claims (2.5.1) and (2.5.2). \fimproof

\noi{\bf Proof of the claims (2.5.1) and (2.5.2)} Take a square
$X$ of $M$. Consider that $X$ has $i$ $v_M$-edges in $eC_S$ and
$j$ $f_M$-edges in $eC_T$. Both $i$ and $j$ are in $\{0,1,2\}$.
Since $b\in V\cap F$, it follows that $|eX\cap(eC_T+eC_S)|$ is
even for every square $X$ of $M$. Hence, we have
$$(i,j)\in\{(0,0),(1,1),(2,0),(0,2)\}.$$

Claim 2.5.1 is directly checked for the four cases above. Thus,
(2.5.1) follows, since $X$ is arbitrary.

Square $X$ corresponds to an edge in $b$ iff $(i,j)=(1,1)$. Direct
observation shows that of the four cases the only one which
implies that exactly one pair of opposite $a_M$-edges is included
in $eR$ is also $(i,j)=(1,1)$. Hence, the proof of (2.5.2) is
complete. \fimproof

Observe that with the proof of these claims every statement used
in the proof of Theorem 2.3 is established.

For a map $M$ we define the {\it facial deficiency}\index{\it
facial deficiency} of $M$ as the dimension of $V^\perp/F$.

\medskip\n{\bf (2.6) Proposition:} The facial deficiency of a connected map $M$ is $\xi(M)$.
\medskip

\inicproof The proof follows from the definition:
$$\begin{array}{rcl}
\dy\dim\left(\frac{V^\perp}{F}\right)&=&\dim(V^\perp)-dim(F)\\
[5pt] &=& (e-v+1)-(f-1) \\ [5pt] &=& (2+e)-(v+f) \\ [5pt] &=&
2-\chi(M)=\xi(M),
\end{array}$$
where $e=|SQ(M)|,\ v=|vG(M)|,\ f=|fG(M)|$.\fimproof

Part (b) of the following corollary appears in Shank [24], proved
by more specific techniques.

\noi{\bf Corollary of (2.6)} (a) Planar maps are rich; (b) the
bicycle space of a planar graph $G$ is generated by the cycles
induced by the zigzag paths of $M^t$, where $G_M=G$.

\inicproof That a planar map\index{planar map} $M$ is rich follows
from the fact that $\xi(M)=0$, which implies, by the previous
proposition that $V^\perp=F$, or $V^\perp=F+Z$.

To establish (b) it is sufficient to show that $V\cap V^\perp=Z$.
By part (a) and Lemma 2.4, it follows that $V+F=Z^\perp$.
Replacing $F$ by $V^\perp$ and taking the perpendiculars on both
sides of the equality yields $V\cap V^\perp=Z$, by 2.0.1(c). This
concludes the proof of the corollary. \fimproof

Planar maps are simple examples of rich maps, simple in the sense
that we do not need $Z$ to generate $V^\perp$. What we propose to
do now is to characterize, in a topological way, rich maps. With
the characterization it is easy to get examples of rich maps in
arbitrary surfaces.

Motivated by Theorem 2.3, concerning the symmetric ternary
relation in cycles and coboundaries, we introduce the following
definition.

A graph is called {\it(cycle) rich} if there are two other graphs
with the same set of edges as the first graph such that a cycle in
one of the three graphs corresponds to the symmetric difference of
two coboundaries, one in each of the other two graphs.

We define also the {\it map-deficiency} of a graph\index{{\it
map-deficiency} of a graph} $G$ to be $\min_M\{{\rm
def}(M):G=G_M\}$. A graph $G$ is called {\it map-rich} if its
map-deficiency is 0, that is, if there exists a rich map $M$ with
$G_M=G$.

Theorem 2.3 implies that all map-rich graphs\index{map-rich
graphs} are rich graphs. Having no counter-examples so far, we
conjecture the converse:

\noi{\bf Conjecture:} All rich graphs are map-rich graphs.

This conjecture if true, would complete the analogy between rich
graph and planar graph.

Since we have not made any progress to settle or disprove this
conjecture, we consider the easier problem for the
characterization of map-rich graphs.

\section*{2.C\quad Surface Characterizations of Map-Rich\newline Graphs}
This section establishes two characterizations (Theorem 2.9) of
map-rich graphs. One of these is used to give a topological
``recipe'' for all map-rich graphs. The other is used to prove
that there are graphs which are not map-rich.

A {\it medial map}\index{\it medial map} is a $t$-map $X=(H,S)$
such that $H$ is 4-regular and the faces of $X$ are partitioned
into two ordered classes called ``black'' and ``white'' such that
every edge of $H$ is incident to faces of different colors.

We define a function acting on the set of maps with image in the
set of $t$-maps: $\psi'(M)$ is the $t$-map obtained from the
embedding of $M$ by contracting to a point in $S_M$ each disc
bounded by a square of $M$. We use the notation
$\psi'(M)=M'=(H_M,S_M)$. The next lemma shows that this notation
is convenient.

\medskip\n{\bf (2.7) Lemma:} The function $\psi'$ is a bijection from the set of  maps
onto the set of medial maps.
\medskip

\inicproof First we prove that $(H_M,S_M)=M'$ is a medial map for
every map $M$. It is evident that $H_M$ is 4-regular. To prove
that $M'$ is properly 2-face-colorable, it is sufficient to show
that $(C_M,S_M)$ is 3-face-colorable; this is straightforward:
color with black the discs bounded by $v$-gons, with white the
discs bounded by $f$-gons, and with red the discs bounded by the
squares of $M$. When we contract the squares to get $M'$, the red
faces disappear and we have the ordered 2-coloring, black and
white, of the faces of $M'$. This shows that $M'$ is a medial map.

Conversely, we show that $(\psi')^{-1}$ is defined for all medial
maps. Given a medial map $X=(H,S)$, expand in $S$ the 4-valent
vertices of $H$ to discs bounded by squares. Of course, this
expansion gives the embedding of a map, denoted $M$: the original
edges of $H$ are the $a_M$-edges; the created edges which bound a
black face are the $v_M$-edges; the created edges which bound a
white face are the $f_M$-edges. It is evident that the $v$-gons,
the $f$-gons, and the squares of $M$ are boundaries of closed
discs in $S$. The lemma is proved.\fimproof

Observe that if $X,Y\in\Gamma(M)$, then $H_X=H_Y$. In the proof of
the above lemma, the black faces of $M'$ correspond to the
$v$-gons of $M$, the white faces to the $f$-gons. Henceforth, we
fix this convention.

Note that the only difference between $M'$ and $D'$ is that the
black and white faces are interchanged.

The set of black faces of $M'$ is denoted by $bF(M')$; its set of
white faces is denoted by $wF(M')$.

Given a map $M$ and a vertex $v$ of valency $2n$ in $G_M$, two
edges incident to $v$ are opposite at $v$ if their indices differ
by $n$ in a successive numbering of the edges incident to $v$ that
follows the cyclic sequence induced by $M^t$.

For an eulerian $G_M$\index{eulerian $G_M$}, a {\it smooth
path}\index{\it smooth path} in $M^t$ is a path in $G_M$ such that
every pair of consecutive edges is opposite at the
vertex\index{opposite at the vertex} between them.

Observe that for $G_M$ eulerian, the smooth paths of $M^t$ define
a partition of $eG_M$. In particular the smooth paths of
$M'=(H_M,S_M)$ define a partition of $eH_M$.

Given $M'$ we define $G'_M$ as follows: the vertices of $G'_M$ are
the smooth paths of $M'$. The edges of $G'_M$ are the vertices of
$H_M$. The two ends of an edge of $G'_M$ are the two smooth paths
(which may be the same) that use the corresponding vertex of
$H_M$.

\medskip\n{\bf (2.8) Lemma:} For every map $M,\ G_P\cong G'_M$, where $P$ is the
phial of $M$.
\medskip

\inicproof Consider $M=(H_M,S_M)$. There exists a natural
bijection, denoted $b_1$, from the squares of $P$ onto the
vertices of $H_M$: the expansion in $S_M$ of the vertices of $H_M$
gives the squares of $M$, which are identifiable with the squares
of $P$. Also there is a natural bijection from the $v$-gons of $P$
onto the smooth paths of $M'$; we denote by $b_2$ this bijection.

Take an edge $x$ of $G_P$ and let $u$ and $w$ be its ends. It
follows, from the definitions, that the ends of $b_1(x)\in eG'_M$
are $b_2(u)$ and $b_2(w)$, which are in $VG'_M$. Therefore,
$G_p\cong G'_M$, proving the lemma. \fimproof

Given a map $M$ we say that $P'$, the medial of the phial of $M$,
is a {\it smooth-path-embedding} of $G_M$, or an {\it
$sp$-embedding}\index{\it $sp$-embedding} of $G_M$. The {\it
surface of the $sp$-embedding}\index{\it surface of the
$sp$-embedding} $P'$ is $S_p$. A graph is {\it eulerian} if every
vertex has even valency.

\medskip\n{\bf (2.9) Theorem.} {\bf (Surface characterizations of map-rich
graphs):}
\begin{itemize}
\item[(a)] The following three conditions are equivalent for
connected graph $G$ with at least one edge:
\begin{itemize}
\item[(i)] $G$ is a map-rich graph. \item[(ii)] $G$ has an
$sp$-embedding in a surface $S$ of Euler characteristic
$\chi(S)=3-|VG|+\gamma$, where $\gamma$ is the dimension of the
bicycle space of $G$.(The value of $OL(S)$ is also determined; see
part (b)). \item[(iii)] There exists an $sp$-embedding
$P'=(H_P,S_P)$ of $G$ such that the cycle space of $H_P$ (over
$GF(2)$) is generated by the boundaries of the faces of $P'$ and
by the cycles induced by the smooth paths of $P'$.
\end{itemize}
\item[(b)] For every pair of phial maps, $M$ and $P$, if $P$ is
orientable, then $G_M$ is eulerian. The converse also holds if the
maps are rich.
\end{itemize}
\medskip

\noi{\bf Proof of part (a)} Initially we prove that (i) implies
(ii). Assume that $G$ is a map-rich graph and that $M$ is a rich
map satisfying $G_M=G$. Denote by $P$, as usual, the phial of $M$.
By the lemma on bicycle spaces\index{lemma on bicycle spaces},
proved below in 2.9.2,
$$Z\cap F=V\cap F\cap Z=V\cap V^\perp$$
where $V,F$, and $Z$ are relative to map $M$. We have also
$$V=F^\perp\cap Z^\perp,$$
since $M$ is rich and thus $V^\perp=F+Z$. We use the two latter
equalities in the expression for $\xi(P)$ which we obtain from
Lemma 2.9.1, proved below, with $P$ replacing $M$:
$$\begin{array}{rcl}
\xi(P)&=&\dim(Z'\cap F^\perp)-\dim(Z\cap F)\\
&=& \dim(V)-\dim(V\cap V^\perp).
\end{array}$$
Since $G$ is connected, $\xi(P)=2-\chi(S_p)$, and $\dim(V\cap
V^\perp)=\gamma$ it follows that
$$\xi(S_P)=3-|VG|+\gamma.$$
Therefore, since $S_p$ is the surface of the $sp$-embedding of
$G$, the proof that (i) implies (ii) is complete provided Lemmas
2.9.1 and 2.9.2 are established.

Now we prove that (ii) implies (i). Assume that $G$ has an
$sp$-embedding $P'$ such that $\xi(S_p)=3-|VG|+\gamma$. Since $G$
is connected the latter equality is equivalent to
$$\xi(P)-\dim(V)-\dim(V\cap V^\perp),$$
where $V$ is the coboundary space of $G_M$, and $M$ the phial of
$P$. By Lemma 2.9.1, with $P$ in the role of $M$,
$$Z\cap F\subseteq V\cap V^\perp\subseteq V\subseteq Z^\perp\cap F^\perp$$
$$\dim(Z^\perp\cap F^\perp)-\dim(V\cap F)=\xi(P).$$
The previous expression for $\xi(P)$ and the above two properties
imply
$$Z\cap F=V\cap V^\perp$$
$$V=Z^\perp\cap F^\perp.$$
Taking perpendicular to both sides of the latter equality we get
$$V^\perp=F+Z,$$
whence $M$ is rich. Since by Lemma 2.8, $G_M=G'_P$; and the latter
is $G$, it follows that $G$ is map-rich. Thus, the proof that (ii)
Implies (i) is complete provided Lemma 2.9.1 is established.

$$\includegraphics[width=148.942mm,height=73.255mm]{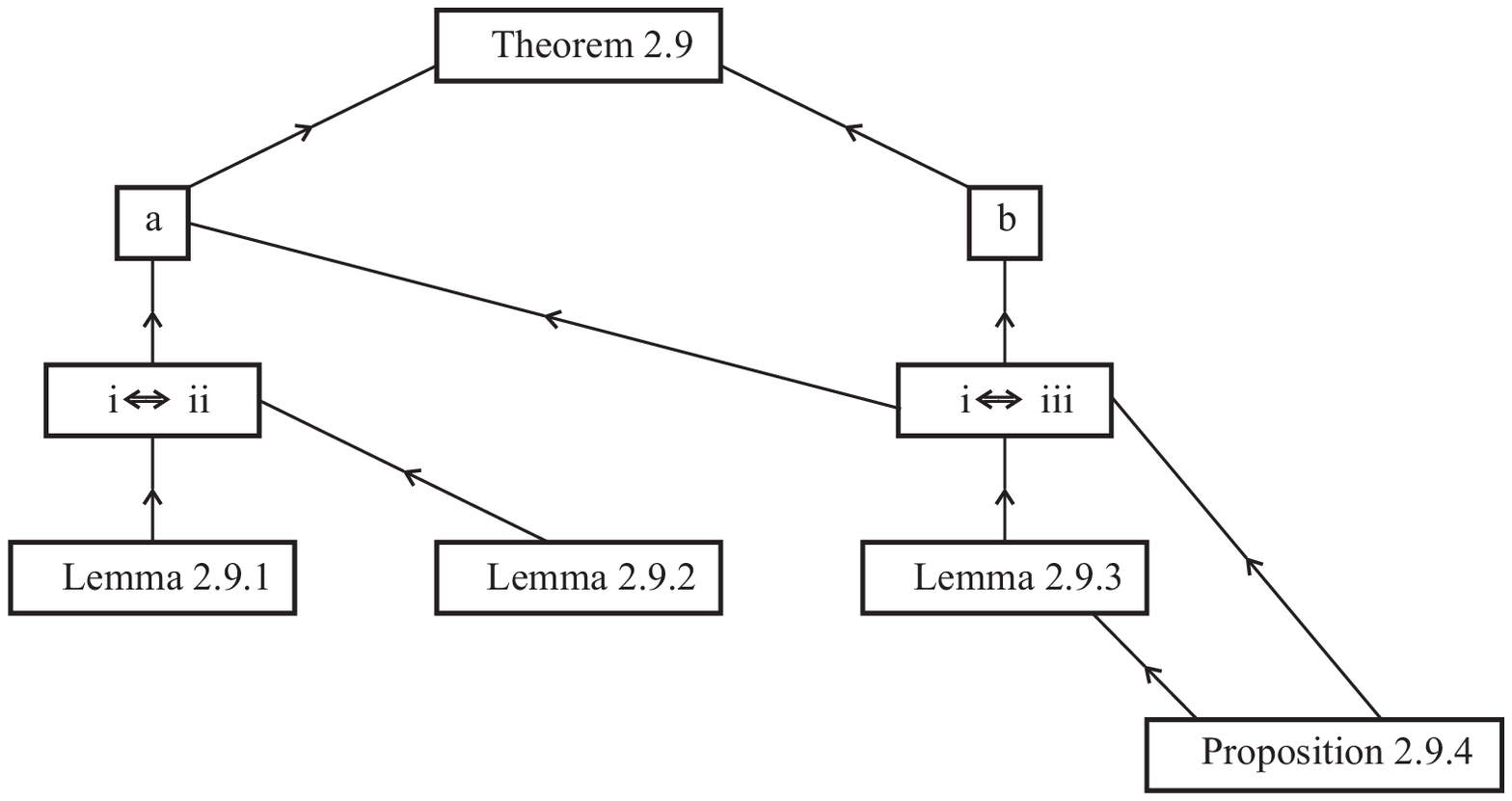}$$%[width=15cm,height=9cm]
\centerline{Dependence Structure for Theorem 2.9}

We show that (i) and (iii) are equivalent. To this end, we define
the {\it $sp$-deficiency} of a connected map\index{{\it
$sp$-deficiency} of a connected map} $P$, or ${\rm def}'(p)$, as
follows. Denote the cycle space of $H_p$ by $CS(H_p)$. Let the
subspaces of $CS(H_p)$ generated by the boundaries of the black
and white faces of $P'$ be denoted, respectively, by $Z'$ and
$F'$; denote by $V'$ the subspace of $CS(H_p)$ generated by the
cycles induced by the smooth paths of $P'$. With this terminology
we have, by definition,
$${\rm def}'(P)=\dim\frac{CS(H_p)}{V'+F'+Z'}.$$
By Lemma 2.9.3, proved below, it follows that $\dim[(F'+Z')\cap
V']=1+\dim(F\cap Z)$.

We also have $\dim[CS(H_p)]=2e-e+l=e+1$, and $\dim(F'+Z')=f+z-l$,
where $f=|vG(D)|,\ z=|vG(P)|,\ e=|eG_p|$; the first of those
equalities is implied by the fact that $P$ is connected; the
second follows from the facts that $\dim(F')=f$ and $\dim(Z')=z$,
because the boundaries of the black faces and the boundaries of
the white faces of $P'$ partition $eH_p$, and from the fact that
$\dim(F'\cap Z')=1$, proved in Proposition 2.9.4.

Since the edges in the smooth paths of $P'$ also partition $eH_p$,
it follows that $\dim(V)'=v$.

The above equalities imply that ${\rm def}(P)={\rm def}'(P)$, for
every map $P$:
$$\begin{array}{rcl}
{\rm def}'(P) &=& \dim[CS(H_p)]-\dim(V'+F'+Z')\\
& = & (e+1)-\dim(F'+Z')-\dim(V')+\dim[(F'+Z')\cap V']\\
& = & (e+1)-(f+z-1)-v+1+\dim(F\cap Z)\\
& = &  e+3+\dim(F\cap Z)-(v+f+z) \\
& = &  {\rm def}(P).
\end{array}$$

Now we can easily conclude the proof that (i) and (iii) are
equivalent. If $G$ is map-rich, let $M$ be a rich map with
$G_M=G$. By Lemma 2.0, this implies that $G'_p=G$. Since ${\rm
def}'(P)={\rm def}(P)$, as shown above, and ${\rm def}(P)={\rm
def}(M)$, by Lemma 2.4, the equality ${\rm def}(M)={\rm def}'(P)$
shows that (iii) follows from (1).

Conversely, suppose that (iii) holds for some $P'$. Let $M$ be the
phial of $P$. Since $0={\rm def}'(P)={\rm def}(P)={\rm def}(M)$,
it follows that $M$ is rich. By Lemma 2.8 we have $G_M=G'_p$;
since $G'_p=G$, it follows that $G_M=G$ and so $G$ is map-rich.

Thus, the proof of the equivalence between (i) and (iii) is
reduced to the proof of Lemma 2.9.3 and Proposition 2.9.4.

\noi{\bf Proof of part (b)} Suppose that $P$ is orientable, whence
$C_p$ is bipartite. Let $\{A,B\}$ be a bipartition of $C_p$. It
follows, from the definition of phial, that for every pair of
corners of $M, X$ and $a_M(v_M(X))$, one is in $A$ and the other
is in $B$. Thus, the number of vertices in every $v$-gon of $M$ is
a multiple of four. Since the valency of a vertex in $G_M$ is half
the number of vertices in the corresponding $v$-gon of $M$, it
follows that $G_M$ is eulerian.

To prove the converse assume that $M$ is rich and $C_M$ is
eulerian. (If we drop the hypothesis that $M$ is rich, then the
converse does not follow in general.) We use the equivalence
between (i) and (iii) established in part (a). Consider the graph
$G$ of part (a) to be $G_M$. Therefore, (iii) holds since $M$ is
rich. Let $P$ be the phial of $M$. By (iii), we have
$CS(H_p)=Z'+F'+V'$, where to recall the terminology, $Z'$ and $F'$
are generated by the boundaries of, respectively, the black and
the white faces of $P'$, and $V'$ is generated by the cycles
induced by the smooth paths of $P'$.

Since $G_M$ is eulerian it follows that each smooth path in $P'$
has an even number of edges. This claim is justified as follows.
Assume that there are $k_1$ loops and $k_2$ non-loops incident to
a vertex $v$ of $G=G_M$. The subgraph of $H_p$ induced by the
corresponding smooth path in $P'$ has $k_1$ 4-valent vertices and
$k_2$ bivalent vertices. The number of edges in this subgraph is,
therefore, $k=\frac{1}{2}(4k_1+2k_2)$. Since $G_M$ is eulerian
$k_2$ is even and so is $k$.

To each smooth path in $P'$ corresponds a pair of
$(v_p,f_p,a_p)$-paths (a cyclic path that uses $v_{p^-},f_{p^-}$,
and a $a_p$-edges in this order) in $C_p$ such that the symmetric
difference of their edge-sets is the edges in the squares that
they meet once. Observe that since each smooth path has an even
number of edges, so does each $(v_p,f_p,a_p)$-path.

We claim that the cycle apace of $C_P$ is generated by the cycles
induced by the squares, the $v$-gons, the $f$-gons, and the
$(v_p,f_p,a_p)$-paths. Observe that if the claim is accepted the
proof is complete: as every generator has an even number of
elements it follows that every cycle in $C_p$ has an even number
of edges, that is, $C_p$ is bipartite and $S_p$ is orientable. The
claim is proved by the inclusion of a bounding (in $S_P$) loop at
each vertex of $H_p$, and by the successive expansion of edges to
form $C_P$. The proof uses the following two easy
observations:\jtt
\begin{itemize}
\item[(I)] The cycle space of a graph $G$ is generated by a
collection $X$ of cycles iff the cycle space of the graph $G^e$
obtained from $G$ by adjoining a new edge $e$ which forms a loop,
is generated by the collection $X^e$ of cycles, obtained from $X$
by including $\{e\}$ as one of its elements. \item[(II)] The cycle
space of a graph $G$ is generated by a collection $X$ of cycles
iff the cycle space of the graph $G_e$, obtained from $G$ by
contracting non-loop $e$, is generated by the collection $X_e$ of
cycles, obtained from $X$ by taking the difference of its elements
with $\{e\}$.
\end{itemize}

It is evident that (I) and (II) hold. Iterated application of them
transforms $(H_p,X)$ into $(C_p,Y)$, where $X$ is the set of
cycles formed by the boundaries of the faces of $P'$ and by the
edges in the smooth paths of $P'$; $X$ generates $CS(H_p)$; $Y$ is
the set of cycles induced by the squares, the $v$-gons, the
$f$-gons of $P$, and by one member of each pair of associated
$(v_p,f_p,a_p)$-paths; $Y$ generates $CS(C_p)$. The proof of
Theorem 2.9 is complete provided statements 2.9.1 to 2.9.4 are
established. \fimproof

\medskip\n{\bf (2.9.1) Lemma.} {\bf (Triple inclusion property):} For every map $M$ we have
$$V\cap F\subseteq Z\cap Z^\perp\subseteq Z\subseteq V^\perp\cap F^\perp$$
and
$$\dim(V^\perp\cap F^\perp)-\dim(V\cap F)=\xi(M).$$
\medskip

\inicproof We restrict the proof, without loss of generality, to
connected maps. To prove the equality, consider the following
equation, with $e=|SQ(M)|,\ v=|vG(M)|,\ f=|fG(M)|$:
$$\dim(V^\perp\cap F^\perp)=(e-v+1)+(e-f+1)-\dim(V^\perp+F^\perp).$$
Using items (d) and (e) of 2.0.1 we have
$$\dim(V^\perp+F^\perp)=e-\dim [(V^\perp+F^\perp)^\perp]=e-\dim(V\cap F).$$
Using the above value of $\dim(V^\perp+F^\perp)$ yields
$$\dim(V^\perp\cap F^\perp)=(e+2)-(v+f)+\dim(V\cap F)=\xi(M)+\dim(V\cap
F),$$ which is the equality that we wish to establish.

To verify the inclusions we use the Absorption Property, Lemma
2.5, which says that
$$V\cap F\subseteq Z.$$
Proposition 2.2 implies that $Z^\perp\supseteq V+F$. Taking
perpendiculars we get, by 2.0.1(c)
$$Z\subseteq V^\perp\cap F^\perp.$$
Thus, we have $V\cap F\subseteq Z\subseteq V^\perp \cap F^\perp$.
To conclude the proof it is enough to show that $V\cap F\subseteq
Z^\perp$. This follows from $V\cap F\subseteq V+F\subseteq
Z^\perp$. The latter inclusion is implied by Proposition 2.2. The
proof of the lemma is complete. \fimproof

\medskip\n{\bf (2.9.2) Lemma on bicycle spaces:} If $M$ is a rich map, then
the bicycle spaces of $G_M,G_D$, and $G_p$ are equal to $V\cap
F\cap Z$.
\medskip

\inicproof Observe that this lemma is a generalization of part (b)
of the corollary of Proposition 2.6. Therefore, here is a property
of planar graphs which generalizes to rich graphs.

The fact that $M$ is rich implies that $V=F^\perp\cap Z^\perp$.
Hence, $V\cap V^\perp=F^\perp\cap Z^\perp\cap V^\perp$. The
symmetry in $V,F,Z$, shows that $F\cap F^\perp$ and $Z\cap
Z^\perp$ are also equal to $V^\perp\cap F^\perp\cap Z^\perp$. We
have $V\cap F\cap Z=V\cap F$, which follows from the Absorption
Property (Lemma 2.5). Also we have $V\cap F\subseteq V\cap
V^\perp$ by property 2.2. Since $V\cap V^\perp=V^\perp\cap
F^\perp\cap Z^\perp$, it follows that $V\cap F\cap Z\subseteq
V^\perp\cap F^\perp\cap Z^\perp$. To prove the reverse inclusion
let $x\in V^\perp\cap F^\perp\cap Z^\perp$. As shown above the
latter space is equal to $V\cap V^\perp=F\cap F^\perp=Z\cap
Z^\perp$. Therefore, $x\in V\cap F\cap Z$ and $V^\perp\cap
F^\perp\cap Z^\perp\subseteq V\cap F\cap Z$, concluding the proof.
\fimproof

\medskip\n{\bf (2.9.3) Lemma:} For every connected map $P=P_M$,
$$\dim[(F'+Z')\cap V']=1+\dim(F\cap Z).$$
\medskip

\inicproof We intend to prove the above equality by presenting an
epimorphism $h$ from $(F'+Z')\cap V'$ onto $F\cap Z$ and by
showing that the kernel of $h$ has dimension 1.

By Proposition 2.9.4, proved below, $F'\cap Z'$ has dimension 1
and hence two elements. It is also shown that the two elements of
$F'\cap Z'$ are the null set of edges and the whole set of edges
of $H_p$. Let $x\in F'+Z'$.  By the preceding argument, there are
two ways to represent $x$ as $\mod\ 2$ sum of boundaries of white
and black faces of $P'$. Let $W_x$ and $B_x$ be subsets of
$wF(P')$ and $bF(P')$ such that
$$x=\partial(W_x)+\partial(B_x).$$
Since the non-null element of $F'\cap Z'$ is $eH_p$, and
$eH_p=\partial(wF(P'))=\partial(bF(P'))$, we have
$$\begin{array}{rcl}
x&=&[\partial(w_x)+\partial(wF(P'))]+[\partial(B_x)+\partial(bF(P'))]\\
[5pt] &=&\partial(wF(P')\backslash W_x)+\partial(bF(P'))\backslash
B_x),
\end{array}$$
which is the other way to represent $x$ as $\mod\ 2$ sum of
boundaries of white and black faces of $P'$.

Note that if $x\in(F'+Z')\cap V'$, then for an arbitrary vertex
$v$ of $H_p$, the intersection of the edges incident to $v$ with
$x$ consists of either zero, one, or two pairs of opposite edges.
Moreover, this intersection is one pair if and only if $v$ (which
is an edge of $G_D$ and $G_p$) is in the coboundaries of
$W_x\subseteq VG_D$, and $B_x\subseteq VG_p$, in the adequate
graphs. Therefore, $\delta_{G_D}(W_x)=\delta_{G_p}(B_x)$.

The function $h$ from $(F'+Z')\cap V'$ into $F\cap Z$ is defined
as follows. Given $x\in(F'+Z')\cap V'$, represent $x$ as
$\partial(W_x)+\partial(B_x)$ for convenient subsets of white and
black faces. Take $h(x)$ to be the coboundary of $W_x$ in graph
$G_D$, which is equal, as explained above, to the coboundary of
$B_x$ in graph $G_p$. Since the coboundaries of a subset and its
complement are the same, it follows that $h(x)$ is independent of
the representation of $x$. Thus, $h$ is well defined.

Now we show that $h$ is a homomorphism. Take $x$ and $y$ in
$(F'+Z')\cap V'$. Find $W_x,B_x,W_y,B_y$ such that
$$x=\partial(W_x)+\partial(B_x),$$
$$y=\partial(W_y)+\partial(B_y).$$
Adding $\mod\ 2$ we get
$$x+y=\partial(W_x+W_y)+\partial(B_x+B_y).$$
Thus,
$$h(x+y)=\partial(W_x+W_y)=\delta(W_x)+\delta(W_y)=h(x)+h(y).$$
Hence, $h$ is a homomorphism. To prove that $h$ is surjective take
an element $u\in F\cap Z$. By definition, there exist
$W_u\subseteq VG_D$ and $B_u\subseteq VG_p$ such that
$$u=\delta(W_u)=\delta(B_u).$$
Consider $W_u$ and $B_u$ as the subsets of the corresponding white
and black faces of $P'$. Define $x\in F'+Z'$ as
$x=\partial(W_u)+\partial(B_u)$. It can be directly observed that
the intersection of $x$ with the four edges incident to an
arbitrary vertex consists of zero, one, or two pairs of opposite
edges. Therefore, $x\in(F'+Z')\cap V'$. Moreover, $h(x)$ is the
coboundary of $W_u$ in $G_D$ that is, $h(x)=u$. Thus, $h$ is an
epimorphism.

Now we prove that the kernel of $h$ has dimension 1. Suppose that
$h(x)=\phi$ and $x$ is not the null vector. Find $W_x$ and $B_x$
such that $x=\partial(W_x)+\partial(B_x)$. One of $W_x$ and $B_x$
is non-empty, for instance, $B_x\ne\phi$. If some black face $b$
is in $B_x$, then all the black faces adjacent to $b$ by a vertex
of $H_p$ are also in $B_x$: where this not true, then
$h(x)=\partial(B_x)\ne\phi$, contradicting the hypothesis. By the
connectivity of $P$, it follows that $B_x=bF(P')$. The same
argument shows that if $W_x\ne\phi$, then $W_x=wF(P')$. But if
$B_x=bF(P')$ and $W_x=wF(P')$, then $x$ is empty, contradicting
the hypothesis. Thus, if $x\ne\phi$ and $h(x)=\phi$, then
$x=\partial(bF(P'))=\partial(wF(P'))=eH_p$. Hence, the kernel of
$h$ has dimension 1. By the fundamental homomorphism theorem, it
follows that
$$\dim[(F'+Z')\cap V']=\dim[{\rm Ker}(h)]+\dim[im(h)]$$
or
$$\dim[(F'+Z')\cap V']=1+\dim(F\cap Z),$$
proving the lemma, subject to the establishment of Proposition
2.9.4, which we do next. \fimproof

\medskip\n{\bf (2.9.4) Proposition:} For connected $P=P_M$ the dimension of $F'\cap Z'$ is 1.
\medskip

\inicproof Assume that $x\in F'\cap Z'$. lt follows that there are
$W_x\subseteq wF(P')$ and $B_x\subseteq bF(P')$ such that
$\partial(W_x=x=\partial(B_x)$. If a face $b\in B_x$, then the
elements of $wF(P')$ which are adjacent to $b$ by an edge of $H_p$
are in $W_x$. Therefore, by the connectivity of $P$, it follows
that if both of $W_x$ and $B_x$ are not empty, then $W_x=wF(P')$
and $B_x=bF(P')$. Hence, if $x\in F'\cap n Z'$ then $x$ is either
the null set or the whole set of edges, $eH_p$. Thus, $\dim(F'\cap
Z')=1$ as stated. \fimproof

Observe that at this point every lemma and proposition used to
prove Theorem 2.9 is established. To discover the surface with
maximum Euler characteristic for which a (usual) embedding of a
given graph exists is, in general, a difficult problem.  In
attacking the famous Heawood Conjecture, Ringel, Youngs, and
others have determined the surface of maximum Euler characteristic
where a (usual) embedding of a complete graph exists. See Ringel
[18], Youngs [29]. With our next theorem we derive a natural upper
hound on the Euler characteristic of a surface where there exists
an $sp$-embedding of a given graph. This upper bound holds as an
equality if and only if the graph is map-rich. In the next section
we show that all complete graphs are map-rich.

\medskip\n{\bf (2.10) Theorem:} Given a connected graph $G$ and an $sp$-embedding
of $G,P'=(H_p,S_p)$, then we have
$$\chi(S_p)\leq 3-v+\gamma,$$
where $\gamma$ is the dimension of the bicycle space of $G$, and
$v=|VG|$. Moreover, equality holds iff $P$ is rich.
\medskip

\inicproof The proof of this theorem is implicit in Theorem 2.9.
However, we give a direct proof. Let $M$ be the phial of $P$.
Hence, $G$ is $G_M$, by Lemma 2.8, and so $G$ inherits the usual
terminology for $M$. Since $G$ is connected, we have
$$\begin{array}{rcl}
{\rm def}(P)&=&\dim\,Z-\dim(V+F)\\ [5pt]
&=&(e-z+l)-(v-l)-(f-l)+\dim(V\cap F)\\ [5pt]
&=&(e+3)-(v+f+z)+\dim(v\cap F).
\end{array}$$
Observe that $V\cap F\subseteq V\cap V^\perp$, by Proposition 2.2.
Since $V\cap V^\perp$ is the bicycle space of $G$ we have
$$0\leq{\rm def}(P)\leq(e+2)-(f+z)-(v-1)+\gamma,$$
which implies $\chi(S_p)\leq 3-v+\gamma$, proving the inequality.
if $P$ is rich, then by the lemma on bicycle spaces, 2.9.2, and by
the Absorption Property, 2.5, we have
$$V\cap V^\perp=V\cap F\cap Z=V\cap F.$$
Moreover, the value of ${\rm def}(P)$ is 0. Thus, equality holds.

Assume that the equality holds, $\chi(S_p)=3-v+\gamma$. Then we
have
$$f+z-e=3-v+\gamma,$$
or
$$0=(e+3)-(v+f+z)+\gamma.$$
Since $\gamma\leq\dim(V\cap F)$, by Proposition 2.2, it follows
that $0\leq{\rm def}(P)$. But the definition of a map is at least
zero. Hence, ${\rm def}(P)=0$, i.e., $P$ is rich. \fimproof

\section*{2.D\quad Map-Rich and Non-Map-Rich Graphs}
We use Theorem 2.9 to show that all complete graphs are map-rich,
and to show that there are graphs which are not map-rich.

For complete graphs with an even number of vertices, $K_{2n}$, it
is an easy matter to construct $sp$-embeddings in the real
projective plane. To see this choose $4n$ points in the boundary
of a disc, and link by chords opposite points, such that any two
chords meet exactly once; if more than two chords meet at the same
point, effect a local deformation of some of them, to make every
crossing point the intersection of two chords. Now, identify every
point in the boundary of the disc with its opposite, or antipodal
point. With this antipodal identification the system of chords is
now a $t$-map, with the real projective plane as its surface. The
crossing points are the vertices of the graph, and its edges are
the segments of chords between two crossings\index{crossings}. By
construction, this graph is 4-regular. We show that the $t$-map is
properly 2-face-colorable. This is evident if $n=1$. Assume that
it is true for $n$ at most $k$, and let $n$ be $k+1$. Observe that
every pair of chords is the boundary of a subset $R$ of faces.
Remove two chorda, $c_1,c_2$, and obtain, by the inductive
hypothesis, a proper 2-face-coloring of the remaining system,
which has $2k$ chords. Put back the two chords. Observe that all
the edges in the boundary of $R$, the set of faces bounded by
$c_1$, and $c_2$, have the same color on both sides. Edges not in
the boundary of $R$ have different colors in each side. Therefore,
to get a proper 2-coloring just interchange the colors in the
faces of $R$. This completes the induction. Therefore, the
projective $t$-map becomes after ordering the bicoloration, a
medial map. A medial map that arises by the above construction is
called a {\it simple $2n$-system} of {\it projective lines}.
Observe that if $P'$ is such an object, then $G_M=G'_p=K_{2n}$.

\medskip\n{\bf (2.11) Theorem:} For every $n$, the complete graph $K_{2n}$ is map-rich.
\medskip

\inicproof To prove that $K_{2n}$ is map-rich it is sufficient to
show that if $P'=(H_p,S_p)$ is a simple $2n$-system of projective
lines, then $M$ is rich: since in this case $G_M-K_{2n}$, it
follows that $K_{2n}$ is map-rich. By Theorem 2.9 (iii) it is
enough to show that the cycle space of $H_p$ is generated by the
boundaries of the faces of $P'$ and the cycles induced by the
smooth paths. Since $S_p$ is the projective plane, the facial
deficiency of $P$ is 1, by Proposition 2.6. Therefore, we just
need one cycle induced by a smooth path not spanned by the
boundaries of the faces of $P'$. Observe that each one of these
cycles is non-bounding, and therefore, none of them is in the span
of the boundaries of the faces. This concludes the proof that
every $K_{2n}$ is map-rich.\fimproof

We now want to extend the previous theorem to complete graphs with
an odd number of vertices.

\medskip\n{\bf (2.12) Theorem:} Assume that $M$ is rich and $G_M=K_{2n}$.
Then deleting in $M^t$ all the edges incident to a vertex of
$K_{2n}$, we have a rich map $N$ with $G_N=K_{2n-1}$.
\medskip

\inicproof We use the labels $\{1,2,\dots,2n\}$ for the vertices
of $G_M=K_{2n}$. Let $[i,j]$ denote the edge with ends $i$ and
$j$. We delete $[2,2n],[3,2n],\dots,[2n-1,2n]$, in this order,
obtaining, respectively, $M_2M_3,\dots,M_{2n-1}$. Let $F_i$ and
$Z_i$ be, respectively, the coboundary spaces of $G_{D_i}$ and
$G_{P_i}$ for $i=2,\dots,2n-1`$. By Proposítion 2.12.1, proved
below, $[i,2n]$ is an edge of an element of $F_{i-i}\cap Z_{i-1}$.
(We take $M_1$ to be $M$.) By Proposition 2.12.2, proved below, it
follows that each one of the $M_i's$, are rich. Observe that
$G_{M_{2n-1}}$ is $K_{2n-1}$ plus a pendent edge $[1,2h]$: at this
point vertex $2n$ is monovalent. Proposition 2.12.3, given below,
shows that the deletion in $M^t$ of a pendent edge does not alter
the deficiency of the map. Therefore, $N$ defined from
$M^t_{2n-1}$ by the deletion of $[1,2n]$ is a rich map with
$G_N=K_{2n-l}$. The proof of the theorem is thus reduced to the
proofs of Propositions 2.12.1, 2.12.2, and 2.12.3. \fimproof

The terminology involved in the statement of the following
proposition is defined in the proof of Theorem 2.12.

\medskip\n{\bf (2.12.1) Proposition:} For $i=2,3,\dots,2n-1$, edge $[i,2n]$ is an
element of a coboundary in $F_{i-1}\cap Z_{i-1}$.
\medskip

\inicproof Since $M_1=M$ is rich, it follows by Lemma 2.9.3 that
$F_1\cap Z_1$ is equal to the bicycle space of $K_{2n}$. Consider
the following sequence, $B_1$, of $2n-2$ coboundaries in the
bicycle space of $K_{2n}$, and hence in $F_1\cap Z_1$:
$$B_1=(\delta(\{1,2\}),\ \delta(\{1,3\}),\dots,\delta(\{1,2n-1\})).$$
We also define $B_j$, for $2\leq j\leq 2n-2$, by deleting from
$B_1$ its first $J-1$ elements.

Observe that the only element in $B_1$ that contains $[2,2n]$ is
$\delta(\{1,2\})$. The effect of deleting an edge in $m^t$ which
is not a loop in $G_D$ nor in $G_p$, as in our case, is the
contraction of this edge in both graphs. Therefore, the elements
of $F_2\cap Z_2$ are the elements of $F_1\cap Z_1$ which do not
contain edge $[2,2n]$. Hence, every element in $B_2$ is an element
of $F_2\cap Z_2$.

The above argument iterates as follows: the only coboundary in
$B_j$ which contains edge $[j+1,2n]$ is its first coboundary,
$\delta(\{1,j+1\})$. Therefore, the elements of $B_{j+1}$ are also
elements of $F_{j+1}\cap Z_{j+1}$. The proposition follows because
for $i=2,3,\dots,2n-1$ edge $[i,2n]$ is an element of coboundary
$\delta(\{1,i\})$, which is a member of $F_{i-1}\cap Z_{i-1}$.
\fimproof

\medskip\n{\bf (2.12.2) Proposition:} Assume that an edge $e$ of $m^t$ belongs to an element
of $F\cap Z$. Then the deletion of $e$ does hot change the
deficiency of $M$.
\medskip

\inicproof Let $M^t$ denote the $t$-map obtained from $M^t$ by the
deletion of $M$. We want to prove that
$$\dim[V^\perp/(F+Z)]=\dim[V^\perp_e/(F_e+Z_e)].$$
To accomplish this, it is enough to show that

(a) $\dim(V^\perp)=1+\dim(V^\perp_e)$;

(b) $\dim(F+Z)=1+\dim(F_e+Z_e)$.

\noi We claim that $\{e\}$ is not a coboundary in $G_M$: if
$\{e\}$ were a coboundary in $G_M$, then $e$ would be a loop in
$D^t$, and so $\{e\}$ would not be in a member of $F\cap Z$,
contradicting the hypothesis. Hence, we have $\dim(V^\perp)=l+\dim
V^\perp_e$, establishing (a).

To prove (b), note that $e$ is a non-loop in both graphs, $G_D$
and $G_P$. The difference between these graphs and
$G_{D_e},G_{P_e}$ is that in the latter two, $e$ is contracted.
Therefore, $\dim(F)=\dim(F_e)+1$, and $\dim(Z)=\dim(Z_e)+1$. To
conclude the proof of the proposition, it suffices to show that
$\dim(F\cap Z)=\dim(F_e\cap Z_e)+1$.

The elements of $F_e\cap Z_e$ are the coboundaries of $F\cap z$
which do not contain edge $e$. Take an arbitrary basis of $F\cap
Z$:
$$\{A_1,\dots,A_k;B_1,\dots,B_\ell\}$$
where the $A_i's$ do not contain $e$, and the $B_i's$ contain $e$.
Observe that $\ell\leq 1$. Replace this basis by
$$\{A_1,\dots,A_k;B_1,B_1+B_2\dots,B_1+B_\ell\}.$$
In the new basis the only member that contains $e$ is $B_1$.
Therefore, if we remove $B_1$ from the above set, we have a basis
for $F_e\cap Z_e$ proving that $\dim(F\cap Z)=1+\dim(F_e\cap
Z_e)$. The proof of the proposition is complete. \fimproof

\medskip\n{\bf (2.12.3) Proposition:} If $e$ is a pendent edge of $G_M$ then the deletion of
$e$ in $m^t$ does not change ${\rm def}(M)$.
\medskip

\inicproof Denote by $M^t$ the $t$-map obtained from $M^t$ by
deleting $e$. Edge $e$ is a loop in both $G_D$ and $G_p$.
Therefore, $F+Z$ is equal to $F_e+Z_e$, where $F_e$ and $Z_e$ are
the coboundary spaces of $D_e$ and $P_e$, dual and phial of $M_e$.
Also, since $\{e\}$ is a coboundary in $G_M,\ V^\perp=V^\perp_e$,
concluding the proof. \fimproof

We observe that with the proof of the latter proposition we have
complete the proof of Theorem 2.12. Thus, all complete graphs are
map-rich.

If we let the torus take the role of the projective plane, a
parallel development shows that all complete bipartite graphs are
map-rich. Also using the torus we can show that the complete
tripartite graphs of the form $K_{n,n,n}$ are map-rich. Those are
straightforward applications of Theorem 2.9 and the techniques
developed in the case of the complete graphs.

With the help of Theorem 2.9 it is an easy matter to get examples
of rich maps in arbitrary surfaces: draw a system of lines which
generates the homology $(\mod\ 2)$ of the surface: draw enough
more lines to have a 2-face-colorable $t$-map; deform locally some
lines incident to a multiple crossing to obtain a 4-regular graph.
Indeed, by the equivalence between (i) and (iii) in part (a) of
theorem 2.9, every rich-map arises in this way.

Now we use the equivalence between (i) and (ii) in Theorem 2.9(a)
to show that not all graphs are map-rich. For a graph $G$ denote
by $G^d$ the graph obtained from $G$ by doubling each one of its
edges.

\medskip\n{\bf (2.13) Lemma:} For any given connected graph $G,G^d$ is
map-rich iff it has an $sp$-embedding in the plane.
\medskip

\inicproof Observe that the coboundary space $V$ of $G^d$ is
contained in its cycle space $V^\perp$. Therefore, $\dim(V\cap
V^perp)=\dim(V)=V-1$, where $v=|VG^d|$. It follows from Theorem
2.9 (ii) that there exists a rich map $P$ with $G'_P=G^d$ iff
$\chi(S_p)=(3-v)+(v-1)$, or $\chi(S_p)=2$. Thus, $G^d$ is map-rich
iff it has an $sp$-embedding in the plane, proving the lemma.
\fimproof

Using the above lemma we can prove easily by exhaustion that if
$G$ is the Petersen graph, then $G^d$ is not map-rich. This
exhaustion is made easy by the symmetry of the Petersen graph and
by the fact that $C^d$ where $C$ is a polygon without chords, has
a unique $sp$-embedding in the plane.

We conclude this chapter with the establishment of an upper bound
for the maximum number of independent vertices of a map-rich graph
$G$ which has all its vertices odd valent.

\medskip\n{\bf (2.14) Theorem:} Assume that all the vertices of a map-rich
graph $G$, connected and without loops, are odd valent. Denote by
$\gamma$ the dimension of the bicycle space of $G$ and let
$v=|VG|$. The maximum number of independent vertices in $G$,
$\alpha$, is at most equal to $v-l-\gamma$. In particular, if $G$
is planar, then $\alpha\leq v-z$, where $z$ is the number of
zigzag paths in a planar realization of $G$.
\medskip

\inicproof Take a rich map $M$ such that $G_M=G$. As usual, denote
by $P$ the phial of $M$; bence, $P'=(H_p,S_p)$ is an
$sp$-embedding of $G$. By Theorem 2.10 $\xi(S_p)=v-l-\gamma$. Take
an independent set of vertices in $VG$, denoted I, of cardinality
$\alpha$. To the vertices in I correspond a collection $J'$ of
disjoint smooth paths of odd length in $P'$. The collection $J'$
corresponds to a collection $J$ of disjoint odd polygons in $C_p$.
(Think in terms of the expansion of the vertices of $H_p$ to form
the discs bounded by the squares of $P$). From Theorem 1.8, it
follows that the $\mod\ 2$ sum of the edges in an arbitrary set of
disjoint odd polygons in $C_p$ is not equal to the $\mod\ 2$ sum
of the edges in some alternately bicolored polygons of $C_p$.
Thus, the members of $J$ are representatives of distinct classes
in the quotient space
$$\frac{CS(C_p)}{\aw SQ(P)\fw+\aw VG(P)\fw+\aw fG(P)\fw},$$
which by Proposition 2.6 (applied to the $t$-map $(C_p,S_p)$) has
$\xi(P)$ classes. Hence,
$$\alpha=|J|\leq\xi(P)=v-1-\gamma.$$
If $G$ is planar, and $M$ is a planar map with $G_M=G$ we
specialize as follows. By the corollary of 2.6, $M$ is rich and
the bicycle space of $G$ is equal to the coboundary space of
$G_P$, where $P$ is the phial of $M$. Thus $\gamma=z-1$ and
applying the bound obtained we get $\alpha\leq v-z$. The proof is
complete. \fimproof

As two examples of planar graphs where the upper bound of Theorem
2.14 holds as equality we mention the cube and the tetrahedron:
the cube has eight vertices, four zigzags (in the planar drawing)
and $\alpha$ equal to 4; the tetrahedron has four vertices, three
zigzags and $\alpha$ equal to 1.

\chapter{Gauss Code Problem Beyond the Plane}\vspace{-1.5cm}
\section*{3.A\quad Preliminary Concepts and Relation with Map-Rich Graphs}\vspace{-1pc}
The Gauss code problem is the following: given a cyclic sequence
of symbols in which each symbol occurs twice, find necessary and
sufficient conditions to embed a closed curve in the plane with
the restriction that the cyclic sequence of self-intersections of
the curve reproduces the given cyclic sequence, when the symbols
are identified with the self-intersection points. Consider, for
instance, the cyclic sequence
$(1,3,4,6,7,2,3,9,6,5,2,1,8,7,5,4,9,8).$
$$\begin{array}{c}
\includegraphics[width=7cm,height=7cm]{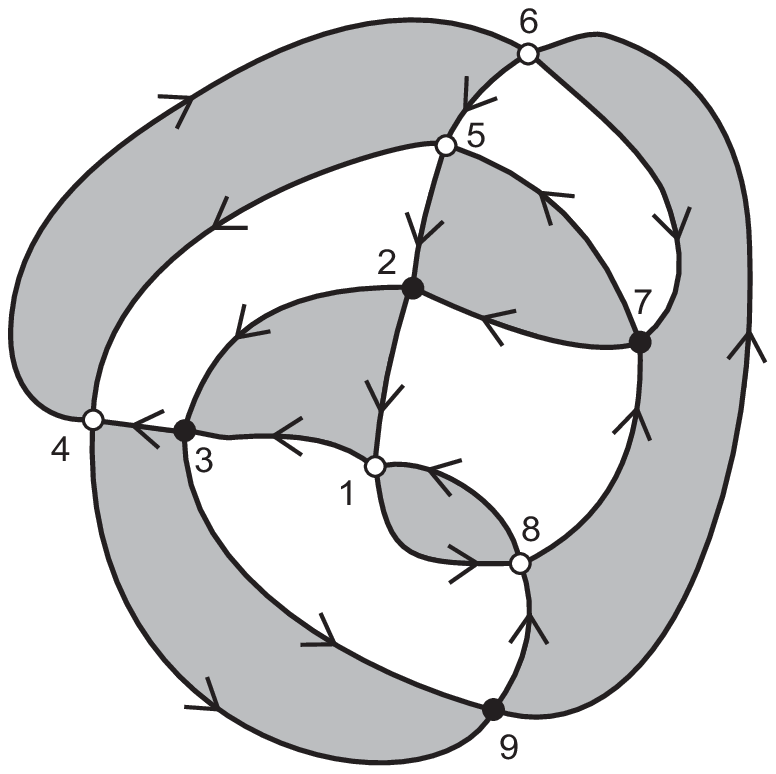}\\
\centerline{Figure 3.1}
\end{array}$$

A planar ``realization'' of this sequence is given in Figure 3.1.
Topological as well as combinatorial solutions for this problem
are known. See a brief history in Gr\"unbaum [10]. Combinatorial
solutions are provided by Lovasz and Marx [17] and by Rosenstiehl
[19]. In this chapter we consider a variant of this problem which
replaces the plane for an arbitrary surface, but introduces the
restriction that the $t$-map formed by the embedding of the curve
is 2-face-colorable. For the plane the 2-face-colorability is
automatically satisfied since the dual of a 4-regular planar graph
is bipartite. We present a necessary and sufficient algebraic
condition (Corollary of (3.3.a),(3.4), and (3.5.a)) for the
realizability of a given sequence in a given surface, in which the
restriction of 2-face-colorability is imposed. The conditions
provide a good combinatorial algorithm for the realizability
problem in the cases of surfaces with connectivity at most one:
plane and projective plane.

In this section we introduce a few more concepts arising from maps
and we show how the solution of the extended version of the Gauss
code problem relates to map-rich graphs\index{map-rich graphs}. A
\index{\it balancing partition} for a map $M$ is a partition of
its corners into two classes satisfying the following
conditions:\jtt
\begin{itemize}
\item[(i)] For every corner $x$ of $M$, $x$ and $v_M(x)$ are in
different classes. \item[(ii)] For every corner $x$ of $M$, $x$
and $a_M(V_M(x))$ are in the same class.
\end{itemize}

If $\{A,B\}$ is a balancing partition for map $M$, then the {\it
$v$-ordering induced by $A$} (by $B$) is the collection of cyclic
sequences formed by the squares containing the corners in the
cycles of permutation $a_M\circ V_M$ acting on $A$ (on $B$).

Observe that the elements in the $v$-ordering induced by $A$ are
the inverse cyclic sequences of the elements in the $v$-ordering
induced by $B$.

A square of $M$ is called \index{\it balanced} relative to a given
balancing partition if opposite corners of the square belong to
the same class of the balancing partition; otherwise the square is
called {\it unbalanced}.

A set of squares of $M$ is called an {\it imbalance} in $M$ if it
is the set of unbalanced squares relative to some balancing
partition.

We note that if $\{A,B\}$ is a balancing partition for $M$, then
interchanging between $A$ and $B$ the corners incident to a
$v$-gon of $M$ produces a balancing partition $(A_1,B_1)$.
Moreover, it is evident that given two balancing partitions for
$M$ we can go from one to the other by the above interchanges.
This observation implies that the set of imbalances of a map $M$
corresponds, under the identification of $SQ(M)$ with $eG_M$ to a
coset of the coboundary space of $G_M$. It implies also that $M$
is orientable iff the null set of edges is an imbalance of $M$; in
the latter case (and only in it) the set of imbalances in $M$ is
the coboundary space of $G_M$.

A {\it descriptor} for a map\index{{\it descriptor} for a map} $M$
is a $v$-ordering relative to some $A$, where $\{A,B\}$ is a
balancing partition for $M$, together with the imbalance relative
to $\{A,B\}$.

The following proposition justifies the terminology
``descriptor''.

\medskip\n{\bf (3.0) Proposition:} Every map is reconstructible from a
descriptor for it.
\medskip

\inicproof To define a map $M$ we simply have to specify how the
corners of its squares are paired by its $a_M$-edges. Use the
$v$-ordering of the descriptor to pair consistently (with the
$v$-ordering) the corners of the squares. Clearly, there are many
different ways to effect a consistent pairing. Choose an arbitrary
such pairing and call the resulting map $N$. Choose a balancing
partition for $N$ and denote by $I$ the induced imbalance.
Denoting by $J$ the imbalance of the given descriptor, what we
have to do to transform $N$ into $M$ is to interchange the
incidences of the pairing edges at the ends of one $v_M$-edge for
each square in the symmetric difference of $I$ and $J$. All those
interchanges are necessary to satisfy the given imbalance. They
are also sufficient because the descriptor constructed is a
descriptor for the resulting map. \fimproof

Edmonds' embedding technique [4], in our language, describes an
orientable map as a $v$-ordering and implicitly lets its imbalance
be the null set.

In the remaining of this section we show how the problem of
deciding whether a given graph with maximum valency 3 is map-rich
is reduced to the problem of deciding what is the minimum
connectivity of a surface where a realization of a Gauss
code\index{realization of a Gauss code} satisfying the
2-face-colorability exists.

The material that follows in the rest of this section is not a
pre-requisite for Section 3.

Given a connected graph $G$, the algorithmic problem of deciding
whether $G$ in map-rich seems to be difficult. A partial answer is
given by Theorem 2.9: $G$ is map-rich iff it admits an
$sp$-embedding in the surface $S$ of Euler characteristic
$\chi(S)=3-|VG|+\gamma$, where $\gamma$ is the dimension of the
bicycle space of $G$, and orientability label $OL(S)=1$ or $-1$
according to whether or not $G$ is eulerian. Observe that $\gamma$
can be computed by an algorithm, polynomial in the size of $G$:
denote by $V$ the coboundary space (over $GF(2)$) of $G$. Whence,
$\gamma=\dim(V\cap V^\perp)=|eG|-\dim(V+V^\perp)$; the latter
equality follows from 2.0.1 (c) and (d). The computation of
$\dim(V+V^\perp)$ reduces to the problem of finding the rank of a
matrix whose rows are formed by the characteristic vectors of
arbitrary bases for $V$ and $V^\perp$.

However, knowing the surface where $G$ must be $sp$-embedded to
establish that $G$ is map-rich still leaves open the problem of
finding the $sp$-embedding.

The situation is analogous to graphs that have an abstract dual.
Given a graph $G$, the easiest way to establish that $G$ has an
abstract dual is, by Whitney's Theorem [26], to construct a
(usual) embedding in the plane. Since there are good algorithms to
test planarity, the problem of deciding whether a graph $G$ has an
abstract dual is well-solved.

Part of the problem of testing $sp$-embeddability of a given graph
in a specific surface lies in the fact that we do not know which
cyclic sequence the points of intersection of a smooth path (that
in the original graph corresponds to a vertex) will define.

In view of the above we were led to consider the problem of
deciding whether a given graph has an $sp$-embedding in a specific
surface, with the restriction that the inverse pairs of cyclic
sequences of intersection points of the smooth paths are fixed in
advance. For graphs of maximum valency 3 this problem is
equivalent to the unrestricted one, since these graphs have just
one possibility for the inverse pairs of cyclic sequences.

Now we show how the restricted problem is algorithmically reduced
to its particular case in which the graph has just one vertex, or
equivalently, in which there is just one smooth path forming the
$sp$-embedding. Assume that a connected graph $G$ is $sp$-embedded
in a surface $S$. Consider that $G$ has more than one vertex. Take
an edge $e$ of $G$ which is not a loop. Edge $e$ is the
intersection point of two distinct smooth paths in the
$sp$-embedding. There are two ways of replacing $e$ by two points,
$e'$ and $e''$, (which become the vertices of a bounding digon in
$S$), such as to merge the two smooth paths into one. See Figure
3.2.\medskip
$$\includegraphics[width=2cm,height=2.5cm]{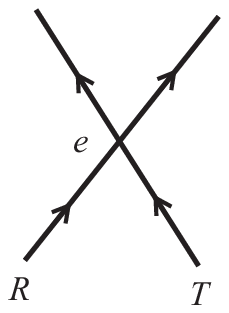}
\qquad\raise 1.15cm\hbox{$\Rightarrow$}\qquad
\includegraphics[width=2cm,height=2.5cm]{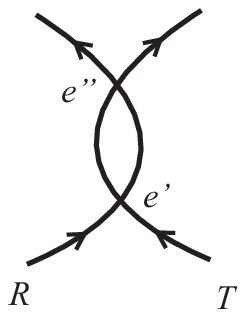}\qquad
\qquad\raise 1.15cm\hbox{or}\qquad
\includegraphics[width=2cm,height=2.5cm]{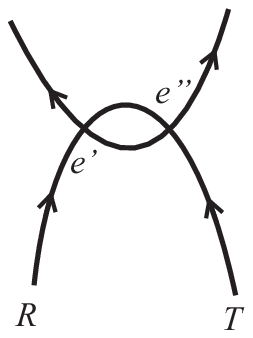}$$
\centerline{Figure 3.2}

Assume that $R$ and $T$ are subsequences of crossings points of
the smooth paths that cross at $e$ which together with $e$ form
the complete cyclic sequences of the crossings points of the same
smooth paths. Combinatorially, the two topologically distinct ways
to merge the smooth paths correspond to the replacement of the
pair of cyclic sequences $\{(e,R),(e,T)\}$ by one of the single
cyclic sequences: $(e',e'',R,e',e'',T)$ or
$(e',e'',R,e',e'',T^{-1})$, where $T^{-1}$ represent the inverse
of sequence $T$.

Observe that there is no possibility of real distinction between
the two ways of merging the smooth paths. But since both of them
work well, in the sense that the operation does not alter the
surface $S$, we simply choose one of the two possibilities
arbitrarily. By repeating this process we reach the point in which
there is just one cyclic sequence. An $sp$-embedding which
realizes this cyclic sequence, in a minimum-connectivity surface,
has the new digons introduced as bounding digons. By shrinking to
a point these bounding digons we obtain an $sp$-embedding of the
original graph.

By the above reduction and by Theorem 2.10, which gives a lower
bound on the connectivity of a surface in which there exists an
$sp$-embedding of a given graph, the problem of determining
whether a graph of maximum valency 3 is map-rich would be solved
if we knew an algorithm to determine what in the minimum
connectivity of a surface in which an $sp$-embedding of a graph
with one vertex and given cyclic sequence of edges around its
vertices exists.

With this reduction we conclude our analysis of the connection
between Gauss codes and map-rich graphs. In the next section we no
longer consider map-rich graphs.

\section*{3.B\quad Main Theorem and Crossing Functions}
We introduce some formal definitions related to the Gauss code
problem.

A {\it Gauss code} $X$ is a cyclic sequence in which each symbol
appears twice. The set of symbols of $X$ is denoted by ${\rm
symb}(X)$.

A Gauss code $X$ is {\it realizable} in a surface $S$ if there
exists a map $P$ with one $v$-gon, such that $X$ is equal to one
of the two $v$-orderings of $P$, and such that $S=S_M$, where $M$
is the phial of $P$.

Observe that if $P$ formally realizes $X$ in the sense of the
above definition, then $M'$ (the medial of the phial of $P$), is
the intuitive topological realization of $X$, as shown in Figure
3.1.

If $S$ is the plane, then the decision of realizability is the
classical Gauss code problem. In this section we extend the
solution of the realizability problem to the case of the
projective plane. We answer the following more general question:
is this (given) Gauss code realizable in a surface $S$ with
$\xi(S)\leq1$? The particularization of our conditions to the
plane reproduces Rosenstiehl's conditions [19]. Indeed, the
methods employed generalize those of [19].

A Gauss code gives the $v$-ordering of a map. In order to have a
complete description of a map we need to have, by Proposition 3.0,
and the definition of descriptor, the associated imbalance. This
indicates that to find a convenient 2-coloration of ${\rm
symb}(X)$ (defined next) is our aim in the (generalized) Gauss
code problem.

A \index{\it 2-colored Gauss code} is a Gauss code $X$ in which
${\rm symb}(X)$ is partitioned into ``black'' and ``white''
elements. We then let the set of white elements in ${\rm symb}(X)$
be the imbalance of the map whose $v$-ordering is given by $X$.
Therefore, we identify a 2-colored Gauss code with a map with one
$v$-gon.

If $P$ is a map with one $v$-gon, let us denote by one of the two
$v$-orderings of $P$. $\overline{P}$ is, therefore, identified
with a (not yet 2-colored) Gauss code. The input for the problem
is $\overline{P}$ and we seek a 2-coloration of ${\rm
symb}(\overline{P})$ such that the map defined $P$ has phial $M$
with the appropriate surface.

Observe that in the definition of realizability we can replace $P$
by $P^\sim$. It is easy to cheek (we do it in Proposition 3.2.1)
that the complement of each imbalance in $P$ is an imbalance in
$P^\sim$. These observations imply that realizability is
independent of the order of terms of the $\{$black,white$\}$
partition.

We also remark that any $\{$black,white$\}$ partition of a
2-colored Gauss code $P$ has an intuitive interpretation in $M'$,
even though we do not use this interpretation in the sequel. In
$M'$ the elements of ${\rm symb}(\overline{P})$ are vertices of
$H_M$. Suppose that we orient the edges of $H_M$ in the direction
of a traversal of the unique smooth path of $M'$. Then, the
2-coloration of the faces of $M'$ induces a 2-coloration of the
vertices of $H_M$: a vertex $x$ of $H_M$ is a black (white)
vertex, whence a black (white) element in ${\rm
symb}(\overline{P})$, if two of the darts directed toward $x$ or
away from $x$ encompasses a black (white) face of $M'$. For
instance, in Figure 3.1 the black elements are $\{2,3,7,9\}$ and
the white elements are $\{1,4,5,6,9\}$. Thus, we see that an
arbitrary 2-coloration of a Gauss code is realizable in the above
sense.

Given a Gauss code $X$, we denote by $i_X(p)$, for $p\in{\rm
symb}(X)$, the subset of ${\rm symb}(X)$ which is {\it
interlaced}\index{\it interlaced} with $p$. An element $p\in{\rm
symb}(X)$ is said to be interlaced with $q(\ne p)$ if $(p,q,p,q)$
is a subsequence of $X$. Note that since $X$ is cyclic, the
relation $i_X$ is symmetric.

The following theorem and its corollary are the main results of
the chapter.

\medskip\n{\bf (3.1) Theorem:} A Gauss code $X$ is realizable in a surface of
connectivity at most 1 if and only if it is possible to color
${\rm symb}(X)$ with two colors such as to satisfy the following
condition:

\noi (EQ)\qquad For every pair $y,z$ of elements in ${\rm
symb}(X)$ the following equivalence holds:
$$|i_X(y)\cap i_X(z)|+|i_X(Y)|\cdot|i_X(z)|$$
is odd iff $y$ and $z$ are interlaced and have the same color.

We want to point that condition (EQ) was suggested by D.H.
Younger. Theorem 3.1 has a corollary which gives a good algorithm
to decide if the 2-coloring exists, and in this case the corollary
also provides all the 2-colorings which makes (EQ) true.

To state the corollary of Theorem 3.1 it is convenient to define
the concept of interlace graph induced by a Gauss code. Given a
Gauss code $X$ the {\it interlace graph} induced by $X$, denoted
by $I_X$, is the simple graph whose set of vertices is ${\rm
symb}(X)$ and whose set of edges is formed by the pairs of
elements in ${\rm symb}(X)$ which interlace.
\medskip

\n{\bf Corollary of (3.1).} {\bf (Good algorithm for the
realizability of Gauss codes in surfaces with $\xi\leq1$):} Given
a Gauss code $X$ there exists a 2-coloration which satisfies (EQ)
of Theorem 3.1 if and only if the following conditions hold:\jtt
\begin{itemize}
\item[(i)] For every pair $\{y,z\}$ of elements in ${\rm symb}(X)$
which is not an edge of $I_X$ the sum $|i_X(y)\cap
i_X(z)|+|i_X(y)|\cdot|i_X(z)|$ is even. \item[(ii)] The subset of
$eI_X$ given by
$$T=(\{y,z\}\in eI_X:|i_X(y)\cap
i_X(z)|+|i_X(y)|\cdot|i_X(z)|\quad\mbox{is even}$$ is a coboundary
in the interlace graph $I_X$.
\end{itemize}

\noi Moreover if (i) and (ii) hold, then the 2-colorings of ${\rm
symb}(X)$ which make (EQ) true are precisely the ones of the form
$\{B,W\}$, where $T=\delta_{I_X}(B)$ and $W={\rm
symb}(X)\backslash B$.

\inicproof Initially we show that if conditions (i) and (ii) are
true, then the rule to form the 2-colorings $\{B,W\}$ makes (EQ)
true. Thus, let $\{B,W\}$ be a suitable 2-coloring. The pairs of
equal elements always satisfy (E0), independently of the
2-coloring. If $\{y,z\}$ is not an edge of $I_X$, then (i) implies
that $|i_X(y)\cap i_X(z)|+|i_X(y)|\cdot|i_X(z)|$ is even; hence,
(EQ) is verified. If $\{y,z\}\in eI_X$, then condition (ii) and
the rule to obtain $\{B,W\}$ imply that $y$ and $z$ are of
different colors iff $\{y,z\}\in T$, i.e.,
$|i_X(y)i_X(z)|+|i_X(y)|\cdot|i_X(z)|$ is even. Observe that also
in this case (EQ) holds.

Suppose now that $\{B,W\}$ is a 2-coloring that satisfies (EQ). If
$\{y,z\}\not\in eI_X$, then (EQ) implies that $|i_X(y)\cap
i_X(z)|+|i_X(y)|\cdot|i_X(z)|$ is even; hence, (i) is verified. By
(EQ), the subset of $eI_X$ whose elements have ends of different
colors is precisely $T$. This implies that the graph obtained from
$I_X$ by contracting the edges not in $T$ is bipartite. Whence,
$T$ is a coboundary in $T_X$, proving that (ii) holds. The proof
of the corollary is complete. \fimproof

Observe that, by the above corollary the number of essentially
distinct ways to choose the 2-coloring, and thus realize the Gauss
code, is $2^{p(I_X)-1}$ where $p(I_X)$ is the number of components
of $I_X$.

By a result that we prove later in the chapter, (Corollary of
3.3(a),3.4 , and 3.5(a), presented immediately after the proof of
Lemma 3.5,) the surface where a Gauss code $X$ is realizable is
orientable iff $I_X$ is eulerian. This property permits the
following specialization for the plane of the conditions of the
above corollary:

\noi{\bf Rosdnstiehl's conditions:} A Gauss code $X$ is realizable
in the plane if and only if the following conditions hold:\jtt
\begin{itemize}
\item[(I)] $I_X$ is eulerian; \item[(II)] for $\{y,z\}\not\in
eI_X,\ |i_X(y)\cap i_X(z)|$ is even; \item[(III)] $T=\{\{y,z\}\in
eI_X:|i_X(y)\cap i_X(z)|\quad\mbox{is even}\}$ is a coboundary in
$I_X$.
\end{itemize}

Conditions (I), (II), (III) appear in Rosenstiehl [19]. The proof
of Theorem 3.1 takes most of the remaining of the chapter. To
establish the proof we develop a theory of ``crossing functions''
which generalizes for arbitrary surfaces the ``projections'
described for the plane by Shank in [24].

Given a 2-colored Gauss code $X$ we now define the {\it crossing
function}\index{\it crossing function} induced by $X$
$$c_X:2^{{\rm symb}(X)}\longrightarrow 2^{{\rm symb}(X)}.$$
Function $c_X$ is given on the singletons by
$c_X(\{a\})=i_X(a)+t_X(a)$, where
$$t_X(a)=\left\{\begin{array}{l}
\phi\quad\mbox{if}\quad a\quad\mbox{is white}\\ [0.5cm]
\{a\}\quad\mbox{if}\quad a\quad\mbox{is black},
\end{array}\right.$$
and extended for arbitrary $A\subseteq{\rm symb}(X)$ as follows:
$c_X(A)=\sum_{a\in A}c_X(\{a\})$, where the sum is $\mod\,2$ sum
of subsets. It is also convenient to extend $t_X$ and $i_X$ to
$2^{{\rm symb}(X)}$ by the analogous formulas: for $A\subseteq{\rm
symb}(X)$,
$$t_X(A)=\sum_{a\in A}t_X(a)$$
$$i_X(A)=\sum_{a\in A}i_X(a)$$

In the proof of Theorem 3.1 we denote the abstract Gauss code $X$
by $\overline{P}$, with the intention of working with the
terminology of maps that we have developed.

Given a Gauss code $\overline{P}$, denote by $ODD(\overline{P})$
and $EVEN(\overline{P})$ the following partition of ${\rm
symb}(\overline{P})$:
$$ODD(\overline{P})=\{y\in{\rm symb}(P):\ |i_P(y)|\quad\mbox{is odd}\}$$
$$EVEN(\overline{P})=\{y\in{\rm symb}(P):\ |i_P(y)|\quad\mbox{is even}\}$$
To simplify the notation, we use $i_P(X),c_P(X),t_P(X)$, instead
of the more correct $i_{\overline{P}}(\{x\}), c_P(\{x\}),
t_{\overline{P}}({x})$. \vfill
$$\includegraphics[width=14.50cm,height=11cm]{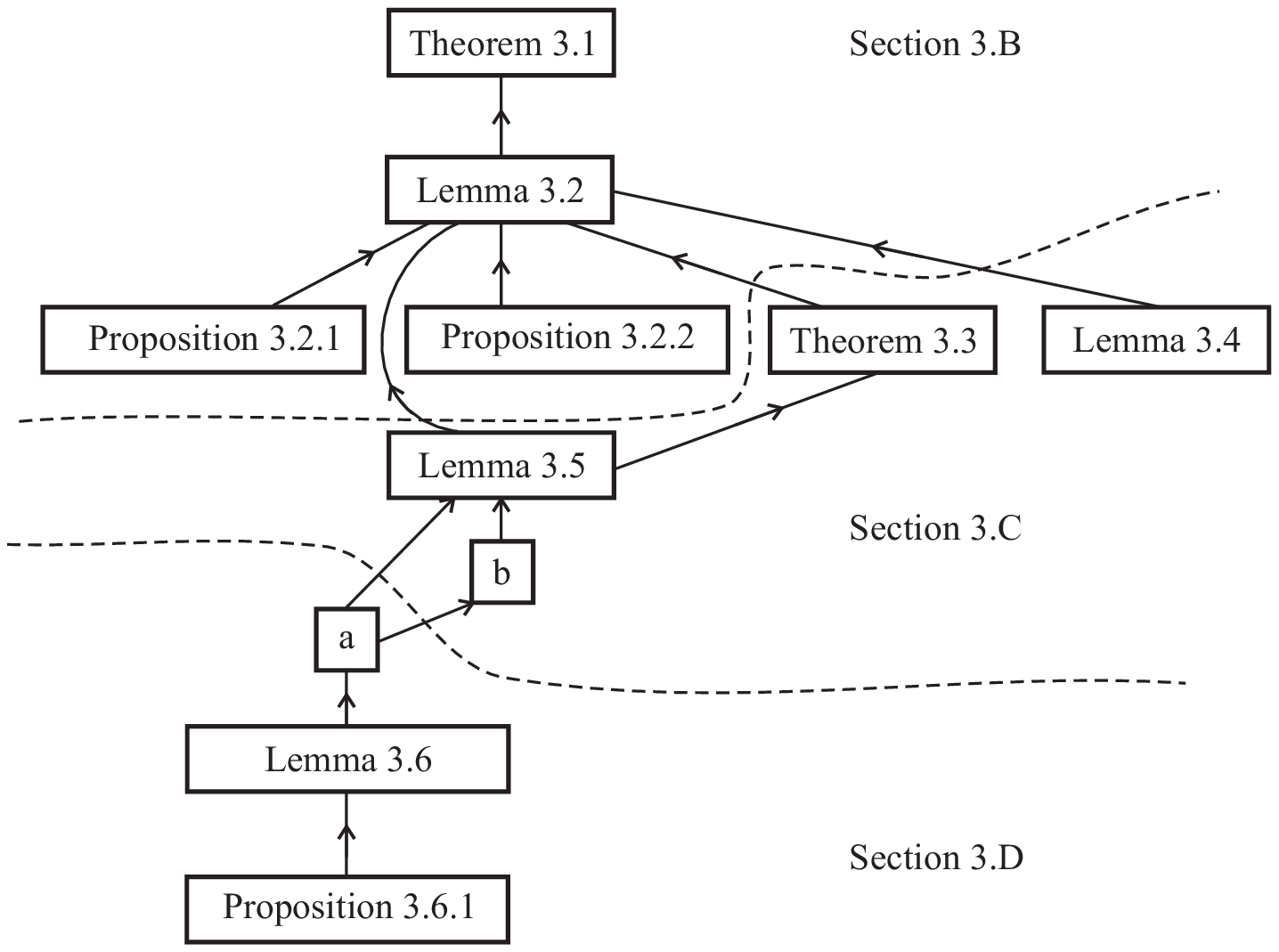}$$
\centerline{Dependence Structure for Theorem 3.1} \eject
\medskip\noi{\bf Proof of Theorem 3.1.} Lemma 3.2, proved below,
states that a Gauss code $\overline{P}$ is realizable in $S$ with
$\xi(S)\leq 1$ iff there exists a 2-coloring of ${\rm
symb}(\overline{P})$ inducing a map $P$ with one $v$-gon, (the
2-coloring gives the imbalance relative to the v-ordering given by
$\overline{P}$), and a crossing function $c_P$ which satisfies the
following equations:

(a)\qquad$(c_P+c^2_P)(x)=\phi,\quad\mbox{for}\quad x\in
EVEN(\overline{P})$;

(b)\qquad$(c_P+c^2_P)(x)=ODD(\overline{P}),\quad\mbox{for}\quad
x\in ODD(\overline{P})$.

We assume Lemma 3.2 and show the equivalence of condition (EQ) of
Theorem 3.1 with conditions (a) and (b) of Lemma 3.2. We begin by
deriving a useful expression for $(c_P+c^2_P)(x)$. For $z\in{\rm
symb}(\overline{P})$ let the variable $\gamma_z\in\{0,1\}=GF(2)$
assume value 1 if $z$ is black and 0 if $z$ is white. Since
$c_P=i_P+t_P$ and $t^2_P=t_P$, we obtain
$$\begin{array}{rcl}
(c_P+c^2_P)(x)&=&i^2_P(x)+i_Pt_P(x)+t_Pi_P(x)+i_P(x)\\ [0.5cm]
&=&\dy i^2_p(x)+\sum_{y\in i_P(x)}(\gamma_x+\gamma_y+1)y.
\end{array}$$

It is convenient to separate the equivalence (EQ) into two
equivalences $(a')$ and $(b')$ stated below:\jtt
\begin{itemize}
\item[$(a')$] if $|i_P(x)|\cdot|i_P(z)|$ is even, then
$|i_P(x)\cap i_P(z)|$ is odd iff $\gamma_x=\gamma_z$ and $x\in
i_P(z)$; \item[$(b')$] if $|i_P(x)|\cdot|i_P(z)|$ is odd, then
$i_P(x)\cap i_P(z)|$ is even if $\gamma_x=\gamma_z$ and $x\in
i_P(z)$.
\end{itemize}

It is clear that statements $(a')$ and $(b')$ together are
equivalent to (EQ) of Theorem 3.1. Our aim in the remaining of the
present proof is to establish the equivalence
$$(a)\quad\mbox{and}\quad (b)\Leftrightarrow(a')\quad\mbox{and}\quad(b').$$

A helpful remark in this direction is that for $A\subseteq{\rm
symb}(\overline{P}),\ z\in i_P(A)$ iff $|i_P(z)\cap A|$ is odd.
Whence, $z\in i^2_P(x)$ iff $|i_P(z)\cap i_P(x)|$ is odd.

With the above remark and the symmetry of $i_P$ and of $i^2_P,\
(a')$ can be restated as follows: if $x\in EVEN(\overline{P})$,
then $z\in i^2_P(x)$ iff $\gamma_x=\gamma_z$ and $x\in i_P(z)$.
This statement, by the expression calculated above for
$(c_P+c^2_P)(x)$, is equivalent to the following one: if $x\in
EVEN(\overline{P})$, then $(c_P+c^2_P)(x)=\phi$, which is
precisely statement $(a)$. Therefore, we have established the
equivalence $(a)\Leftrightarrow (a')$.

To conclude, we show the equivalence
$$(a)\quad\mbox{and}\quad (b)\Leftrightarrow(a)\quad\mbox{and}\quad(b').$$
Statement $(b')$ can be rewritten as follows: if $x$ and $z$ are
both in $ODD(\overline{P})$, then $z\in i^2_P(x)$ iff
$\gamma_x=\gamma_z$ and $x\in i_P(z)$. By the expression for
$(c_P+c^2_P)(x)$, this shows that $(b')$ is equivalent to the
following statement: if $x$ and $z$ are both in
$ODD(\overline{P})$, then $z\in(c_P+c_P)(x)$. It is then clear
that $(b)$ implies $(b')$ and that $(a)$, the symmetry of the
relation $c_P+c^2_P$, and $(b')$ imply (b).

This concludes the Proof of Theorem 3.1, provided Lemma 3.2 is
established. \fimproof

\medskip\n{\bf (3.2) Lemma:} A Gauss code $\overline{P}$ is realizable in
surface $S$ of connectivity at most one if and only if there
exists a 2-coloring of ${\rm symb}(\overline{P})$ inducing a map
$P$ with a crossing function $c_P$ which satisfies the following
equations:

(a)\qquad $(c_P+c^2_P)(x)=\phi,\quad\mbox{for}\quad x\in
EVEN(\overline{P})$;

(b)\qquad $(c_P+c^2_P)(x)=ODD(\overline{P}),\quad\mbox{for}\quad
x\in ODD(\overline{P})$;
\medskip

\inicproof Assume that there exists a 2-coloring of ${\rm
symb}(\overline{P})$ such that equations (a) and (b) hold. These
equations 2 imply that the dimension of the image of $c_P+c^2_p,\
\dim[Im(c_P+c^2_P)]$, is at most one. The corollary of Proposition
3.2.1, which we prove next, states that $c_P+c_{P^\sim}$ is the
identity. Whence, $c_{P^\sim}\circ c_P=c_P+c^2_P$ and thus,
$\dim[Im(c_{P^\sim}\circ c_P)]\leq1$. By Theorem 3.3(a), proved
below, we have $\xi(M)=\dim[Im(c_{P^\sim}\circ c_P)]$, where $M$
is the phial of $P$. Therefore, $\xi(M)\leq1$, and this proves
that $P$ is realizable in a surface of connectivity at most one.
This part of the proof is complete, provided the corollary of
Proposition 3.2.1 as well as Theorem 3.3(a), are established.

Conversely, assume that the Gauss code $\overline{P}$ is
realizable in the plane or In the projective plane. This means
that there exists a 2-coloring of ${\rm symb}(\overline{P})$ such
that the map $P$ formed has phial $M$ satisfying $\xi(M)\leq1$.
For such maps, Proposition 3.2.2, which we prove below, states
that the $s$-cycles in $M^\perp$ are precisely the elements of the
coboundary space $F$ of $G_D$ (where $D$ is the dual of $M$).
Lemma 3.4, also proved below, states that for maps $M$ whose phial
$P$ has one $v$-gone $c_P(x)$ is an $s$-cycle in $M^t$ if and only
if $|i_P(x)|$ is even, i.e., $x\in EVEN(\overline{P})$. Putting
together 3.2.2 and 3.4 we obtain the following condition for maps
$P$ with one $v$-gon and $\xi(M)\leq1:\ c_P(x)$ is in $F$ iff
$x\in EVEN(\overline{P})$.

Lemma 3.5(b), proved below, when applied to map $P$ states that
${\rm Ker}(c_{P^\sim})=F$. Therefore, $c_P(x)$ is in ${\rm
Ker}(c_{P^\sim})$ iff $x\in EVEN(\overline{P})$. This condition
implies $(c_{P^\sim}\circ c_P)(x)=\phi$, for $x\in
EVEN(\overline{P})$. Since, as mentioned before $c_{P^\sim}\circ
c_P=c_P+c^2_P$, this concludes (a), modulo the unproved
statements, whose proofs are supplied in the sequel.

To prove (b) consider $x\in ODD(\overline{P})$. Recall that
$$(c_P+c^2_P)(x)=i^2_P(x)+\sum_{y\in i_P(x)}(\gamma_x+\gamma_y+1)y,$$
an expression that we computed in the Proof of Theorem 3.1. Also
recall that for $z\in{\rm symb}(\overline{P})$, $\gamma_z=1$ iff
$z$ is a black element of ${\rm symb}(\overline{P})$. Therefore,
$y\in(c_P+c^2_P)(x)$ iff one of the three following conditions
holds:\jtt
\begin{itemize}
\item[(i)] $y\in i^2_P(x),\quad y\in i_P(x)$\quad
and\quad$\gamma_x\ne\gamma_y$; \item[(ii)] $y\in i^2_P(x),$\quad
and\quad $y\not\in i_P(x)$; \item[(iii)] $y\not\in i^2_P(x),\quad
y\in i_P(x)$\quad and\quad$\gamma_x=\gamma_y$.
\end{itemize}

If $y\in EVEN(\overline{P})$ we show that
$y\not\in(c_P+c^2_P)(x)$. Note that $y\in(c_P\circ c^2_p)(x)$ iff
$x\in(c_P+c^2_P)(y)$, property easily deduced from (i), (ii), and
(iii) above. But $(c_P+c^2_P)(y)$ is the null subset, by part (a).
Thus, if $y\in EVEN(\overline{P})$, then $y\not\in(c_P+c^2_P)(x)$.
To establish (b) it suffices, then, to prove that
$y\in(c_P+c^2_P)(x)$, for $y\in ODD(\overline{P})$. To this end
initially observe that $z\in(c_P+c^2_P)(z)$, for every $z\in
ODD(\overline{P})$. This follows from $z\in i^2_P(z)$ and
$z\not\in i_P(z)$. Since $\xi(M)\leq1$ and
$c_P+c^2_P=c_{P^\sim}\circ c_P$, Lemma 3.3(a) proved below implies
that $\dim[Im(c_P+c^2_P)]\leq1$. Whence, there are one or two
elements in $Im(c_P+c^2_P)$ one of which is the null vector. Thus,
we conclude that for $y\in ODD(\overline{P}),\
(c_P+c^2_P)(x)=(c_P+c^2_P)(y)$, because both images are non-null.
Now to finish the proof assume an arbitrary $y\in
ODD(\overline{P})$. We show that at least one of the conditions
(i), (ii), (iii) above is verified. Suppose, as a first case that
$y\in i^2_P(x)$. If $y\not\in i_P(x)$, then (ii) is verified.
Therefore, assume that $y\in i_P(x)$. It $\gamma_x=\gamma_y$, then
(i) is verified. Suppose then, that $\gamma_x=\gamma_y$. We are
under the hypotheses $y\in i^2_P(x),\ y\in i_P(x)$ and
$\gamma_x=\gamma_y$. By the expression for $(c_P+c^2_P)(x)$ this
implies that $y\not\in(c_P+c^2_P)(x)$. But we know that
$y\in(c_P+c^2_P)(y)$ and that $(c_P+c^2_P)(x)=(c_P+x^2_P)(y)$.
Therefore, $y\in(c_P+c^2_P)(x)$, which is a contradiction,
concluding the analysis in the case that $y\in i^2_P(x)$. Assume
the second case, $y\not\in i^2_P(x)$. Since $y\in(c_P+c^2_P)(x)$
we conclude that $y\in i_P(x)$ and that $\gamma_x=\gamma_y$. This
is precisely condition (iii). This concludes the proof of (b). The
proof of the lemma is, then, reduced to the proofs of the
Corollary of 3.2.1 and of the following statements: 3.2.2, 3.3(a),
3.4, and 3.5(b). \fimproof

For the next proposition and its corollary we consider, as usual,
$SQ(X)\cong SQ(M)$ if $X\in\Gamma(M)$.

\medskip\n{\bf (3.2.1) Proposition:} Given any map $M$, $T\subseteq SQ(M)$ is an
imbalance in $M$ iff $SQ(M)\backslash T$ is an imbalance in
$M^\sim$.
\medskip

\inicproof Assume that $T$ is an imbalance in $M$. By definition
there exists a balancing partition $B_T$ for $M$ which has $T$ as
its imbalance. Since $M^\sim$ has the same $v$-gons as $M$, $B_T$
is also a balancing partition for $M^\sim$. Given a corner $X$ of
$M$ or of $M^\sim$ , $X$ and $v_M(X)= v_{M^\sim}(X)$ are not in
the same class of $B_T$. This implies that two corners of a square
which are in the same class of $B_T$ are adjacent in $C_M$ iff
they are opposite in the corresponding square of $M^\sim$. The
converse follows by the same argument, since the operations of
taking the antimap of a map and taking the complement of a subset
are involutions. \fimproof

\noi{\bf Corollary of Proposition 3.2.1:} For a map $P$ with one
$v$-gon, the linear transformation
$$c_P+c_{P^\sim}:\ 2^{SQ(P)}\longrightarrow 2^{SQ(P)}$$
is the identity.

\inicproof The corollary follows from the definitions of $c_P$ and
$c_{P^\sim}$, because the above proposition implies that a symbol
of the associated 2-colored Gauss codes is black in $P$ iff it is
white in $P^\sim$. \fimproof

\medskip\n{\bf (3.2.2) Proposition:} Given a map $M$, if $\xi(M)\leq1$, then the $s$-cycles of $M^t$
are precisely the elements of the coboundary space $F$ of $G_D$,
where $D$ is the dual of $M$.
\medskip

\inicproof If $\xi(M)=0$, the proposition follows from (2.6),
since the hypothesis implies that $V^\perp=F$. If $\xi(M)=1$, then
by Theorem 1.7 $M$ is non-orientable, since the parity of $\xi(M)$
and $\chi(M)$ agree. Take an odd polygon in $C_M$ and consider its
image under $\psi^M_c$. This image is an $r$-cycle in $M^t$, which
we denote by 4x$. The cycle 4x$ is not in $F$, since the elements
of $F$ are $s$-cycles. Therefore, by Proposition 2.6, it follows
that
$$V^\perp/F=\{F,x+F\}.$$

Consider an $s$-cycle $y$. Either $y\in F$ or else $y\in x+F$. In
the latter case, $x+y\in F$. Since $\psi^M_c$ is a homomorphism,
by Proposition 1.4, the symmetric difference of an $r$-cycle and
an $s$-cycle is an $r$-cycle. Hence, $x+y$ is an $r$-cycle which
is contained in $F$. This is not possible. Therefore, $y\in F$,
establishing the proposition. \fimproof

\section*{3.C\quad Properties of Crossing Functions}
In this section we continue the proof of the main theorem (Theorem
3.1). At the same time we establish some theorems about general
crossing functions associated with $t$-maps with one vertex. Items
(b) and (c) of the next theorem are not used in connection with
Theorem 3.1, but are included for completeness.

\medskip\n{\bf (3.3) Theorem:} if $P$ is a map with one $v$-gon, then we have

(a)\qquad $\dim[Im(c_{P^\sim}\circ c_P)]=\xi(M)$;

(b)\qquad $Im(c_{P^\sim}\circ c_P)=V^\perp\cap F^\perp$;

(c)\qquad ${\rm Ker}(c_{P^\sim}\circ c_P)=V+F$,

\noi where $M$ is the phial of $P$, $D$ is the dual of $M$, and
$V$ and $F$ are the coboundary spaces of $G_M$ and $G_D$,
respectively.
\medskip

\inicproof To prove this theorem we make use of Lemma 3.5, proved
below. By Lemma 3,5(a) we have $Im(c_P)=V^\perp$. Therefore,
$$Im(c_{P^\sim}\circ c_P)=Im(c_{P^\sim}|V^\perp),$$
where, ``$|V^\perp$'' stands for restriction to $V^\perp$. By the
fundamental theorem for homomorphisms applied to
$c_{P^\sim}|V^\perp$ we have
$$\dim[Im(c_{P^\sim}|V^\perp)]+\dim({|rm Ker}(c_{P^\sim}|V^\perp)]=\dim V^\perp.$$
By lemma 3.5(b), applied to map $P^{\sim}$, it follows that ${\rm
Ker}(c_P)=F$. By Proposition 2.2, $F\subseteq V^\perp$, implying
${\rm Ker}(c_{P^\sim}|V^\perp)=F$. Therefore,
$$\dim[{\rm Ker}(c_{P^\sim}|V^\perp)]=f-1,$$
where $f=|VG_D|$ (Observe that since $P$ has one $v$-gon, $D$ is
connected). Let $v=|VG_M|$ and $e=|eG_M|$. By the above equations,
it follows that
$$\begin{array}{rcl}
\dim[Im(c_{P^\sim}\circ c_P)]&=& (e-v+1)-(f-1)\\ [6pt] &=&
(e+2)-(v+f)\\ [6pt] &=& \xi(M)
\end{array}
$$

The proof of part (a) is complete, provided Lemma 3.5 is
established.

We proceed to prove (b). It follows, from the Corollary of
Proposition 3.2.1, that $c_{P^\sim}\circ P=c_{P^\circ}$. This
commutativity and the facts that, by Lemma 3.5(b),
$Im(c_P)=V^\perp$ and $Im(c_{P^\sim})=F^\perp$ imply the inclusion
$Im(c_{P^\perp}\circ c_P)\subseteq V^\perp \cap F^\perp$. By the
triple inclusion property, lemma 2.9.1,
$$\dim(V^\perp\cap F^\perp)-\dim(V\cap F)=\xi(M).$$

By the absorption property, Lemma 2.5, it follows that $V\cap
F=V\cap F\cap Z$, with $Z$ coboundary space of $G_P$. Since
$|VG_P|=1,\ Z$ is the null space. Therefore,
$$\dim(V^\perp\cap F^\perp)=\xi(M).$$
By part (a) $\dim[Im(c_{P^\sim}\circ c_P)]=\xi(M)$. Thus,
$Im(c_P\circ c_{P^\sim})$ is a subspace of $V^\perp\cap F^\perp$
which has the same dimension as the whole space. Therefore,
$$Im(c_{P^\sim}\circ c_P)=V^\perp\cap F^\perp;$$
this concludes (b).

To prove (c) we also use Lemma 3.5. Take $x\in V+F$; say that
$x=u+w$ with $u\in V$ and $w\in F$. We have
$$(c_{P^\sim}\circ c_P)(u+v)=c_{P^\sim}(c_P(u)+c_P(w))=C_{P^\sim}(c_P(w)).$$

The latter equality follows because $c_P(u)=\phi$ by Lemma 3.5(b).
As $c_P$ and $c_{P^\sim}$ commute we have
$$c_{P^\sim}c_P(w)= c_P\circ c_{P^\sim}(w)=c_P(\phi)=\phi.$$
That $c_{P^\sim}(w)$ follows from Lemma 3.5(b) applied to map
$P^\sim$. We conclude then that $V+F\subseteq{\rm
Ker}(c_{P^\sim}\circ c_P)$.

The value of $\dim[{\rm Ker}(c_{P^\sim}\circ c_P)]$ is $e-\xi(M)$,
where $e=|SQ(M)|$. This follows from the fundamental theorem for
homomorphisms and from part (a). The value of
$\dim(V+F)=(v-1)+(f-1)$, where $v=|VG_M|,\ f=|VG_D|$. This follows
from 2.0.1(a), and from the fact that $V\cap P$ is the null space.
Recall that $V\cap F=V\cap F\cap Z$, and $Z$ is the null space,
since $Z$ is the coboundary space of $G_P$ which has only one
vertex. We can rewrite $(v-l)+(f-1)$ as $e-\xi(M)$, by the
definition of $\xi(M)$. Since $V+F$ is a subspace of ${\rm
Ker}(c_{P^\sim}\circ c_P)$ and has the same dimension, it follows
that it is equal to ${\rm Ker}(c_{P^\sim}\circ c_P)$. This
concludes (c) and the proof of the Theorem 3.3, provided Lemma 3.5
is established. \fimproof

\medskip\n{\bf (3.4) Lemma:} For every map $M$ whose phial $P$ has one
$v$-gon, $c_P(x)$ is a cycle in $G_M$, for every $x\in{\rm
symb}(\overline{P})$ and
$$c_P(x)\quad\mbox{is an $s$-cycle in}\quad
M^t\quad\mbox{iff}\quad|i_P(x)|\quad\mbox{is even}.$$
\medskip

\inicproof To prove this lemma, which is of fundamental importance
in the chapter, we define a function $c'_P$ on the corners of $M$
with image in the set of cycles of $C_M$. For a corner $X$ of $M$,
the cycle $c'_P(X)$ is defined by the edges of $C_M$ occurring
once in the reentrant path which starts at $z_M(X)$ and proceeds
by using $a_{M^-},\ v_{M^-}$, and $f_M$-edges (in this order)
until it reaches another corner of $x$, which denotes the square
to which $x$ belongs. This corner can be $v_M(X)$, in which case
we close the path by using the $f_M$-edge which links $v_M(X)$ to
$z_M(X)$; it also can be $f_M(X)$, in which case we close the path
by using the $v_M$-edge which links $f_m(X)$ to $z_m(X)$. The fact
that $x$ is a loop in $G_P$ and a parity argument show that the
first corner of $x$ reached by the path is not $X$, and is not
$z_M(X)$. This completes the definition of $c'_P$.

We claim that the image of $c'_P(X)$ under $\psi^M_c$ is $c_P(x)$,
independently of the corner $X$ of $x$. To prove this claim
observe that there is just one edge of square $x$ in $c'_P(X)$,
for $X$ corner of $x$. If $x$ is an $s$-loop in $p^t$ ($x$ is a
black element of ${\rm symb}(\overline{P})$), it follows from the
definitions that the edge of $x$ present in $c'_P(X)$ is an
$f_M$-edge; otherwise, if $x$ is an $r$-loop in $p^t$, this edge
is a $v_m$-edge. Since the $f_M$-edges are not contracted under
$\psi$ (the canonical contraction that transforms $(C_M,S_M)$ into
$(G_M,S_M)$), square $x$ is an edge in $c_P(x)$, in the case that
$x$ is an $s$-loop in $p^t$. Since the $v_M$-edges are contracted
under $\psi$, square $x$ is not an edge in $c_P(x)$, in the case
that $x$ is an $r$-loop in $p^t$. The other edges belonging to the
image of $c'_P(X)$ under $\psi^M_c$ correspond to the squares
whose intersection with $c'_P(X)$ is a path of length 2, that is,
the squares which are met once by the reentrant path defining
$c'_P(X)$. Those squares are precisely the squares in $i_P(x)$.
The claim that $\psi^M_c(c'_P(X))=c_P(x)$ for every corner $X$ of
square $x$ is, therefore, proved.

By the definition of $s$- and $r$-cycles, to complete the proof of
the lemma it is enough to show that the congruence
$$|c'_P(X)|=|i_p(x)| \quad\mbox{mod}\ 2$$
holds for every square $x$ with $X$ corner of $x$.

Since $C_M$ is cubic, the cycle $c'_P(X)$ induces a subgraph of
$C_M$ which consists of a certain number of disjoint polygons,
whose set is denoted by $\Omega$. We count the vertices of these
polygons. In square $x$ there are two corners which are vertices
of a polygon in $\Omega$. If a square is met twice by the
reentrant path which defines $c'_P(X)$, then all its four corners
are vertices of polygons in $\Omega$. If a square is met once by
the reentrant path, then three of its corners are vertices of a
polygon in $\Omega$. Evidently, if a square is not met by the
reentrant path, none of its corners is in a polygon in $\Omega$.
Hence, the parity of the number of vertices of the polygons in
$\Omega$, which is the same as the parity of $|c'_P(X)|$, is equal
to the parity of $|i_P(x)|$. This establishes the Lemma. \fimproof

\medskip\n{\bf (3.5) Lemma:} Let $P=P_M$ be a map with one $v$-gon, and
$V$ the coboundary space of $G_M$. We then have the following
equalities for the image and the kernel of the crossing function
$c_P: 2^{SQ(P)}+2^{SQ(P)}$ induced by P:

(a)\quad $Im(c_P)=V^\perp$

(b)\quad ${\rm Ker}(c_P)=V$
\medskip

\inicproof By Lemma 3.4 above, $c_P(x)$ is in $V^\perp$ for every
$x\in{\rm symb}(P)$. Since the image of $c_P$ in generated by the
members of
$$\{c_P(x):\ x\in{\rm symb}(\overline{P})\},$$
it follows that $Im(c_P)\subseteq V^\perp$.

The proof of the converse inclusion needs some additional
preparation. We give it in the next section, Lemma 3.6. We assume
that (a) holds for the proof of (b).

Take $s\in{\rm symb}(\overline{P})$ and $T\subseteq{\rm
symb}(\overline{P})$. We show that the scalar product
$\aw\{s\},c_P(T)\fw$ is equal to $\aw c_P(\{s\}),T\fw$ as follows:
$$\begin{array}{rcl}
\aw\{s\},c_P(T)\fw &=&\dy\sum_{x\in
T}\aw\{s\},i_P(x)\fw+\sum_{z\in T}\aw\{s\},t_P(x)\fw\\ [0.5cm]
&=&\dy\sum_{x\in T} \aw i_P(s),x\fw+\sum_{z\in T}\aw t_P(s),x\fw\\
[0.5cm] &=&\aw c_P(s),T\fw.
\end{array}$$

Let $T\subseteq{\rm symb}(\overline{P})$ be an element of ${\rm
Ker}(c_P)$. Since $\aw x,c_P(T)\fw=\aw c_P(x),T\fw$ and
$c_P(T)=\phi$, it follows that $\aw c_P(x),T\fw=0$ for every
$x\in{\rm symb}(\overline{P})$. As we are assuming (a), $V^\perp$
is generated by the members of $\{c_P(x): x\in{\rm
symb}(\overline{P})\}$. It follows that $T$ belongs to the space
perpendicular to $V^\perp$ namely, $V$. Hence ${\rm
Ker}(c_P)\subseteq V$.

We first learned how to prove that $V\subseteq{\rm Ker}(c_P)$ from
H. Shank and subsequently we discover the following shorter proof.

Take a coboundary $T\subseteq{\rm symb}(\overline{P})\cong
e(G_M)$; thus $T$ belongs to $V$. Since by part (a), for $s\in{\rm
symb}(\overline{P}),\ c_P(s)$ is a cycle in $G_M$, it follows that
$\fw c_P(s),T\fw=0$ for every $s\in eG_M$. This is the same as
$\aw s,c_P(T)\fw=0$ for every $s\in eG_M$, by the equality
established above. Hence, $c_P(T)=\phi$, proving that $T\in{\rm
Ker}(c_P)$. This concludes the proof of the lemma, provided the
inclusion $V^\perp\subseteq Im(c_P)$ is established. We treat this
inclusion in Lemma 3.6. \fimproof

The following corollary is an algebraic-combinatorial
characterization of realizability of a Gauss code in a given
arbitrary surface:

\noi{\bf Corollary of (3.3.a),(3.4) and (3.5.a):} A Gauss code is
realizable in a surface of connectivity $\xi$ iff there exists a
$\{$black,white$\}$ partition of ${\rm symb}(\overline{P})$ such
that the induced map with one $v$-gon $P$ satisfies
$$\dim[Im(c_{P^\sim}\circ c_P)]=\xi.$$
Moreover, the surface where $\overline{P}$ in realizable is
orientable iff $I_{\overline{P}}$ is eulerian.

\inicproof The first part of the corollary, treating the
connectivity, is a straightforward application of the definition
of realizability and of the Theorem 3.3(a).

To establish the second part, observe that for every $P$ induced
by a $\{$black, white$\}$ partition Lemma 3.5(a) states that the
members of $\{c_P(x): x\in{\rm symb}(\overline{P})\}$ generate the
cycle space $V^\perp$ of $G_M$ where $M$ is the phial of $P$. By
Lemma 3.4.2, $c_P(x)$ is an $s$-cycle in $M^t$ iff $x\in
EVEN(\overline{P})$. Note that the $\mod\,2$ sum of cycles in
$G_M$ is an $s$-cycle iff the number of $r$-cycles is even.
(Recall that $\psi^M_c:CS(M)\to CS(G_M)=V^\perp$ is a
homomorphism, result proved in Proposition 1.4.) Hence, $V^\perp$
has only $s$-cycles iff $x\in EVEN(\overline{P})$ for every
$x\in{\rm symb}(\overline{P})$. But $V^\perp$ has only $s$-cycles
iff $C_M$ is bipartite, or iff $S_M$ is orientable. Since $S_M$ is
the surface where $\overline{P}$ is realized (once the
$\{$black,white$\}$ partition is chosen), the corollary follows.
\fimproof

\section*{3.D\quad Surjectivity of the Crossing Functions}
In this section we conclude the proofs of the statements which are
still open with the proof that $V^\perp\subseteq Im(c_P)$.
Initially we define some concepts, leading to the idea of
``reflexivity\index{reflexivity}''.

Consider a graph $L$ and an arbitrary, but fixed, orientation for
its edges. Let $\pi$ be a cyclic path in $L$ such that each edge
of $L$ appears twice in $\pi$. Given a pair $(L,\pi)$ under these
conditions, we define cyclic sequence $\pi'$ as follows. The
symbols in $\pi'$ are labels for the edges of $G$. The order in
which they occur is given by cyclic path $\pi$. Moreover each
symbol $x$ in $\pi'$ is signed ``$+$'' or ``$-$'' according to
whether or not edge labelled $x$ is traversed in the right
direction (according to the fixed orientation for the edges of
$L$).

We identify $\pi'$ with a 2-colored Gauss code, where the black
elements are the symbols whose two occurences in $\pi'$ are with
the same sign and the white elements are the others. Functions
$i_\pi$, and $c_\pi$, are, therefore, defined.

Given an abstract cyclic sequence $\pi'$ in which each symbol is
signed and occurs twice, we define a digraph $\Omega(\pi')$ as
follows. The darts of $\Omega(\pi')$ are the symbols which occur
in $\pi'$. To obtain the vertices of $\Omega(\pi')$ we start with
a head and a tail for each dart and effect the following
identifications, one for each successive pair of signed symbols in
$\pi'$:

if $(+x,+y)$ appears in $\pi'$, then identify head $(x)$ and tail
$(y)$;

if $(+x,-y)$ appears in $\pi'$, then identify head $(x)$ and head
$(y)$;

if $(-x,+y)$ appears in $\pi'$, then identify tail $(x)$ and tail
$(y)$;

if $(-x,-y)$ appears in $\pi'$, then identify tail $(x)$ and head
$(y)$.

Given $(L,\pi)$ we say that cyclic path $\pi$ is {\it reflexive}
in digraph $L$ if $\Omega(\pi')=L$.

Observe that for any pair $(L,\pi),\ \Omega(\pi')$ can be
transformed into $L$ by identifying some of its vertices. $\pi$ is
reflexive in $L$ if no identification is necessary.

Given $(L,\pi)$ the crossing function $c_\pi$, has its image
contained in the cycle space of $L$, $CS(L)$. (To get $CS(L)$ the
orientation of the edges of $L$ plays no role.) To show this
inclusion, observe that the edges between the two occurrences of
an edge $e$ in $\pi$ correspond to a path $\pi_e$ which is
reentrant if $e$ is traversed in opposite directions in $\pi$;
$\pi_e$ is made reentrant by the composition with the path
(tail$(e)$,$e$,head$(e)$), if $e$ is traversed in the same
direction in $\pi$.

\medskip\n{\bf (3.6) Lemma:} If $P$ is a map with one $v$-gon, then
$\{c_P(x)\!:\!x\!\in\!{\rm symb}(\overline{P})\}$
generates~$V^\perp$.
\medskip

\inicproof Consider the pair $(G_M,\pi_z)$ where $\pi_z$ is the
unique zigzag path in $M^t$. Consider also a fixed orientation of
the edges of $G_M$. Let $\pi'_z$ be induced by $(G_M,\pi_z)$. We
claim that $c_P=c_{\pi'_z},$ when $SQ(P)$ is identified with
$eG_M$. The cyclic sequence of squares met by the $v$-gon of $P$
and the cyclic sequence of the edges of $G_M$ given by $\pi'_z$
are equal. Thus, to prove the claim, we have only to show that an
edge $x$ of $G_M$ is traversed in the same direction twice in
$\pi_z$ iff the corresponding square $x$ is not in the imbalance
of $P$, i.e., if edge $x$ is an $s$-loop in $p^t$. But $x$ is an
$s$-loop in $P^t$ iff for $X$ arbitrary corner of $x$,
$(X,z_M(X))$ and $(v_M(X),f_M(X))$ are two ordered pairs of
consecutive vertices in the $v$-gon of $P$ traversed in one of the
two directions. Observe that $(x,z_M(X))$ and $(v_M(X),f_M(X))$
are such pairs iff edge $x$ of $G_M$ is traversed in the same
direction twice by the zigzag. Therefore, the claim is proved.

By Proposition 3.6.1, proved below, if $\pi_z$ is reflexive in
$G_M$, then the proposition follows. Therefore, assuming 3.6.1, to
conclude the proof it is enough to show that for arbitrary $M^t$
with just one zigzag path, $\pi_z$ is reflexive in $G_M$.

If $\pi_z$ is not reflexive in $G_M$, then, by definition,
$\Omega(\pi'_z)\ne G_M$. This means that $\Omega(\pi'_z)$ has more
vertices than $G_M$, and to obtain $G_M$ from $\Omega(\pi'_z)$
some identifications of vertices of $\Omega(\pi'_z)$ are
necessary. It follows that there are two darts, $x$ and $y$, which
satisfy $\mbox{head}(x)\mbox{tail}(y)$ in $G_M$ but not in
$\Omega(\pi'_z)$. We can choose $x$ and $y$ verifying this
condition and, moreover, the condition that $x$ and $y$ are
consecutive edges in the cyclic sequence of edges around the
vertex $\mbox{head}(x)=\mbox{tail}(y)$ in $M^t$. However, note
that such consecutive edges correspond to a pair of consecutive
edges in $\pi_z$. Whence, by definition of
$\mbox{head}(x)=\mbox{tail}(y)$, which contradicts the assumption.
Thus, we conclude that $\pi_z$ is reflexive in $G_M$. The proof of
the lemma is reduced to the proof of Proposition 3.6.1, which we
give next. \fimproof

The main idea in the proof the the following proposition was
obtained from Read and Rosenstiehl [21].

\medskip\n{\bf (3.6.1) Proposition:} If $\pi$ is a cyclic path in a digraph $L$,
such that $\pi$ uses each edge of $L$ twice, then
$c_{\pi'}\!:\!2^{eL}\!\to\!CS(L)$ is an epimorphism if and only if
$\pi$ is reflexive in~$L$.

\inicproof To prove the proposition it is enough to show that the
image of $c_\pi$, is the cycle space of $\Omega(\pi')$. If this is
proved, it follows that if $\pi$ is not reflexive in $L$, then
$c_\pi$, is not an epimorphism; also conversely, if $\pi$ is
reflexive in $L$, then $c_\pi$, is an epimorphism.

We construct $Q(\pi)'$ from $\pi'$ as follows. We start by drawing
an arbitrary dart representing an occurrence of a signed symbol in
$\pi'$. We proceed drawing darts in the order that they appear,
without lifting the pencil from the paper, as long as this is
possible. Each time that an edge occurs for the first time we must
draw it as a pendent edge (the final vertex must be a new one).
The rule hot to lift the pencil is not possible to obey when an
edge $e$ occurring for the second time is not incident to the last
vertex reached. At each such occurrence, we make a copy of the
graph drawn so far, denoting it by $L_e$. Next, we make the
necessary identification of the two vertices and proceed. At the
end we have a sequence of partial graphs
$$L_{e_1},L_{e_2},\dots,L_{e_k},\Omega(\pi'),$$
where edge $e_i$ forced the $i$-th identification. Suppose we have
$k$ identifications. If no identifications were made $(k=0)$, then
we would have $e+1$ vertices, where $e$ is the number of edges of
$\Omega(\pi')$. Each identification reduces the number of vertices
by one. Therefore, if $v$ is the number of vertices of
$\Omega(\pi')$ we have
$$e+1-k=v\qquad\mbox{or}\qquad k-e+1-v.$$

Observe that this value of $k$ is the dimension of $CS(L)$, since
$L$ is connected.

For every $i=1,2,\dots,k$ the two facts that follow are evident:

(a)\quad for $j<i\quad c_\pi,(e_j)$ is a cycle in $L_i$;

(b)\quad $c_\pi,(e_i)$\quad is not a cycle in $L_i$.

Fact (a) follows from the observation that a cycle remains a cycle
when identifications of vertices are performed. Fact (b) follows
from the equivalence: the set of edges occurring an odd number of
times in a path $\pi$ is a cycle iff $\pi$ is reentrant.

From (a) and (b) follow that each $c_\pi,(e_i)$ is not expressible
as $\mod\,2$ sum of a subset of $\{c_\pi,(e_j):j<i\}$. Hence,
$$c_\pi,(e_1),\dots,c_\pi,(e_k)$$
are independent cycles in $\Omega(\pi')$. Since $k$ is the
dimension of the cycle space of $\Omega(\pi')$, the $c_\pi(e_i)$
is form a basis for this space. The proof is complete. \fimproof

Observe that with the proof of this proposition every statement
made so far in Chapter 3 is established.

We conclude the chapter with an example of a class of arbitrarily
large Gauss codes which are realizable only in surfaces of
connectivity at most two.

Denote by $\overline{P}_n$ the Gauss code
$$(1,2,\dots,n,1,2,\dots,n).$$

\n{\bf (3.7) Proposition:} if $n$ is odd there is one realization
of $\overline{P}_n$ in the plane, and every other realization of
$\overline{P}_n$ is in the torus. If $n$ is even there is one
realization of $\overline{P}_n$ in the projective plane, and every
other realization of $\overline{P}_n$ is in the Klein bottle.
\medskip

\inicproof Suppose we are given an arbitrary $\{$black,white$\}$
partition, $\{B,W\}$, of ${\rm symb}(\overline{P}_n)$. Let us
denote by $M_n$ the phial of the map $P_n$ induced by the
2-colored Gauss code.

Suppose that $n$ is odd. In this case $i^2_{P_n}(x)$ is the whole
of ${\rm symb}(\overline{P}_n)$. (Recall that $y\in i^2_{P_n}(x)$
and $i_{P_n}(y)$ have an odd number of common elements.) As we
have
$$(c_{P^\sim}\circ c_{P_n})=i^2_{P_n}(x)+\sum_{y\in i_{P_n}(x)}(1+\gamma_x+\gamma_y)y,$$
where $\gamma=1$ iff $z\in B$, the following implications hold:

(i)\quad if $x\in B$, then $(c_{P^\sim_n}\circ c_{P_n}(x)=W$;

(ii)\quad if $x\in W$ then $(c_{P^\sim_n}\circ c_{P_n}(x)=B$;

Therefore, except in the case where one of $W$ and $B$ is empty,
we have $\dim[Im(c_{P^\sim_n}\circ c_P)]=2$. Since
$ODD(\overline{P})$ is empty, there are no $r$-circuits in
$G_{M_n}$, by Lemmas 3.4 and 3.5(a). By Theorem 3.3(a), it follows
that $S_{M_n}$ is the torus. In the case where one of $B$ and $W$
is empty we have (i) and (ii) implying
$$\dim[Im(c_{P^\sim_n}\circ c_{P_n})]=0.$$
Therefore, in this case $S_{M_n}$ is the plane, by Theorem 3.3(a).

Suppose that $n$ is even. Now $i^2_{P_n}(x)$ is empty for every
$x$ in ${\rm symb}(\overline{P}_n$). Hence, we have the
implications:

(i')\quad if $x\in B$, then $(c_{P^\sim_n}\circ c_{P_n}(x)=B$;

(ii')\quad if $x\in W$, then $(c_{P^\sim_n}\circ c_{P_n}(x)=W$.

It follows that if one of $B$ and $W$ is empty, then
$\dim[Im(c_{P^\sim_n}\circ c_{P_n})]=1$ and $S_{M_n}$ is the
projective plane. If both of $B$ and $W$ are non-empty, then
$\dim[Im(c_{P^\sim_n}\circ c_{P_n})]=2$, and $S_{M_n}$ is the
Klein bottle. \fimproof

\chapter{Transversals of Orientation -- Reversing Circuits}
\section*{4.A\quad Proliminary Concepts and Relation with Imbalances}
Given a finite collection $\Omega$ of subsets of a finite set $X$,
a {\it transversal} of $\Omega$ is a subset $T$ of $X$ such that
$T\cap Y$ is non-empty for every element $Y$ of $\Omega$. A {\it
minimal} transversal is a transversal not properly contained in
another transversal. A {\it minimum} transversal is a transversal
of minimal cardinality. As usual, the adjectives minimal and
minimum, as well as maximal and minimum, have this distinction,
when applied to some concept.

This chapter establishes a minimax relation (Theorem 4.5), which
includes in its proof an algorithm, polynomially bounded in the
``size'' of the map, to exhibit a minimum transversal of
$r$-circuits in arbitrary projective maps.

Our aim in this first section is to relate transversals of
$r$-circuits and imbalances (Lemma 4.1). Also we show (Proposition
4.2) that the problem of deciding whether a subset of $eG_M$ is a
minimum transversal of $r$-circuits in an arbitrary map $M$ is
$NP$-complete, in the sense introduced by Karp [13].

With the next lemma we switch from transversals of $r$-circuits to
imbalances, which is technically convenient.  Once mote, because
of the natural 1-1 correspondence between $eG_M$ and $SQ(M)$, we
do not distinguish explicitly their members or subsets.

\medskip\noi{\bf (4.1) Lemma:} For every map $M$, the minimal transversals of
$r$-circuits in $M^t$ are precisely the minimal imbalances in $M$.
\medskip

\inicproof Let $I$ be a minimal imbalance in $M$. Denote by $f_I$
the set of $f_M$-edges in the squares in $I$ and let $\{A,B\}$ be
the balancing partition inducing $I$. By the definition of
balancing partition and of induced imbalance, $f_I$ is precisely
the set of edges with both ends in $A$ or both ends in $B$. Hence,
$\delta A=eC_M\backslash f_I=\overline{f}_I$. Therefore, we
conclude the following equivalence: $I$ is a minimal imbalance in
$M$ iff $\overline{f}_I$ is a maximal coboundary in $C_M$.

Proposition 4.1.1, proved below, states that a subset of edges in
$eG$, for arbitrary map $G$, is a maximal coboundary in $G$ if and
only if its complement is a minimal transversal of odd circuits in
$G$. Thus, applying Proposition 4.1.1, we conclude that, the
complement $f_I$ of a maximal $\overline{f}_I$, as above, is a
minimal transversal of odd circuits in $C_M$. To conclude the
proof it is enough to observe that $I$ is a minimal transversal of
$r$-circuits in $M^t$ if and only if $f_I$ is a minimal
transversal of odd circuits in $C_M$. This equivalence follows
from the definition of $r$-circuit. Provided Proposition 4.1.1 is
established, the proof of the lemma is complete. \fimproof

\medskip\noi{\bf (4.1.1) Proposition:} For an arbitrary graph $G$, a subset of
edges in $eG$ is a maximal coboundary in $G$ if and only if its
complement is a minimal transversal of odd circuits in $G$.
\medskip

\inicproof It is evident that the complement in $eG$ of any
coboundary in $G$ is a transversal of odd circuits in $G$. Thus,
to prove the proposition, it is enough to show that a minimal
transversal of odd circuits is the complement in $eG$ of a
coboundary.

Let $T$ be a minimal transveral of odd circuits in $G$. Hence the
graph $G\backslash T$, obtained by deleting from $G$ all the edges
in $T$, is bipartite. Let $\{A,B\}$ be a bipartition of
$G\backslash T$. We claim that every edge in $T$ must have both
ends in $A$ or both ends in $B$. If edge $e\in T$ has one end in
$A$ and the other in $B$, then $T\backslash\{e\}$ is yet a
transversal of odd circuits and so, $T$ would not be a minimal
transversal. This contradiction proves the claim. Therefore, the
set of edges in the complement of $T$, each linking a vertex in
$A$ to a vertex in $B$, is the coboundary of $A$ in $G$.\fimproof

\medskip\noi{\bf (4.2) Proposition:} To determine whether a given subset of $eG_M$ is a
minimum transversal of $r$-circuits in a general map $M$ is an
$NP$-complete problem.
\medskip

\inicproof We show that the problem of determining whether a given
subset of $eG$ is a maximum coboundary in an arbitrary graph $G$
is reducible to the problem stated in the proposition. Since the
maximum coboundary problem is $NP$-complete, see Garey and Johnson
[6], this is enough to establish (4.2). By Proposition 4.1, it is
enough to reduce the maximum coboundary problem to the minimum
imbalance problem.

Assume we are given a graph $G$ for which we want to discover a
maximum coboundary. Let us embed $G$ arbitrarily in an orientable
surface such as to obtain a $t$-map $M^t$ with $G_M\cong G$.
Observe that $M$ is obtainable with an amount of work which is
bounded by a polynomial in the number of edges of $G$.
Specifically choose, for each vertex $x$ of $G$, an arbitrary
cyclic sequence of edges incident to $x$ to go around $x$ in
$M^t$. These cyclic sequences form a $v$-ordering for $M$. Take
the empty set to be the imbalance corresponding to this
$v$-ordering. By Proposition 3.0 a unique map $M$ is specified.
Since the set of imbalances in a map $M$ is a coset of the
coboundary space of $G_M$, and since $\phi$ is an imbalance in
$M$, we conclude that the set of imbalances in $M$ is precisely
the coboundary space of $G_M\cong G$. By Proposition 3.2.1,
$\subseteq SQ(M)$ is an imbalance in $M$ if and only if
$SQ(M)\backslash T$ is an imbalance in $M^\sim$. Therefore, to
find a maximum coboundary in $G$ is the same as to find a maximum
imbalance in $M$, which is the same as to find a minimum imbalance
in $M^\sim$. This concludes the proposition. \fimproof

\section*{4.B\quad Imbalances and Cycles in the Dual}
In this section we prove that the imbalances of a map are cycles
of a specific homology in the dual. We also prove that the set of
minimal imbalances in a projective $t$-map is precisely the set of
$r$-circuits in the dual $t$-map.

\medskip\noi{\bf (4.3) Theorem:} The set of imbalances in a map $M$ is a class of
homologous $(\mod\,2)$ cycles in the dual $D^t$ of $M^t$.
\medskip

\inicproof We already know that the set of imbalances in $M$ is a
coset of the coboundary space $V$ of $G_M$. It follows that the
symmetric difference of two arbitrary imbalances is a coboundary
in $G_M$, and hence a boundary in $D^t$.

Therefore, to conclude the proof it is sufficient to establish
that the set of squares in an imbalance $I$ corresponds to a set
of edges in $G_D$ which is a cycle in $G_D$.

Let $I$ be an arbitrary imbalance in a map\index{imbalance in a
map} $M$ induced by a balancing partition $\{A_I,B_I\}$. Let us
impose an orientation on the $v_M$- and the $a_M$-edges of $C_M$
such that every $v$-gon of $M$ becomes a directed polygon and such
that each $v_M$-dart has tail in $A_I$ and head in $B_I$.

Consider now the traversal of an arbitrary $f$-gon $K$, and an
$f_M$-edge $x$ in $K$. The two $a_M$-edges, preceding and
following $x$, are traversed one in the same direction and the
other in the opposite direction (according to the imposed
orientation of the $a_M$-edges) if and only if the square
containing $f_M$-edge $x$ is in the imbalance $I$. Observe that in
the complete traversal of $K$ the number of times that there is a
reversal of direction when we go from an $a_M$-edge to the next is
even. Therefore, the number of squares in $I$ met by the traversal
of $K$ is even. Thus, the number of squares in $I$ met once by the
traversal of $K$ is even. Observe that this set of squares
corresponds to a coboundary of the vertex $K$ in $G_D$ (the
vertices of $G_D$ are $f$-gons). Since $K$ is arbitrary, the
imbalance $I$ is orthogonal to every element in the coboundary
space $F$ of $G_D$. That is to say, I belongs to $F^\perp$ the
cycle space of $G_D$. This establishes the theorem. \fimproof

Let $M$ be a projective map. By Lemma 2.6 and the previous
theorem, the set of imbalances in $M$ is either the set $V$ of
boundaries in $D^t$ or a set of the form $x+V$, where $x$ is an
$r$-cycle in $D^t$. The set $V$ of coboundaries in $G_M$ is
discarded, since $V$ is the set of imbalances in $M$ iff $M$ is
orientable, which is not the case by Theorem 1.7. Therefore, we
obtain the following lemma:

\medskip\noi{\bf (4.4) Lemma:} The minimal imbalances in a projective map $M$
are precisely the $r$-circuits in $G_D$ where $D$ is the dual of
$M$.
\medskip

\inicproof We proved that the set of imbalances in $M$ is of the
form $X+V$, where $x$ is an $r$-cycle in $G_D$ and $V$ is the
coboundary space of $G_M$. The proof of the lemma follows, then,
by the definition of $r$-circuit as a minimal $r$-cycle. \fimproof

\section*{4.C\quad Main Theorem}
A collection $Q$ of subsets of $X$ is a {\it disjoint}
collection\index{{\it disjoint} collection} if the members of
$\Omega$ have no elements in common. The main theorem of the
chapter is the following:

\medskip\noi{\bf (4.5) Theorem:} For an arbitrary projective map $M$ with $G_M$ eulerian,
a minimum transversal of $r$-circuits in $M^t$ is  equal in
cardinality to a maximum disjoint collection of $r$-circuits in
$M^t$.

By Lemma 4.1, the minimal transversals of $r$-circuits of a map
$M$ are precisely the minimal imbalances of $M$. Hence, in the
statement of (4.5) we can replace ``transversal of $r$-circuits in
$M^t$'' by ``imbalances in M''.

As we show in the sequel, the hypothesis that $G_M$ is eulerian is
not essential, and the algorithmic proof of (4.5) provides us with
a polynomial algorithm to determine the size of a minimum
imbalance in general projective maps. As we showed in Proposition
4.2, the minimum imbalance problem is $NP$-complete. It is,
therefore, of interest to know exactly for which classes of maps
(4.5) holds. In our conclusions of Chapter 5 we present another
class of maps in which (4.5) holds.

By Lemma 4.4, the minimal imbalances in a projective map $M$ are
precisely the $r$-circuits of $G_D$, where $D$ is the dual of $M$.
This provides us with yet another equivalent form of statement
4.5: we can replace ``transversal of $r$-circuits in $M^t$'' by
``$r$-circuit in $D^t$'', where $D$ is the dual of $M$.

It is this form of (4.5), making use of $r$-circuits in the dual,
which is suitable for the proof that we present in Section 4.
Before we introduce, in the next section, the concept of
``reducer'', which is the most important concept in the proof of
(4.5) we give a corollary avoiding the hypothesis that $G_M$ is
eulerian.

Given a map $M$ let us denote by $M_d$ the map which corresponds
to the $t$-map $M^t$ in which each edge of $G_M$ is replaced by a
pair of edges forming a bounding digon in the surface $S_M$. Let
$D_d$ be the dual of $M_d$. It is evident that each circuit in
$G_{D_d}$ corresponds to a circuit in $G_D$ half as long in
length. In particular, each $r$-circuit in $D^t$ corresponds to an
$r$-circuit in $D^t_d$ half as long in length. Therefore, the
cardinality of a minimum $r$-circuit in $D^t_d$ is twice the
cardinality of a $r$-circuit in $D^t$. Consequently, the
cardinality of a minimum transversal of $r$-circuits in $M^t_d$ is
twice the cardinality of a minimum transversal of $r$-circuits in
$M^t$. This observation implies the following
\medskip

\medskip\noi{\bf (4.5) Corollary:} For an arbitrary projective map $M$, a minimum
transversal of $r$-circuits in $M^t$ has cardinality equal to half
the cardinality of a maximum collection of $r$-circuits in $M^t$
with the property that each edge of $G_M$ is used at most twice in
the collection.
\medskip

\inicproof The corollary follows from the application of Theorem
4.5 to the map $M_d$, describe above. \fimproof

In the next section we introduce the concept of ``reducers'' The
elimination of reducers is the basis for a polynomial algorithm to
establish the minimax results stated. In Figure 4.1 below, we show
the result of the application of this algorithm, illustrating the
Corollary of 4.5. Figure 4.1 presents a projective $t$-map $M^t$,
for which a minimum transversal of $r$-circuits has cardinality 5:
the set of five edges which appear into two pieces is clearly a
transversal of $r$-circuits in $M^t$. The proof that the set of
these five edges is a minimum transveral is established by
exhibiting ten $r$-circuits in which each edge is used at most
twice by the ten circuits.

$$\includegraphics[width=8cm,height=8cm]{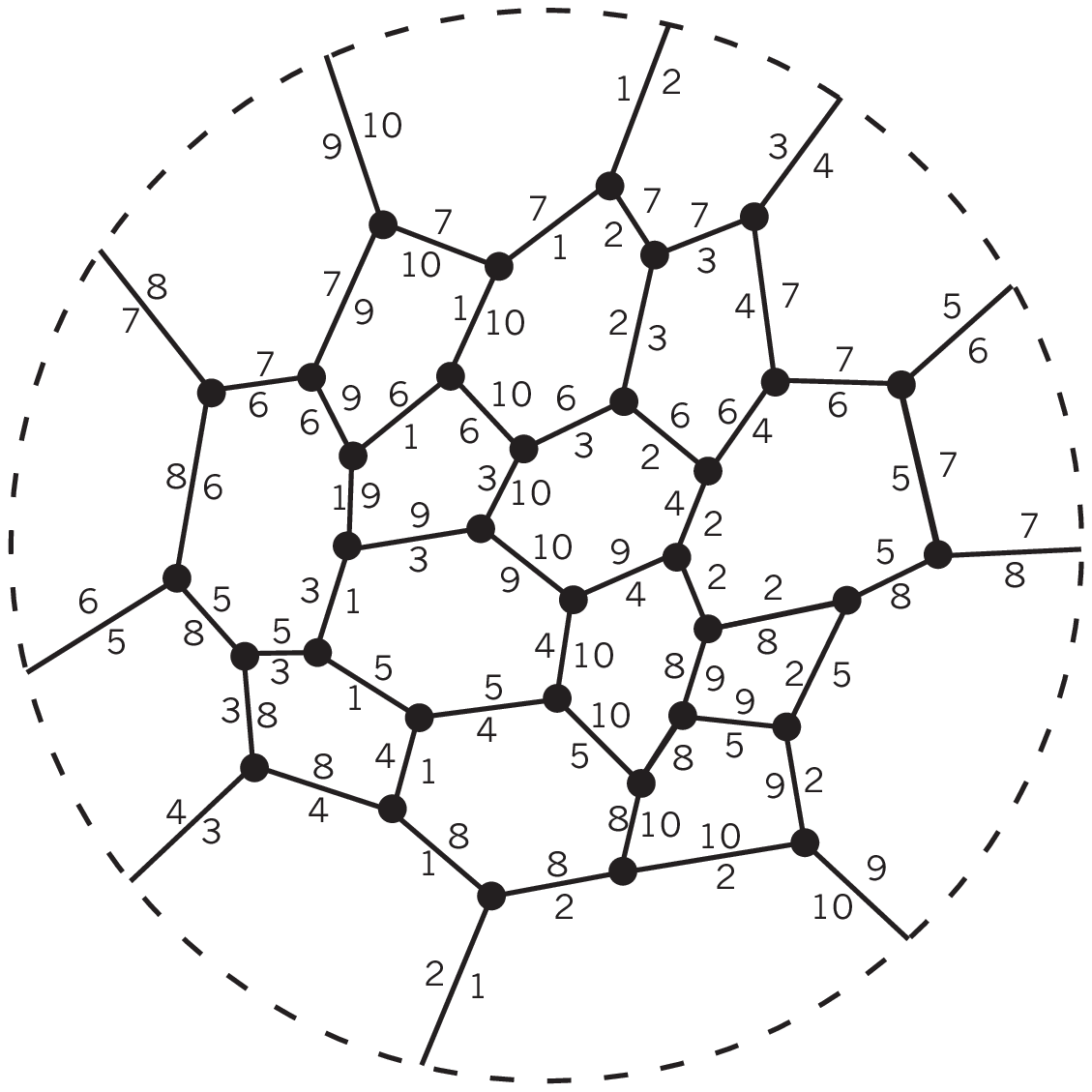}$$
\centerline{Figure 4.1}

\section*{4.D\quad Splitting Reducers}
A pre-requisite for the definition of a reducer is the notion of
contractibility. The definition that we give is suitable for our
purposes. This becomes apparent in Corollary 4.8.1.

If $R$ is a subset of faces of $M^t$, we denote by $eR$ the set of
edges of $G_M$ that appear (one or two times) in the facial paths
of the faces in $R$.

A subset $R$ of faces of $M^t$ is contractible if its boundary is
a circuit and if every circuit in $G_M$ included in $eR$ is a
boundary in $M^t$.

Observe that $R$ is contractible iff the union of the faces in $R$
together with the (topological) subgraph of $G_M$ generated by
$eR$ is homeomorphic to a closed disc.

We call the attention of the reader to the fact that for a subset
of faces $R$ whose boundary is a circuit it is not the case, in
general, that either $R$ or its complement $\overline{R}$ is
contractible. However, it is the case if $M$ is projective, as we
establish later in Corollary 4.8.1. In the Klein bottle we may
have a subset $R$ of faces whose boundary is a circuit but such
that neither $R$ nor $\overline{R}$ is contractible.

For a map $M$ in which $G_M$ is eulerian, we denote by $SP(M^t)$
the set of smooth paths of $M^t$.

For a map $M$ in which $G_M$ is eulerian, a non-empty subset $R$
of faces of $M^t$ is called a {\it reducer} for $M^t$ if $R$ is
contractible and the circuit which constitutes the boundary of $R$
is induced by at most two subpaths of members of $SP(M^t)$.

A {\it type 0 reducer} is a reducer whose boundary is induced by a
(complete) smooth path. Note that a type 0 reducer cannot be
minimal: if $R_0$ is a type 0 reducer there is a reducer R, (not
necessarily of type 0) such that $R$ is a subset of faces properly
contained in $R_0$.

A {\it type 1 reducer} is a reducer whose boundary is induced by a
unique proper reentrant subpath of a smooth path.

A {\it type 2 reducer} is a reducer whose boundary is induced, not
only by one but by two subpaths of smooth paths with the same
ends.

Observe that a reducer $R$ of type 1 which is also a minimal
reducer, in the sense that $R$ does not properly contain a subset
of faces $S$ which is also a reducer of type 0, 1 or 2, is
necessarily a singleton $\{f\}$ where $f$ is a face whose boundary
is a loop in $G_M$ (monovalent face).

We identify now an important property of any minimal reducer $R$
of type 2. Denote by $p$ and $q$ the two subpaths of members of
$SP(M^t)$ which induce the boundary of $R$. A {\it
segment}\index{\it segment} of $R$ is a maximal subpath of a
member of $SP(M^t)$ with at least one edge whose edges are
contained in $eR\backslash\partial R$, where $\partial R$, recall,
is the boundary of $R$. It is evident that $R$ has the following
property TRANS: for an arbitrary segment $s$, of $R$, one end of
$s$ is an interior vertex of $p$ and the other is an interior
vertex of $q$.

A {\it transreducer}\index{\it transreducer} is either a
monovalent face or else a type 2 reducer satisfying property
TRANS, described above. The transreducers are important for two
reasons: first, given a reducer $R$ it is algorithmically simple
to find a transreducer $S$ contained in $R$; second, property
TRANS in a reducer is enough for our purposes, which are described
in the sequel.

Figure 4.2 below shows examples of transreducers.

$$\includegraphics[width=3cm,height=3.5cm]{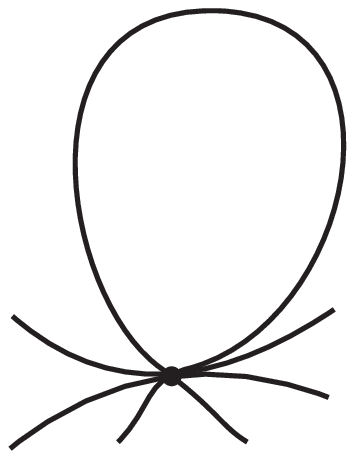}\qquad\qquad
\includegraphics[width=3cm,height=3.5cm]{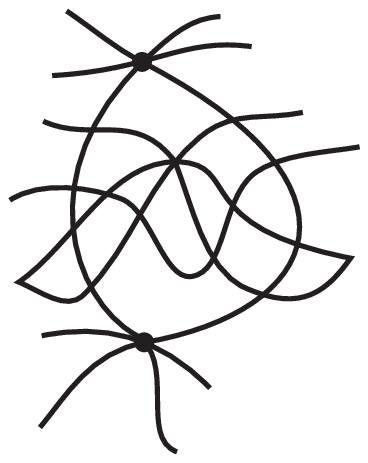}$$
\centerline{Figure 4.2}

An \index{\it angular point} in a reducer $R$ it a vertex which is
a common end of the path or paths which induce the boundary of
$R$. Hence, a reducer has zero, one or two angular points,
according as to whether it is of type 0, type 1 or type 2.

The splitting of a reducer\index{splitting of a reducer} $R$ (not
of type 0) is a local topological operation effected in $M^t$ near
an angular point $v$ of $R$, defined as follows. Denote by $e_1$
and $e_2$ the two edges in $\partial(R)$ which are incident to
$v$. Let $f_1$ and $f_2$ be the opposite edges to $e_1$ and $e_2$
at $v$, respectively. Separate the topological paths formed by
$f_1$ and $e_2$ from the one formed by $f_2$ and $e_1$, so that
both paths miss $v$. If $v$ has valency 4, then it is replaced by
two bivalent vertices after the splitting. See Figure 4.3.

$$\includegraphics[width=3cm,height=3.5cm]{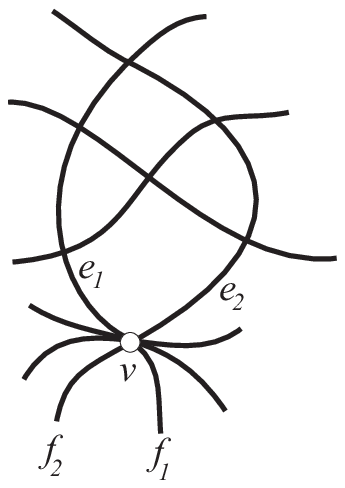}\quad
\includegraphics[width=3cm,height=3.5cm]{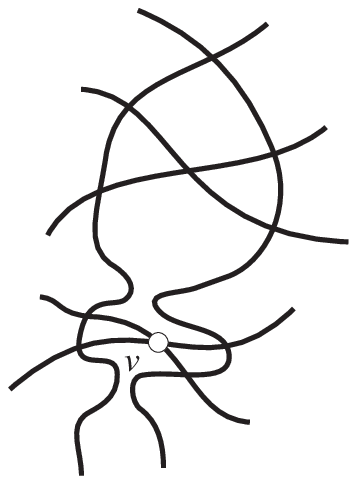}\quad
\includegraphics[width=3cm,height=3.5cm]{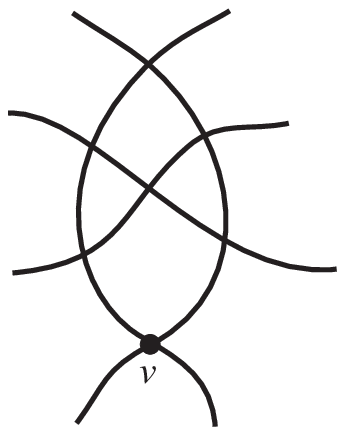}\quad
\includegraphics[width=3cm,height=3.5cm]{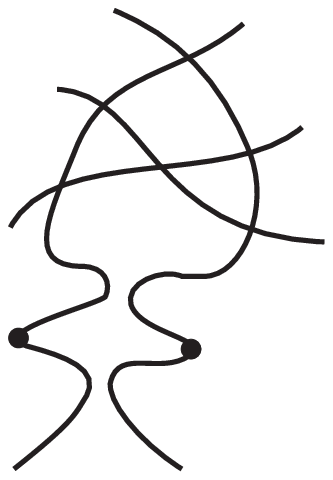}$$
\centerline{Figure 4.3}

Let us denote by $N^t$ the $t$-map obtained from $M^t$ by the
splitting of a transreducer. Let $E$ be the dual of $N$. It is
evident that the classes of homology in $CS(G_D)/V_M$ and in
$CS(G_E)/V_N$ are in natural 1-1 correspondence. We use this
correspondence to identify those classes and to state the
following important

\medskip\noi{\bf (4.6) Lemma:} For a map $M$ with $G_M$ eulerian, the splitting of a
transreducer does not change the cardinality of a minimum circuit
in a non-null homology class of $D^t$, where $D$ is the dual of
$M$.
\medskip

\inicproof To prove this lemma, consider a non-bounding circuit
$c$ in $E^t$, where $E^t$ is the dual of $N^t$, the $t$-map
obtained from $M^t$ by the splitting of a transreducer $R$. If $c$
is also a circuit in $G_D$, there is nothing to be proved.
Therefore, we may consider that the polygon corresponding to $c$
has among its vertices some which are new (triangular) faces in
$N^t$ introduced by the splitting. See Figure 4.4. (In the case
that the angular point of $R$ is 4-valent, then two faces of $M^t$
becomes a (new) face of $N^t$.) \vfill\eject
$$\includegraphics{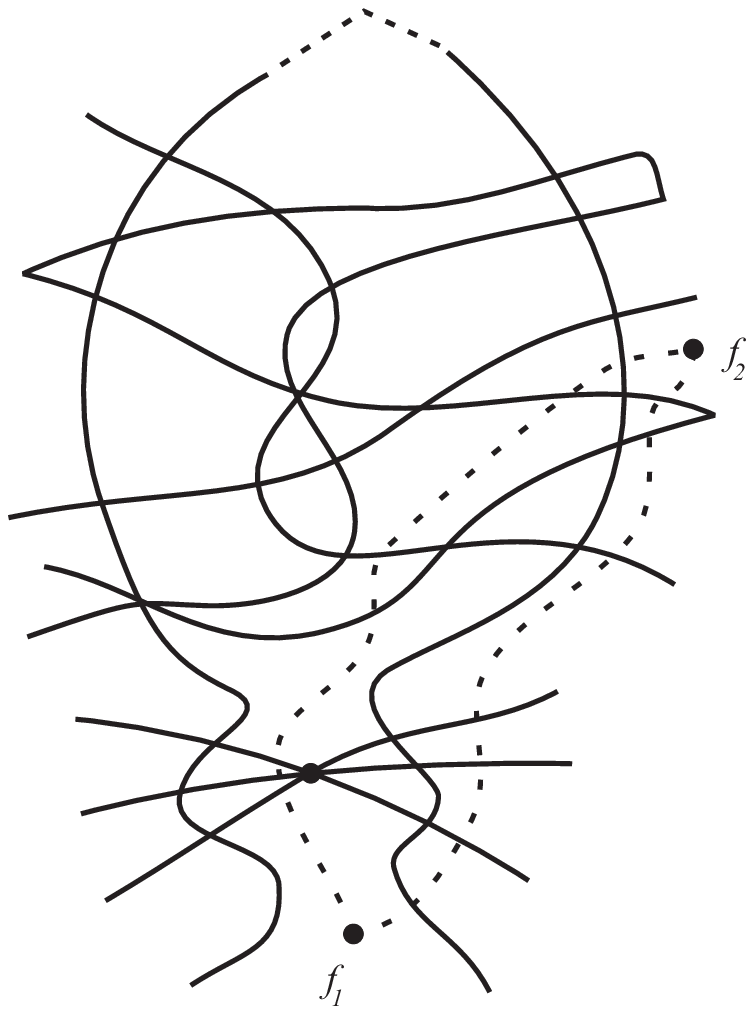}$$
\centerline{Figure 4.4}

Let $f_1$ be the face of $M^t$ incident to the two edges not in
the boundary of $R$ which are incident to the angular point of $R$
near which the splitting is effected. If $c$ is not a circuit in
the dual graph, $G_D$' then the polygon corresponding to $c$ has
an \index{arc} $a_1$ (a graph induced by a path without repeated
vertex-terms) which goes from vertex $f_1$ uses some new vertices
(new faces in $N^t$) and finishes in $f_2$ which is the first old
face not in $R$.

We claim that arc $a_1$ can be replaced by an arc $a_2$ which also
links $f_1$ to $f_2$ has the same number of edges as $a_1$, and
does not have among its vertices, either new faces of $N^t$ or
faces in $R$. See Figure 4.4. The new arc $a_2$ in $G_E$ starts at
$f_1$, crosses an edge of $G_N$ which corresponds to an edge in
the extension of the boundary of $R$ incident to the angular
point, crosses some edges of $G_N$ which have a vertex incident to
the boundary of $R$, until it reaches $f_2$. The crossing points
between the (topological) paths corresponding to $a_1$ and $a_2$
and the topological edges in $G_E$ are evidently in 1-1
correspondence with $eG_E$ and $eG_N$. Denote by $d_1$ the set of
edges of $a_i,\ i=1,2$. Observe that $d_1\subseteq c,\ d_2\cap
c=\phi$, and that $d_1+d_2$ is a boundary in $E^t$. Observe also
that $c'=(c\backslash d_1)+d_2$ is a circuit in $G_N$ as well as
in $G_E$, and it is homologous to c. By the fact that c is minimum
in its homology class, $|c'|\geq|c|$, that is $|d_2|\geq|d_1|$. It
is evident by construction, that $a_2$ links $f_1$ to $f_2$ and
does not have among its vertices either new faces of faces in $R$.
To conclude the claim, we show that the number of crossings of
$a_2$ with edges in $G_N$, i.e., $|d_2|$, is at most equal to the
number of crossings of $a_1$ with edges in $G_N$, i.e., $|d_1|$.
The first crossing of $a_2$ is compensated for by the last
crossing of $a_1$, with an edge in the boundary of $R$. A crossing
between $a_2$ and an edge in the extension of a segment of a
smooth path whose edges are contained in $eR$, is compensated for,
since $R$ satisfies TRANS, by an odd number of crossings between
$a_1$ and edges in the same segment. Thus $|d_2|\leq|d_1|$,
proving that $|c'|-|c|$. Thus, we have shown that the minimum
circuit in a non-null homology class of $D^t$ is at most equal in
cardinality to a minimum circuit in the corresponding homology
class of $E^t$. The converse inequality is trivial, because any
circuit in $G_D$ is also a circuit in $G_E$, (with the obvious
identifications of edges we have $eG_D\subseteq eG_E$), and the
homologies classes are maintained. This conclude the proof of the
lemma. \fimproof

It is possible that the minimum cycle of a given homology class is
not a circuit. For this reason we include the following
strengthening of Lemma 4.6.

\medskip\noi{\bf (4.7) Theorem:} For a map $M$ with $G_M$ eulerian, the splitting
of a transreducer does not change the cardinality of a minimum
cycle in a non-null homology class of $D^t$,  where $D$ is the
dual of $M$.
\medskip

\inicproof Let $N$ and $E$ be the pair of dual maps corresponding
to $M$ and $D$ after the splitting is performed. Let $c$ be a
minimum cycle in a non-null member $H$ of $CS(G_E)/V_N$. Let us
decompose $c$ into pairwise disjoint circuits
$$c=c_1+c_2+\cdots+c_k.$$
Each circuit $c_i$ is not a member of $V_N$, the coboundary space
of $G_N$; otherwise $c\backslash c_i$ would be a cycle in $H$ with
fewer edges that $c$, contradicting the assumption.

By Lemma 2.6, for each circuit $c_i$ in $G_E$ which is not a
circuit in $G_D$, we can obtain a circuit $c'_1$, homologous to
$c_i$, of the same cardinality as $c_i$, and which is a circuit in
$G_D$. Observe that, since $c'=(c\backslash c_i)+c'_i$ is a cycle
in $G_E$,  the minimality of $|c|$ in $H$ implies that $c_j\cap
c'_i=\phi$, for $j\ne i$. Hence, $|c'|=|c|$.  By replacing in this
way each $c_i$ which is not in $G_D$ by a corresponding $c'_i$ we
get a cycle $d$, homologous to $c$, of the same cardinality as
$c$, and which is a cycle in $G_D$. Thus, we have shown that the
minimum cycle in the homology class of $D^t$ identified with $H$
is equal in cardinality to a minimum cycle in $H$. The converse
inequality is trivial, because any cycle in $G_D$ is also a cycle
in $G_E$, and the homologies are maintained. This concludes the
proof of the theorem. \fimproof

With the above theorem we end this section. In the next section we
particularize it to projective maps so as to prove the main result
of the chapter, Theorem 4.5.

\section*{4.E\quad Proof of the Main Theorem}
In this section we give a proof of Theorem 4.5. Before, however,
we establish some lemmas. The following lemma is implicit in the
proof of Theorem 2.14. Now we include it as an explicit statement.

\medskip\noi{\bf (4.8) Lemma:} For every map $M$, the maximum number of pairwise
disjoint odd polygons in $C_M$ is at most $\xi(M)$.
\medskip

\inicproof By Proposition 2.6, applied to $t$-map $(C_M,S_M)$, the
quotient space
$$\frac{CS(C_M)}{\aw SQ(M)\fw+\aw vG(M)\fw+\aw fG(M)\fw}$$
has exactly $\xi(M)$ elements (classes of homology). By Theorem
1.8, no subset of circuits induced by pairwise disjoint odd
polygons is in the space generated by the set of edges in the
squares, $v$-gons, and $f$-gons of $M$. Hence, the set of edges in
each polygon of an arbitrary disjoint collection of odd polygons
in $C_M$ represents a distinct class of homology. \fimproof

\medskip\noi{\bf (4.8.1) Corollary:} Assume that $M$ is a projective map
and $R$ is a subset of faces in $M^t$ whose boundary is a circuit.
Then, either $R$ or $\overline{R}$, the complement of $R$ in the
set t of faces of $M^t$ , is contractible.
\medskip

\inicproof Assume indirectly that both $R$ and $\overline{R}$ are
not contractible. By definition, this implies that there exist
non-bounding circuits $S_1$ and $S_2$ (circuits which are not
boundaries), for which $S_1\subseteq eR,\ S_2\subseteq
e\overline{R}$. By Proposition 3.2.2., $S_1$ and $S_2$ are
$r$-circuits. Since $S_1$ and $S_2$ can met only along the common
boundary of $R$ and $\overline{R}$, it follows that we can find
disjoint odd (in cardinality) circuits, $S'_1$ and $S'_2$, in
$C_M$ such that $\psi^M_c(S'_i)=S_i$, $i=1,2$. But this
contradicts lemma 4.8, since it implies the existence of two
disjoint odd polygons in $C_M$, and $\xi(M)=1$. \fimproof

For disjoint cycles $S_1$ and $S_2$ in $G_M$ we say that $S_1$ and
$S_2$ {\it cross} if for every pair $(S'_1,S'_2)$ of cycles in
$C_M$ such that $\psi^M_c(S'_1)=S_1$ and $\psi^M_c(S'_2)=S_2,\
S'_1\cap S'_2\ne \phi$.

This formal definition of crossing is a set-theoretical
description of a crossing at a vertex of two edge-disjoint
polygons in $M^t$.

\medskip\noi{\bf (4.8.2) Corollary:} For any projective map $M$, if $S_1$ and $S_2$ are
disjoint $r$-cycles in $M^t$, then $S_1$ and $S_2$ cross.
\medskip

\inicproof Assume that $S_1$ and $S_2$ do not cross. Thus, there
exists disjoint odd cycles $S'_1$ and $S'_2$ in $C_M$ such that
$\psi^M_c(S'_i)=S_i,\ i=1,2$. Therefore, there exists disjoint odd
polygons $Q_1$ and $Q_2$ in $C_M$ such that $eQ_i\subseteq S'_i$.
But this contradicts lemma 4.8, because $\xi(M)=1$. \fimproof

A path $\pi$ is {\it simple} if no vertex appears twice as
vertex-term in $\pi$. A path $\pi$ is {\it semi-simple}\index{\it
semi-simple} if it is simple or if it is reentrant and the only
vertex which appears twice as vertex-term in $\pi$  is its origin
which is the same vertex as its terminus.

The technique used to prove 4.5 is to effect the splitting of
transreducers. To find a transreducer it is enough to be able to
locate reducers. To locate a reducer we first find a
pre-reducer\index{pre-reducer}, defined as follows.

A pre-reducer in a $t$-map $M^t$ with $G_M$ eulerian is an ordered
pair $(\pi_1,\pi_2)$ of subpaths of members of $SP(M^t)$, such
that $\pi_1$ is semi-simple, $\pi_2$ is non-degenerate, $\pi_1$
and $\pi_2$ have only their ends as common vertices, and
$e\pi_1\cup e\pi_2=e\pi_1+e\pi_2$ is an $e$-cycle in $M^t$.

We remark that in a pre-reducer $(\pi_1,\pi_2)$, $\pi_1$, may be
degenerate, being formed by just one vertex.

A {\it straight} pre-reducer $(\pi_1,\pi_2)$ is a pre-reducer in
which $\pi_2$ is semi-simple.

\medskip\noi{\bf (4.9) Proposition:} For a projective map $M$ with $G_M$ eulerian, if
$(\pi_1,\pi_2)$ is a straight pre-reducer\index{straight
pre-reducer}, then there is a reducer $R$ for $M^t$ such that
$\partial(R)=e\pi_1$ or $\partial(R)=e\pi_1+e\partial_2$.
\medskip

\inicproof Assume that $\pi_1$ and $\pi_2$ are both reentrant. By
definition of straight pre-reducer, $e\pi_i$ and $e\pi_2$ are both
$r$-circuits or both $s$-circuits. Observe, however, that $e\pi_1$
and $e\pi_2$ do not cross. Hence, by Corollary 4.8.2, they are
both $s$-circuits. By Proposition 3.2.2, $e\pi_1$ is a boundary of
a subset $X$ of faces. By Corollary 4.8.1, either $X$ or its
complement $\overline{X}$ is contractible. Let $R$ be a
contractible member of $\{X,\overline{X}\}$. By definition, it
follows that $R$ is a reducer, and since $\partial R=e\pi_1$, the
proposition follows.

Assume now that $\pi_1$ (and $\pi_2$) are not reentrant. By
definition of straight pre-reducer $e\pi_1+e\pi_2$ is an
$s$-circuit. Applying Proposition 3.2.2 and Corollary 4.8.1 as
before, with $e\pi_1+e\pi_2$ in the place of $e\pi_1$, we define a
reducer $R$ such that $\partial R=e\pi_1+e\pi_2$ concluding the
proof. \fimproof

\medskip\noi{\bf (4.10) Lemma:} For a projective map $M$ with $G_M$ eulerian, if $(\pi_1,\pi_2)$ is a
pre-reducer which is not straight, then there exists a pre-reducer
$(\pi_3,\pi_4)$ such that $e\pi_3\cup e\pi_4$ is a proper subset
of $e\pi_1\cup e\pi_2$.
\medskip

\inicproof Since $\pi_2$ is not semi-simple, we can write
$\pi_2=\alpha\circ\beta\circ\gamma$, where $e\beta$ is a circuit.
If $e\beta$ is an $s$-circuit in $M^t$, then take $\pi_4$ to be
the degenerate path whose only vertex is the origin of the
reentrant path $\beta$, and take $\pi_4$ to be $\beta$; it follows
that $(\pi_3,\pi_4)$ satisfies the condition of the lemma.
Therefore, we may assume that $e\beta$ is an $r$-circuit in $M^t$.
Since $e\pi_1+e\pi_2$ is an $s$-cycle and since
$e\pi_2=e\alpha+e\beta+e\gamma$, it follows that
$e\pi_1+e\alpha+e\gamma$ is an $r$-cycle. By Corollary 4.8.2,
$e\pi_i+e\alpha+e\gamma$ crosses $e\beta$.

It follows that at least one of the following vertex-terms must
exist:\jtt
\begin{itemize}
\item[(i)] a first vertex-term of $\alpha^{-1}$ other than its
origin which is also a vertex-term of $\beta$; \item[(ii)] a first
vertex-term of $\gamma$ other than its origin which is also a
vertex-term of $\beta$.
\end{itemize}
The situation is symmetric, up to taking inverse paths, and so we
may consider that (ii) holds. Denote by $v$ the vertex-term
described in (ii), and by $\gamma_1$ the proper subpath of
$\gamma$ that goes from its origin up to $v$. Let us write $\beta$
as $\beta_1\circ\beta_2$ where the terminus of $\beta_1$ and the
origin of $\beta_2$ coincide with $v$. One of $\beta_1$ and
$\beta_2$ may be degenerate. The reentrant paths
$\gamma_1,\beta_1$, and $\beta^{-1}_2$ all have the same ends.
Since $e\beta_1+e\beta_2$ is an $r$-cycle (indeed an $r$-circuit),
it follows that $e\beta_1+e\gamma_1$ is an $r$-cycle iff
$e\beta_2+e\gamma_1$ is an $s$-cycle. By the symmetry of the
situation, we may assume that $e\beta_1+e\gamma_1$ is an
$e$-cycle. We can now easily verify that $(\beta_1,\gamma_1)$ is a
pre-reducer. Since $e\beta$ is a circuit, $\beta_1$ is
semi-simple. By definition of $v,\gamma_1$ is non-degenerate and
$\beta_1$ and $\gamma_1$ have in common only their ends. Since
$e\beta_1+e\gamma_1$ is an $s$-cycle, it follows that
$(\beta_1,\gamma_1)$ is a pre-reducer. Put
$\pi_3=\beta_1,\pi_4=\gamma_1$, and the proof of the lemma is
complete. \fimproof

The following lemma is the last one before we procee to the proof
of Theorem 4.5.

\medskip\noi{\bf (4.11) Lemma:} For a projective map $M$ with $G_M$ eulerian, if
$\pi\in SP(M^t)$ is such that $e\pi$ is an $s$-cycle in $M^t$, or
if $\pi\in SP(M^t)$ is not semi-simple, or if two elements of
$SP(M^t)$ have more than one vertex-term in common, then there
exists a pre-reducer in $M^t$.
\medskip

\inicproof Suppose that $\pi\in SP(M^t)$ is such that $e\pi$ is an
$s$-cycle in $M^t$. Take $\pi_0$ to be a degenerate path whose
only vertex~term is a vertex-term of $\pi$. By definition,
$(\pi_0,\pi)$ is a pre-reducer in $M^t$.

Assume next that every member of $SP(M^t)$ induces an $r$-cycle
and that $\pi\in SP(M^t)$ is not semi-simple. Let us write $\pi$
as $\alpha\circ\beta$, where $e\beta$ is a circuit. (In particular
this implies $|e\beta|\geq1$.) If $e\beta$ is an $s$-circuit, then
consider $\delta$ to be the degenerate path whose only vertex is
the origin of $\beta$. The pair $(\delta,\beta)$ is a pre-reducer
in $M^t$. Thus, we may assume that $e\beta$ is an $r$-circuit. But
since $e\alpha+e\beta=e\pi$ is an $r$-cycle, it follows that
$e\alpha$ is an $s$-cycle. Take $\delta$ to be the degenerate path
whose only vertex-term is the origin of $\alpha$. The pair
$(\delta,\alpha)$ is a pre-reducer in $M^t$.

Finally, assume that for every $\pi\in SP(M^t)e\pi$ is an
$r$-cycle and $\pi$ is semi-simple. Suppose that $\pi_1$ and
$\pi_2$ have two distinct vertices $v$ and $w$ as common
vertex-terms. Write $\pi_1$ as $\alpha\circ\beta^{-1}$,\ $\pi_2$
as $\gamma\circ\delta^{-1}$, where $\alpha,\beta,\gamma$ and
$\delta$ all have origin $v$ and terminus $w$. If
$\alpha+\gamma^{-1}$ induces an $r$-cycle, then
$\alpha\circ\delta^{-1}$ must induce an $s$-cycle, because
$e\gamma+e\delta=e\pi_2$ is an $r$-cycle, and $e\alpha+e\gamma=
e\alpha+e\delta+e\pi_2$. Therefore, we may suppose, without loss
of generality, that $e\alpha+e\gamma$ is an $s$-cycle. (If not,
take $e\alpha+e\delta$ in place of $e\alpha+e\gamma$.) Let $u$ be
the first vertex-term of path $\alpha$ such that $u$ is also a
vertex-term of $\gamma$. Vertex-term $u$ must exist because $w$ is
common to both $\alpha$ and $\gamma$. Let $\alpha_1$ be the
subpath of $\alpha$ which goes from $v$ until $u$. Let $\gamma_1$
be the subpath of $\gamma$ which goes from $v$ until $u$. We want
to show that $(\alpha_1,\gamma_1)$ is a pre-reducer in $M^t$. By
definition of $u,\ \alpha_1$ is semi-simple, $|e\gamma_1|\geq1$,
and paths $\alpha_1$ and $\gamma_1$ have only their ends as common
vertices. To conclude that $(\alpha_1,\gamma_1)$ is a pre-reducer,
it is enough to show that $e\alpha_1+e\gamma_1$, is an $s$-cycle.
Suppose that $e\alpha_1+e\gamma_1$, is an $r$-cycle. Let
$\alpha_2$ and $\gamma_2$ be defined by the equations
$\alpha=\alpha_1\circ \alpha_2,\ \gamma=\gamma_1\circ\gamma_2$.
Observe that the cycle $e\alpha_2+e\gamma_2$ is an $r$-cycle,
since $e\alpha_1+e\gamma_1$ is an $r$-cycle and $e\alpha+e\gamma$
is an $s$-cycle. Observe also that, again by the definition of
$u,\ e\alpha_1+e\gamma_1$ does not cross $e\alpha_2+e\gamma_2$.
But this contradicts Corollary 4.8.2. Therefore,
$e\alpha_1+e\gamma_1$, is an $s$-cycle and so,
$(\alpha_1,\gamma_1)$ is a pre-reducer, concluding the proof of
the lemma. \fimproof

Now we are in a position to prove Theorem 4.5 in the following
equivalent form, involving $r$-circuits in the dual:

\noi{\bf (4.5 restated) Theorem} For an arbitrary projective map
$M$ with $G_M$ eulerian, a minimum $r$-circuit in $D^t$ is equal
in cardinality to a maximum disjoint collection of $r$-circuits in
$M^t$.

\inicproof The proof of the theorem is divided into two parts. In
the first part we show how to reduce the proof of the general case
to the special case in which $M^t$ is a system of projective
lines\index{system of projective lines}, which we define nest. In
the second part we give a proof for this special case.

A projective $t$-map $M^t$ with $G_M$ eulerian is called a {\it
system of projective lines} if every member of $SP(M^t)$ induces
an $r$-circuit and for two arbitrary members, $\pi_1$ and $\pi_2$
of $SP(M^t)$ the polygons defined by $\pi_1$ and $\pi_2$ have
precisely one vertex in common.

Now we show how to reduce the Proof of 4.5 from an arbitrary $M^t$
with $G_M$ eulerian to a system of projective lines. Suppose we
are given a projective $t$-map $M^t$ with $G_M$ eulerian. If there
exists a transreducer $R$ for $M^t$, then Lemma 4.6 shows that the
splitting of $R$ does not change the cardinality of a minimum
$r$-circuit in the dual. Moreover, a disjoint collection of
$r$-circuits in the $t$-map after splitting corresponds to a
disjoint collection of $r$-circuits in the original $t$-map before
splitting. Therefore, a proof of the theorem for the split $t$-map
implies a proof for the original $t$-map.

We now introduce a measure relative to which the $t$-map obtained
by splitting a transreducer is a smaller $t$-map. For a $t$-map
$M^t$, denote by $\mu(M^t)$ the ordered pair of non-negative
integers $(\lambda(M^t),\ f(M^t)$), where,
$$\lambda(M^t)=\sum_{v\in VG_M}\max\{0,{\rm val}_{G_M}\{v\}-4\},$$
and $f(M^t)$ is the number of faces of $M^t$. We say that
$\mu(N^t)<\mu(M^t)$ if $\lambda(N^t)<\lambda(M^t)$ or if
$\lambda(N^t)=\lambda(M^t)$ and $f(N^t)<f(M^t)$. (Lexicographic
order.)

Observe that if $N^t$ is obtained from $M^t$ by splitting a
reducer, then $\mu(N^t)<\mu(M^t)$ this follows because if $v$, the
angular point near which the splitting is performed, has valency
greater than 4, then $\lambda(N^t)<\lambda(M^t)$; if $v$ has
valency 4, then $\lambda(Nt)=^\lambda(M^t)$, and
$f(N^t)=f(M^t)-1$. This implies that the number of successive
splittings of transreducers that we can perform in a $t$-map $M^t$
before reaching a $t$-map, $N^t$ free of transreducers is at most
$2|eG_M|+f(M^t)$.

The above analysis implies that we may restrict ourselves to prove
4.5 for maps free of transreducers. We now show that a map free of
transreducers is system of projective lines. If $M^t$ is free of
transreducers, $M^t$ is free of reducers, because a reducer always
contains a transreducer. Since a straight pre-reducer implies a
reducer, by Proposition 4.9, it follows that $M^t$ has no straight
pre-reducer. By Lemma 4.10 a pre-reducer which is not straight
implies a ``smaller'' pre-reducer. Hence, any pre-reducer implies
a straight pre-reducer. As a consequence, $M^t$ is free of
pre-reducers. By lemma 4.11, this implies that:\jtt
\begin{itemize}
\item[(i)] for $\pi\in SP(M^t)$, $e\pi$ is an $r$-cycle;
\item[(ii)] if $\pi\in SP(M^t)$, then $\pi$ is semi-simple;
\item[(iii)] if $\pi_1$ i and $\pi_2$ are distinct members of
$SP(M^t)$, then $\pi_1$ and $\pi_2$ have at most one vertex-term
in common.
\end{itemize}

Items (i) and (ii) imply that if $\pi\in SP(M^t)$, then $e\pi$ is
an $r$-circuit. Item (iii) implies that the two polygons induced
by two arbitrary members of $SP(M^t)$ have at most and vertex in
common. By Corollary 4.8.2, two $r$-circuits must cross.
Therefore, every pair of polygons as above has precisely one
vertex in common. This proves that if $M^t$ is free of
transreducers, then $M^t$ must be a system of projective lines.
The first part of the proof is complete.

Now we conclude the theorem for a system of projective lines
$M^t$. Explicitly, we exhibit a disjoint collection $\Omega$ of
$r$-circuits in $M^t$ and an $r$-circuit $R$ in $D^t$ such that
$|\Omega|=|R|$. The collection $\Omega$ is simply $\{e\pi: \pi\in
SP(M^t$)\}. Let $n=|\Omega|$. To obtain $R$, let $p$ be the closed
topological path corresponding to an arbitrary member of $SP(M^t)$
and a take point $P$ near $p$ in $S_M$. Refer to Figure 4.5 for
this description. Proceed from $P$ parallel to $p$ until the point
$Q$ opposite to $P$ relative to $p$ is reached. Close the
topological path $p'$ which we are describing by linking $Q$ to
$P$, crossing $p$. The closed topological path $p'$ which we
described corresponds to an $r$-circuit $R$ in $D^t$. The
crossings of $p'$ with the points in $G_M$ are in 1-1
correspondence with the edges in $R$. Each crossing of $p'$, with
the exception of the last described, corresponds to a crossing of
$p$ with another member of $SP(M^t)$; hence, they are $n-1$ in
number. Therefore, $|R|=n=|\Omega|$, proving the Theorem.
\fimproof

$$\includegraphics[width=8cm,height=8cm]{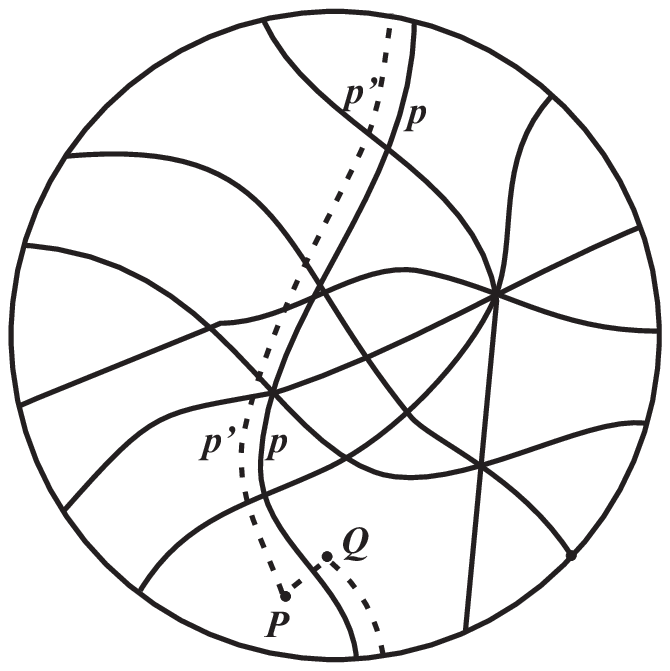}$$
\centerline{Figure 4.5}
\bigskip

To conclude this section we show with an example that even for the
torus the maps free of reducers can be quite complicated.
Consider, for example, the toroidal $t$-map free of reducers shown
in Figure 4.6. There are three smooth paths as pointed by
different thickness on the edges of the illustration. Any two such
smooth paths meat each other seven times, illustrating the point
that straight pre-reducers do not imply reducers in general.

$$\includegraphics[width=8cm,height=8cm]{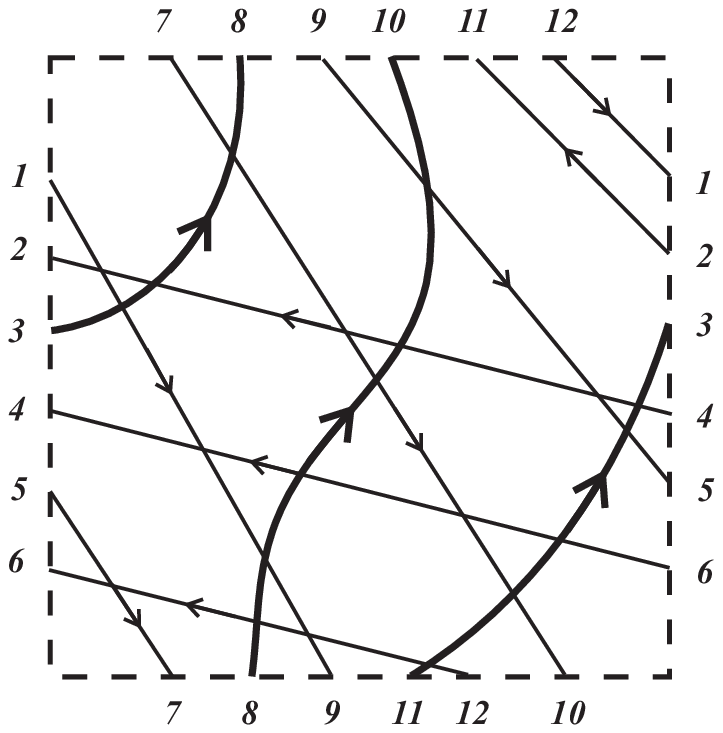}$$
\centerline{Figure 4.6}

\section*{4.F\quad Algorithmic Considerations}
In this section we present an algorithm, PROJMINMAX, which has as
its input a projective map $M_0$ with $G_{M_0}$ eulerian and as
its output a pair $(\Omega_0,R_0)$, defined as follows:\jtt
\begin{itemize}
\item[(i)] $\Omega_0$ is a collection of cyclic paths in
$G_{M_0}$, such that the corresponding cycles in $M^t$ form a
disjoint collection of $r$-circuits. \item[(ii)] $R_0$ is an
$r$-circuit in $G_{D_0}$, where $D_0$ is the dual of $M_0$.
\end{itemize}
Moreover, the pair $(\Omega_0,R_0)$ satisfies $|\Omega_0|=|R_0|$.
This is a pair whose existence is guaranteed in theorem 4.5.

We also show that the amount of space (space complexity) and the
number of basic steps (time complexity) of PROJMINMAX are bounded
by polynomials $n$, the ``size'' of the map, where
$n=|eG_{M_0}|=|SQ(M_0)|$.

We say that an algorithm A has {\it time (space) complexity}
$0(n^k)$, if there exists a constant $C$ such that the number of
basic steps (bytes of memory required) in the execution of $A$ to
solve a problem whose input has size $n$ is at most $Cn^k$, except
for finitely many values of $n$.

The algorithm that we present is suggested by the theory for
finding transreducers and by the inductive step consisting in
splitting a transreducer, which form the basis for the proof of
Theorem 4.5. However, in the algorithmic realization, we do not
actually split transreducers.

In order to avoid splittings in PROJMINMAX we introduce the
following $\Omega$-parametrized concepts, where $\Omega$ is an
edge-disjoint collection of cyclic paths in $G_M$ whose edge-sets
partition $eG_M$, for map $M$ with $G_M$ eulerian. The following
concepts generalize the ones already described in which $\Omega$
is defined by the smooth paths of $M^t$.

An {\it $\Omega$-reducer} for $M^t$ is a subset of faces $R$ such
that $R$ is induced by at most two subpaths of members of $\Omega$
and such that $eR$ contains only bounding circuits.

If $R$ is an $\Omega$-reducer for $M^t$, then a {\it segment} of
$R$ is a maximal subpath of a member of $\Omega$ with at least one
edge, whose edges are contained in $eR\partial R$.

An {\it angular point} in an $\Omega$-reducer $R$ is a vertex in
which the extension of the path or paths forming the boundary of
$R$ cross.

Property {\it TRANS} in a $\Omega$-reducer $R$ with two angular
points and boundary induced by two subpaths $p$ and $q$ of members
of $\Omega$ is defined as follows: for every segment $s$ of $R$,
the extension of a in $\Omega$ crosses $p$ and not $q$ at one end
of $s$ and $q$ and not $p$ at the other end.

An {\it $\Omega$-transreducer} is an $\Omega$-reducer $R$ with two
angular points for which property TRANS holds or else an
$\Omega$-reducer $R$ with only one angular point that is formed by
a single face.

An {\it $\Omega$-pre-reducer} is an ordered pair $(\pi_1,\pi_2)$
of sub-paths of members of $\Omega$ such that $\pi_1$ does not
self-cross, $\pi_2$ is non-degenerate, $\pi_1$ and $\pi_2$ cross
only at their ends (when extended), and $e\pi_1\cup
e\pi_2=e\pi_1+e\pi_2$ is an $s$-cycle in $M^t$. An
$\Omega$-pre-reducer is {\it straight} if $\pi_2$ does not
self-cross.

The proof of the three following ``$\Omega$-extensions'' are
repetitions of the ones given in 4.9, 4.10, and 4.11, where
$SP(M^t)$ takes the role of $\Omega$.

\noi{\bf(4.9)* (Extension of 4.9)} For projective map $M$ with
$G_M$ eulerian and collection $\Omega$ of cyclic paths whose
edge-sets partitions $eG_M$, if $(\pi_1,\pi_2)$ is a straight
$\Omega$-pre-reducer, then there is an $\Omega$-reducer $R$ for
$M^t$ such that $\partial R=e\pi_1$ or $\partial R=e\pi_1+e\pi_2$.
\fimproof

\noi{\bf(4.10)* (Extension of 4.10)} For $(M,\Omega)$ as above, if
$(\pi_1,\pi_2)$ is an $\Omega$-pre-reducer which is not straight,
then there exists an $\Omega$-pre-reducer $(\pi_3,\pi_4)$ such
that $e\pi_3\cup e\pi_4$ is a proper subset of $e\pi_1\cup
e\pi_2$\fimproof

\noi{\bf(4.11)* (Extension of 4.11)} For a pair $(M,\Omega)$ as
above, if $\pi\in\Omega$ is such that $e\pi$ is an $s$-cycle in
$M^t$, or if $\pi\in\Omega$ self-crosses, or if two elements of
$\Omega$ cross more than once, then there exists an
$\Omega$-pre-reducer in $M^t$. \fimproof

In Figure 4.7 we give the flow-graph for PROJMINMAX. It operates
on the map $M=(C_M,V_M,f_M)$ rather than on $M^t$. The collection
$\Omega$ at every stage of the algorithm is a collection of cyclic
paths free of repeated edges whose edge-sets partition $eG_{M_0}$.
We observe that a convenient way to keep track of $\Omega$ is, for
every $v$-gon at $M_0$, to specify a pairing in its $v_M$-edges.
This pairing is possible, since $G_{M_0}$ is eulerian. Given
$\Omega$, we define this pairing as follows. For each subsequence
of a member of $\Omega$ of the form $(e_1,V,e_2)$, where $e_1$ and
$e_2$ are edges of $G_M$ and $v$ is a vertex of $G_M$, denote by
$E_1$ and $E_2$ the $M$-squares that correspond to edges $e_1$ and
$e_2$; denote by $v$ the $v$-gon of $M$ which corresponds to the
vertex $v$. We say that the $v_M$-edge which is in $E_1$ and in
$V$ and the $v_M$-edge which is in $E_2$ and in $V$ are {\it
$\Omega$-mates}. Observe that the complete collection of pairs
consisting of $\Omega$-mates is enough to specify $\Omega$ itself.

Splitting a transreducer is replaced by the suitable change of two
pairs of $\Omega$-mates: $\{a,b\},\ \{c.d\}$ become $\{a,c\},\
\{b,d\}$ or $\{a,d\},\ \{c,d\}$. Observe that, in this way,
$\Omega$ changes but not the map $M$, as would be the case if
splittings were performed.

$$\includegraphics{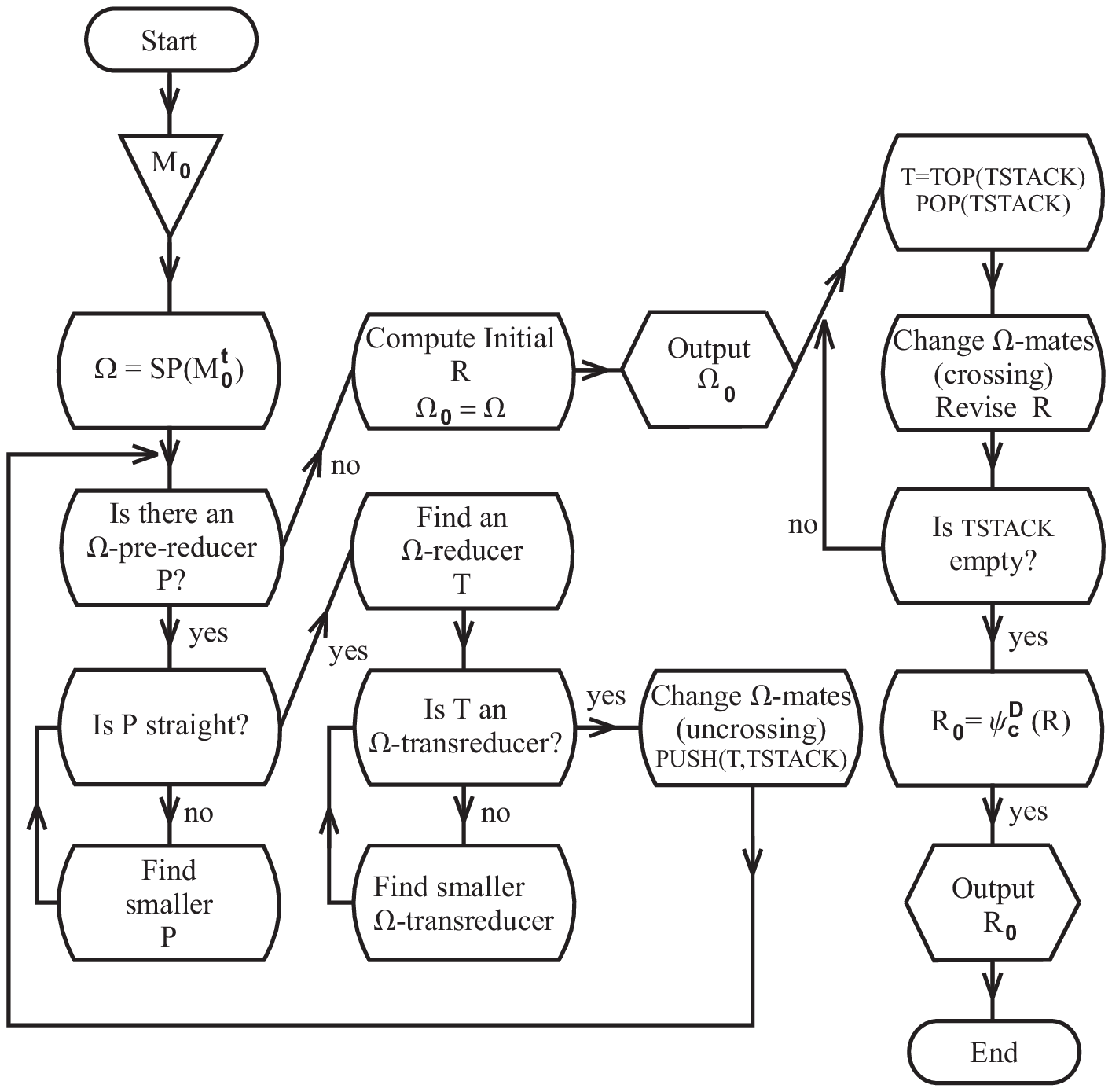}$$%%[width=3cm,height=3.5cm]
\centerline{Figure 4.7: {\em Flow-graph for PROJMINMAX}}
\bigskip

In the proof of 4.5 restated, we show that the number of
splittings is bounded by $2|eG_M|+f(M_0)\leq 0(n)$, where $f(M_0)$
is the number of $f$-gons in $M_0$. Modifying $\Omega$-mates is a
suitable way to perform splittings; it follows that $0(n)$ is the
order of the number of times that we answer the question ``Is
there an $\Omega$-pre-reducer P?''. See Figure 4.7.

For each positive answer of the above question, we answer in the
negative at most $0(n)$ times the following question: is $P$
straight? This follows from (4.10)*, and from the fact that an
$\Omega$-pre-reducer can get smaller at most $n$ times.

For each positive answer to the first question in Figure 4.7, only
once is the second question to be answered in the positive. For
each such event, the third question ``Is $T$ an
$\Omega$-transreducer?'' is answered at most $0(n)$ times. This is
so because a test to decide whether TRANS is satisfied ends in a
negative answer only if a smaller $Omega$-reducer is found. Hence,
this can happen at most $n$ times.

From the above analysis we conclude that the second and third
questions are answered in the execution of PROJMINMAX at most
$0(n^2)$ times and that the first question is answered at most
$0(n)$ times.

We use (4.11)* to bound the time complexity for an answer to the
first question by $0(n^2)$. The specific recipe on how to decide
whether there exists an $\Omega$-pre-reducer is given by an
adaption of the proof of (4.11) for the corresponding
$\Omega$-extension (4.11)*. The bound $0(n^2)$ is sufficient for
the decision.

An answer for the second and third questions can be obtained in
$0(n)$ basic steps. To find a smaller $\Omega$-pre-reducer, refer
to (4.10)*, modelling its proof according to the proof of (4.10).
This procedure is seen to take at most $0(n)$ steps. To find an
$\Omega$-reducer from a straight $\Omega$-pre-reducer it also
takes $0(n)$ steps. The same is true for finding a smaller
$\Omega$-reducer, in the case that the $\Omega$-reducer at hand
does not satisfy TRANS.

When an $\Omega$-transreducer $T$ is identified, we suitably
change the two pairs of $\Omega$-mated $v_M$-edges, an operation
which corresponds to split $T$. This is followed by placing $T$
(as a set of $f$-gons, for instance) in the top of a stack called
TSTACK. This structure is important for the revision of $R$
(explained below) in the second phase of the algorithm. The time
complexity of the above two operations is each at most $0(n)$. The
space required to keep all the structures so far is dominated by
the space to keep TSTACK, which takes at most $0(n^2)$ bytes of
memory.

This analysis implies that the time complexity for the first part
of PROJMINMAX, which terminates with a negative answer to the
first question, is $0(n^3)$.

The symbol $R$ represents throughout a circuit in $C_M$. Given
$\Omega$, the {\it $\Omega$-weight} of $R$ is defined as the
number of $\Omega$-mated pairs with the property that exactly one
member of the pair is a $v_M$-edge in $R$. The Initial $R$ is any
circuit in $c_M$ with the property that $\psi^M_c(R)$ is the
circuit induced by an arbitrary member of $\Omega_0$ and such that
its $\Omega_0$-weight is $|\Omega_0|$. See Figure 4.5. When we
change $\Omega$ (crossing), the $\Omega$-weight of $R$ can
increase by 2. The purpose of the revision of $R$ is to find
another $R$ with $\Omega$-weight equal to $|\Omega_0|$. Given the
$\Omega$-transreducer, this revision is a slight adaptation in
terms of $C_M$ of the Proof of Lemma 4.6. Refer to Figure 4.4. The
$v_M$-edges of $C_M$ correspond precisely to the crossings
presented in this figure. The question asking whether TSTACK is
empty is made at most $0(n)$ times. Instruction $T={\rm
TOP(TSTACK)}$ is self-explanatory; instruction POP(TSTACK) removes
the topmost element of TSTACK. Each individual procedure of the
second phase has complexity $0(n)$. Therefore, the overall
execution of the second phase of PROJMINMAX has complexity
$0(n^2)$. The space complexity of this phase is also $0(n^2)$
since TSTACK requires more memory than everything else.

When TSTACK is empty, we have that $\Omega=SP(M^t)$, and that the
$\Omega$-weight of the current $R$ is equal to $|\Omega_0|$. For
every pair of $f_M$-edges in $R$ incident to a $v$-gon, we adjust
$R$ so as to use the shorter of the two possibilities to go around
the $v$-gon. This does not change the $\Omega$-weight of $R$.
Observe that, since $\Omega$ is now $SP(M^t_0)$, the
$\Omega$-weight of $R$, which is $|\Omega_0|$, is equal to the
number of $v_M$-edges in $R$. But this is the cardinality of
$R_0=\psi^D_c(R)$. Therefore, $|\Omega_0|=|R_0|$. Observe that the
parity of $R$ is not changed in the revisions. Hence, $R_0$ is an
$r$-cycle (indeed an $r$-circuit) in $G_{D_0}$.

By comparing the two phases of PROJMINMAX we see that its overall
time complexity is $0(n^3)$ and its overall space complexity is
$0(n^2)$.

\chapter{Final Comments}
\section*{5.A\quad Remarks on Chapter 1}
In Chapter 1 we describe a graph model for the study of graphs
embedded in surfaces. We show in this thesis that, with this
language, certain problems of topological flavor can be formulated
in combinatorial terms. The principal advantage of the graph
approach is that a 2-dimensional structure is naturally encoded by
a colored 1-dimensional structure. An important point in favor of
this methodology is that the corresponding concept to topological
orientability is the simpler concept of parity; in consequence,
orientation-reversing circuit corresponds to odd circuit. The
graph definition also provides a geometric-algebraic bridge via
Cayley color graphs and Tutte's axioms for a map [25]. Elementary
properties of permutations and the definitions introduced yield a
simple combinatorial proof that the Euler characteristic of an
orientable surface is even; see Theorem 1.7. In Theorem 1.8 we
prove a fact about cubic 3-edge-colored graphs which corresponds
topologically to the property that no union of a disjoint
collection of polygons in $G_M$ whose induced circuits are
$r$-circuits in $M^t$ forms the boundary of a subset of faces in
$M^t$. The model used also leads to the symmetric theory in
vertices, faces, and zigzags which motivates the definition of the
phial map and of the anti-map, which are maps analogous to the
well-known dual map. The dual and phial maps form the support for
the extension of planar duality presented in Chapter 2. To
conclude this section, we describe two generalizations of the
graph model for a map which merit investigation:

(I) In the definition of a map, we drop the restriction that two
of the colors must induce squares. The faithful embedding of a
3-edge-colored cubic graph, with the colors ordered is a natural
embedding of a hypergraph: the hyperedges are the polygons induced
by the two first colors. They hypervertices are the polygons
induced by the two last colors. Each polygon induced by the first
and last color is a ``face'' which, as in the case of graphs, is
not an intrinsic concept, because it depends on the particular
embedding.

(II) We can try to move one step further in dimension. Consider a
4-regular 4-edge-colored graph, in which the faithful embedding of
each component induced by any three of the four colors is in a
2-sphere. We can ``fill with 3-space'' each one of these
2-spheres, obtaining solid 2-balls. Observe that each bicolored
polygon, ``filled with 2-space'', occurs twice as 2-faces in the
collection of 2-spheres. Thus, we have a natural method for the
identification of the boundaries of the 2-balls, thus forming a
3-manifold without boundary. Hence, we can talk about faithful
embeddings for the members of the above class of 4-regular,
4-edge-colored graphs, which we call for short {\it $4$-graphs}.
By, analogy let us call a cubic 3-edge-colored graph a {\it
$3$-graph}.

The Proof of Theorem 2.0.2 implies that the faithful embedding of
a 3-graph is in a 2-sphere iff the cycle space of the 3-graph is
generated by its bicolored polygons. An important open question is
the following: are there necessary and sufficient
graph-theoretical conditions for the faithful embedding of a
4-graph to be in a 3-sphere?

\section*{5.B\quad Remarks on Chapter 2}
In chapter 2 we introduce rich maps and graphs, as a
generalization of planar duality. Loosely speaking, rich graphs
are the graphs which have two abstract duals. Planar graphs are
the class of graphs whose members have one abstract dual. Since
rich graphs are a natural abstract extension of planar graphs,
perhaps some important theorems for the latter extend to the
former.

At present, however, more needs to be understood about the nature
of rich graphs, perhaps along the topological lines of Theorem
2.9. To us, the more important questions about rich-graphs are the
following:\jtt
\begin{list}{}{\parsep-2pt}
\item[(i)] Are all rich graphs map-rich? \item[(ii)] is it
possible to strengthen Theorem 2.9 so as to provide an
algorithmically good characterization of non-map-rich graphs?
\end{list}

To conclude this section we state, without proof, the following
graph-theoretic characterization of imbalances in rich maps.

\n{\bf Theorem} For a rich map $M$, a subset $T\subseteq SQ(M)$ is
an imbalance in $M$ iff $T$ corresponds to a subset of edges which
is a cycle in $G_D$ and the complement of a cycle in $G_P$. ($D$
and $P$ are the dual and the phial of $M$, respectively).

\section*{5.C\quad Remarks on Chapter 3}
In Chapter 3 we solve a generalized version of the Gauss code
problem, treating together the cases in which the surface involved
is either the plane or the projective plane, under the restriction
of 2-face-colorability.

As a generalization we can remove the restriction of
2-face-colorability and try to solve the problem for the
projective plane. But as it becomes clear from the text of Chapter
3, this is a completely different problem. More relevant for the
present context would be to find an algorithm to determine what is
the minimum connectivity of a surface where a realization of a
given Gauss code (satisfying the 2-face-colorability restriction)
exists. As we showed in Section 3, this would provide an algorithm
to determine whether or not a graph with maximum valency 3 is
map-rich.

Recently we have asked ourselves how one might algorithmically
decide realizability of a Gauss code in the Klein bottle. It seems
that the problem is transformed into a linear system of equations
in $GF(2)$. However, it looks like the methods which we are using
do not extend for the torus. If so, the transition from the Klein
bottle to the torus is the breaking point for the techniques that
we have available at present.

Another line of generalization of the Gauss code problem is to
permit the same symbol to occur an arbitrary number of times in
the input sequence. Thus, multiple crossings appear in the curve.
This generalization was suggested by Gr\"unbaum [10]. This version
of the problem is related to the embedding of hypergraphs.

\section*{5.D\quad Remarks on Chapter 4}
In Chapter 4 we show that to decide whether a given set of squares
of a map\index{squares of a map} is a minimum imbalance is an
$NP$-complete problem. We provide a solution for this problem for
the class of projective maps. The natural question is: how far can
we go? Is it possible to solve the same problem for Kleinian maps?
More restrictively, which are the classes of maps, beyond the
projective maps, for which the minimax relation described in
Theorem 4.5 holds? An apparently unrelated class where it holds is
the class of antimaps of planar maps with eulerian graphs. This
follows from a result by Seymour [23], strengthening a theorem by
Lovasz [16]. This result implies that a minimum transversal of odd
coboundaries in bipartite graph has the same cardinality as a
maximum disjoint collection of odd coboundaries. For planar graphs
the minimal coboundaries correspond to circuits in the planar
dual. Thus, Seymour's result implies: for planar eulerian graphs,
the minimum transversal of odd circuits is equal in cardinality to
a maximum disjoint collection of odd circuits.

To show that the above minimax relation implies (4.5) for the
class of antimaps of planar maps with eulerian graphs, observe
that it is enough to establish that for $m^t$, where $M^\sim$ is
an orientable map, the concept of $r$-circuit is synonymous to odd
circuit in $G_M=G_{M^\sim}$. To prove that the latter two concepts
are equivalent we show, initially, that for an arbitrary map $M$,
a circuit $X$ in $G_M$ is an $r$-circuit in $M^t$ iff for some
imbalance I (and thus an arbitrary imbalance I), $|X\cap I|$ is
odd, where $I$ is considered a subset of $eG_M$. Let $X'$ be a
circuit in $C_M$ such that $\psi^M_c(X')=X$. Let $\{A,B\}$ be a
balancing partition for $M$. Denote by I the imbalance induced by
$\{A,B\}$ and by $f_I$ the $f_M$-edges in the squares of $I$.
Observe that $eG_M\backslash\delta A=f_I$ and that the parity of
$|X'\cap f_I|$ is equal to the parity of $|X'|$. Also, $|X'\cap
f_I|\equiv|X\cap I|\mod\,2$. By definition of $r$-cycle, it
follows that $X$ is an $r$-circuit iff $|X\cap I|\equiv|X'|$ is
odd.

Now we conclude that for the antimap of an orientable map,
$r$-circuit is a synonym for odd circuit. If $M^\sim$ is
orientable, then $\phi$ is an imbalance in $M^\sim$ and, by
(3.2.1), $eG_M$ is an imbalance in $M$. Put $eG_M$ in the role of
$I$ in the above equivalence. The proof is complete.

Two classes of maps for which 4.5 holds are, therefore, the class
of projective maps $M$ with $G_M$ eulerian, and the class of maps
$M$, where $G_M$ is eulerian and $M^\sim$ is planar. These classes
are, apparently, unrelated. It would be interesting to find a
common generalization of these classes for which (4.5) holds.

\bigskip

\printindex

\end{document}